\documentclass[publmath,final,numbook]{svjour}
\usepackage{amsmath,amsfonts}
\usepackage{graphicx}

\usepackage{hyperref}

\reversetheomheadings 
\subsectiondisplayed

\renewcommand{\geq}{\geqslant}
\renewcommand{\leq}{\leqslant}
\newcommand{\thetaval}{\frac{3}{26}}
\newcommand{\bash}{\backslash}
\newcommand{\dti}{d^{\times}}

\newcommand{\ilim}{\mathop{\varprojlim}\limits}

 \newcommand{\whitintertwinerconstant}{\frac{2 \xi_F(2) (\disc F)^{1/2}}{\xi_F^*(1) }}

\DeclareMathOperator{\GL}{GL}
\DeclareMathOperator{\SL}{SL}

\DeclareMathOperator{\modu}{mod}

\def\stacksum#1#2{{\stackrel{{\scriptstyle #1}}{{\scriptstyle #2}}}}
\def\peter#1{\langle #1\rangle}
\def\ov#1{\overline{#1}}

\newcommand{\muP}{\mu_{\mathrm{P}}}
\newcommand{\Eisreg}{\Eis^*}

\newcommand{\defparam}{z}
\newcommand{\Knorm}{\mathscr{U}}
\newcommand{\scrP}{\mathscr{P}}
\newcommand{\geod}{\mathscr{G}}

\newcommand{\mcG}{\mathcal{G}}
\newcommand{\mcH}{\mathcal{H}}
\newcommand{\knconstant}{q^{-d_{\psi}/2} \frac{\zeta_k(1)}{\zeta_k(2)}}

\newcommand{\knconst}{\iota_{k}}
\newcommand{\whitconst}{\eta_k}

\newcommand{\scE}{\mathscr{E}}

\newcommand{\mfg}{\mathfrak{g}}

\newcommand{\Sob}{\mathcal{S}}
\newcommand{\mcN}{\mathcal{N}}

\newcommand{\Id}{\mathrm{Id}}

\newcommand{\shift}{T}

\newcommand{\Sobcond}{C_{\mathrm{Sob}}}

\newcommand{\mcB}{\mathcal{B}}
\newcommand{\mcI}{\mathcal{I}}
\newcommand{\mcK}{\mathcal{K}}
\newcommand{\mcW}{\mathcal{W}}

\newcommand{\bfX}{\mathbf{X}}
\newcommand{\bfE}{\mathbf{E}}

\newcommand{\horrid}{\Upsilon}

\newcommand{\bfB}{\mathbf{B}}
\newcommand{\bfG}{\mathbf{G}}
\newcommand{\bfM}{\mathbf{M}}

\newcommand{\data}{\mathcal{X}}

\newcommand{\paramset}{\mathcal{N}}
\newcommand{\refs}{\eqref}
\newcommand{\Space}{\mathcal{V}}
\newcommand{\ra}{\rightarrow}

\newcommand{\intc}{\frac{1}{2\pi i}\mathop{\int}\limits}

\newcommand{\Cc}{{\mathbb{C}}}

\newcommand{\Zz}{{\mathbb{Z}}}
\newcommand{\linftyexp}{d_0}
\newcommand{\Proj}{\Pi}
\newcommand{\Rr}{{\mathbb{R}}}
\newcommand{\Qq}{{\mathbb{Q}}}
\newcommand{\Aa}{{\mathbb{A}}}

\newcommand{\Ht}{\mathrm{ht}}

\newcommand{\rmB}{{\mathrm{B}}}
\newcommand{\rmN}{{\mathrm{N}}}
\newcommand{\rmA}{{\mathrm{A}}}
\newcommand{\rmZ}{{\mathrm{Z}}}

\newcommand{\supp}{{\mathrm{supp}}}

\newcommand{\mcO}{\mathcal{O}}
\newcommand{\mfo}{\mathfrak{o}}

\newcommand{\eps}{\varepsilon}

\newcommand{\automorphicdual}{\hat{G}_{Aut}}

\numberwithin{equation}{section}
\newcommand{\tr}{\mathrm{tr}}

\renewcommand{\sc}{\mathrm{sc}}

\newcommand{\Eis}{\mathrm{Eis}}

\newcommand{\adconstant}{\kappa}
\newcommand{\adconstantX}{\kappa_2}

\newcommand{\Ad}{\mathrm{Ad}}

\newcommand{\Whit}{\mathcal{W}}
\newcommand{\disc}{\mathrm{disc}}

\newcommand{\C}{\mathbb{C}}

\newcommand{\vol}{\mathrm{vol}}

\newcommand{\G}{\mathbf{G}}

\newcommand{\Cond}{\mathrm{Cond}}

\newcommand{\height}{\mathrm{ht}}

\newcommand{\PGL}{\mathrm{PGL}}

\newcommand{\kirill}{\mathscr{K}}
\newcommand{\order}{\mathfrak{o}}

\newcommand{\R}{\mathbb{R}}
\renewcommand{\H}{\mathbb{H}}

\newcommand{\Z}{\mathbb{Z}}
\newcommand{\Q}{\mathbb{Q}}

\newcommand{\adele}{\mathbb{A}}

\newcommand{\quot}{\mathbf{X}}
 
\newcommand{\mcF}{\mathcal{F}}

\DeclareFontFamily{OT1}{rsfs}{}
\DeclareFontShape{OT1}{rsfs}{n}{it}{<-> rsfs10}{}
\DeclareMathAlphabet{\mathscr}{OT1}{rsfs}{n}{it}

\newcommand{\beq}{\begin{displaymath}}
\newcommand{\eeq}{\end{displaymath}}

\newcommand{\be}{\begin{equation}}
\newcommand{\ee}{\end{equation}}

\spnewtheorem*{Rem}{Remark}{\it}{\rm}
\spnewtheorem{princ}[subsection]{Ergodic principle}{\it}{\it}
\spnewtheorem*{hypo}{Hypothesis $H_{\theta}$}{\it}{\it}
\spnewtheorem*{Lemma*}{Lemma}{\it}{\it}
\spnewtheorem{Lemmat}[subsubsection]{Lemma}{\it}{\it}
\spnewtheorem{Prop}[subsubsection]{Proposition}{\it}{\it}
\spnewtheorem{Defi}[subsubsection]{Definition}{\it}{\it}
\spnewtheorem{Exam}[subsubsection]{Example}{\it}{\it}
\spnewtheorem*{Theorem*}{Theorem}{\it}{\it}
\spnewtheorem*{Corollary*}{Corollary}{\it}{\it}
\spnewtheorem*{Question*}{Question}{\it}{\it}
\spnewtheorem*{Proposition*}{Proposition}{\it}{\it}

\usepackage{publmath}  
\include{diagram}
\newarrow{Int}{boldhook}{-}{-}{-}{>}

\begin{document}

\title{THE SUBCONVEXITY PROBLEM FOR $\GL_{2}$.}
\author{\firstname{Philippe} MICHEL\thanks{Partially supported by the advanced research grant n. 228304 from the European Research Council and the  SNF grant 200021-12529.}\and \firstname{Akshay} VENKATESH\thanks{Partially supported by the Sloan foundation, the Packard Foundation and by an NSF grant.}}
\institute{Ph. M.\\ Ecole Polytechnique F\'ed\'erale de Lausanne,\\philippe.michel@epfl.ch \and A. V.\\ Stanford University,\\ akshay@math.stanford.edu}

\setcounter{tocdepth}{2}

\maketitle
   \begin{abstract} Generalizing and unifying prior results, we solve the subconvexity problem for the $L$-functions of $\GL_{1}$ and $\GL_{2}$ automorphic representations over a fixed number field, uniformly  in all aspects.  
A novel feature of the present method is the softness of our arguments; this is largely due to a consistent use of  canonically normalized period relations, such as those supplied by the work of Waldspurger
and Ichino--Ikeda.
\end{abstract}

\tableofcontents

\section{Introduction}

\subsection{The subconvexity problem}
We refer the reader who is not familiar with $L$-functions to \S \ref{classical}
for an introduction, in explicit terms, to some of the ideas of this paper.

Throughout this paper, $F$ denotes a fixed number field  and $\Aa$
its ring of ad\`eles. For $\pi$ an automorphic
representation of $\GL_n(\Aa)$ (with unitary central character, not necessarily of finite order),
Iwaniec and Sarnak have attached an {\em analytic conductor} $C(\pi) \in \R_{\geq 1}$; 
it is the product of the usual (integer) conductor with a parameter measuring
how large the archimedean eigenvalues are. More intrinsically,
the logarithm of the conductor is proportional to the density of zeros
of the corresponding $L$-function.   See \S \ref{notation-global}. 

The {\em subconvexity problem} is concerned with the {\em size} of $L(\pi,s)$ when $\Re s=1/2$: it consists in improving over the so-called
{\em convexity bound} (see \cite{GAFA2000} for instance)
$$L(\pi,1/2)\ll_{n,F,\eps}C(\pi)^{1/4+\eps}$$
for any\footnote{Recently, Heath-Brown has established a general convexity bound \cite{HBconvex}, which together with the work of Luo, Rudnick, Sarnak \cite{LRS2} implies the clean convexity bound $L(\pi,1/2)\ll_{n,F}C(\pi)^{1/4}$.} $\eps>0$.  
The main result of the present paper is the resolution of this problem
for $\GL_{1}$ and $\GL_{2}$-automorphic representations:

\begin{Theorem}\label{GL12}  There is an absolute constant $\delta>0$ such that: for $\pi$ an automorphic representation of $\GL_1(\Aa)$ 
or $\GL_2(\Aa)$ (with unitary central character), one has
$$L(\pi,1/2)\ll_F C(\pi)^{1/4-\delta}.$$
\end{Theorem}
\begin{Rem} Contrary to appearance, this also includes the question of growth along
the critical line, i.e. what is called the $t$-aspect, because $$L(\pi,1/2+it)=L(\pi\otimes|.|_{\Aa}^{it},1/2).$$
For example, an interesting corollary is a subconvex bound for the $L$-function
of a Maass form with eigenvalue $1/4 +\nu^2$ at the point $t=1/2+i \nu$. Another corollary is a subconvex bound (in the discriminant) for the
 Dedekind $L$-function of a cubic extension of $F$ (cf. \cite{ELMV3} for an application of the latter).
\end{Rem}

The above result is a specialization (by taking $\pi_{2}$ to be a suitable Eisenstein series) of the following more general result

\begin{Theorem}\label{RSthm}  There is an absolute\footnote{independent of the number field $F$.} constant $\delta>0$ such that: 
for $\pi_1, \pi_2$ automorphic representations on $\GL_2(\adele_F)$ we have\footnote{More precisely we prove the bound $L(\pi_{1}\otimes\pi_{2},{1}/{2} )\ll_{F,\pi_2}
  C(\pi_1)^{1/2-2\delta}$. That the latter implies the former is a consequence of the bounds in \cite{BHe}.}
:
\begin{equation}\label{eq:RSthm}
L(\pi_{1}\otimes\pi_{2},{1}/{2} )\ll_{F,\pi_2}
  C(\pi_1\otimes\pi_{2})^{1/4-\delta};
\end{equation}
more precisely, the constant implied depends polynomially
  on the discriminant of $F$ (for $F$ varying over fields of given degree) and on $C(\pi_{2})$.
 \end{Theorem}

The value of $\delta$ is easily computable. We have not attempted to optimize any exponent, our goal in this paper being of  giving clean proofs in a general context. 

\begin{Remark} The bound generalizes (up to the value of $\delta$) a variety of subconvex bounds \cite{Weyl,Burgess,HB,Good,Meurman,DFI1,DFI2,Sar2,FoIw,Ivic,DFI8,KMV,LiuYe,M04,Bl,MotJut,HM,BlHM,BH,Ve,DG}.  Its main feature, however, is its uniformity in the various possible parameters (the so-called ``conductor'', $t$-aspect, or  ``spectral'' aspects): such bounds are sometimes called ``hybrid''.  The first such hybrid bound is that of Heath-Brown \cite{HB} for Dirichlet character $L$-functions;  recent hybrid bounds are to be found in the work of Jutila and Motohashi  \cite{MotJut} and in the work of Blomer and Harcos \cite{BH}. 
The present work generalizes these also, again, up to the value of $\delta$. 

However, 
in saying this we have done an injustice to some of those papers; 
in some cases, the very point was to obtain the best exponent,
whereas our emphasis is quite different. For example, 
let us compare the present result to that of \cite{MotJut}. This paper gives, in particular, the uniform  bound for the value
$L(\varphi_{\lambda},\frac{1}{2}+it) \ll (|t|+|\lambda|)^{1/3}$ where $\varphi_{\lambda}$ is a Maass form.
This bound is very strong -- the same exponent as the classical Weyl bound; on the other hand, 
it fails to be subconvex when $t \sim \lambda$, where the conductor drops.
The present work fills this lacuna and provides a subconvex bound for that ``critical point''; as far as we are aware this subconvex bound is new
even over $\Qq$.
 On the other hand, while our method presumably 
leads to a respectable $\delta$, it would not be so strong as the result of \cite{MotJut}, one reason being that we are using the {\em amplification method}. 
 It is also worth observing that the phenomena of the conductor dropping
often leads to major difficulties, both in our methods and other treatments. 

    \end{Remark}
    
    \begin{Remark}
 It is reasonable to ask what one hopes by studying hybrid settings -- especially 
given that many applications of subconvexity do not require them. 

  It is generally believed that the analytic
 behavior of $L$-functions are ``universal,'' in that they are controlled by a single scaling parameter,
 the analytic conductor $C$. Taking $C \rightarrow \infty$
 in different ways can correspond to analysis of eigenfunctions of large eigenvalue;
 analysis of eigenfunctions on a surface of large volume; or sections of a highly ample holomorphic bundle, and the reasons why these should all have similar asymptotic behavior is not clear. 

We may hope to achieve some insight into this ``universality'' by studying hybrid phenomena.
Indeed, at many points in the text, the reader will note the close parallel between increasing
the level at finite primes, and increasing the eigenvalue at archimedean primes.  
\end{Remark}

Our methods were outlined in our ICM announcement \cite{MV-ICM}.
 In particular we do not use trace formulas of any kind.  
 An additional feature (origi\-nating in \cite{M04}), is that a {\em special case} of Theorem \ref{RSthm} (the case of $\pi_{1}$ a character) enters the proof of the full theorem. That special case is proven in Theorem \ref{thmcharacter}, and is based (following \cite{Ve}) on a study of the equidistribution  of cycles on adelic quotients.

We have also tried to make use of the following ideas to simplify the proof (both notationally and conceptually):
\begin{itemize}
\item[-] Sobolev norms (cf. \cite{BR}; in the adelic context \cite{Ve}); 
\item[-] Canonically normalized period formulas (see \cite{Wald,IchIk});
\item[-] Regularization of integrals of automorphic forms (we give a self-contained treatment that avoids truncation).  
  \end{itemize} 
  
 There remain many interesting questions related to the subconvexity story even for $\GL_2$. 
For instance, the ``approximate functional equation'' gives a way
to numerically compute any given $L$-function at the central point, in time $C^{1/2+\eps}$,
where $C$ is the analytic conductor. An interesting question is whether
some of the ideas that enter into the proof of subconvexity can be interpreted
to give faster algorithms. There is some suggestion of this in existing fast algorithms
for computation of $\zeta(1/2+it)$.   

\subsubsection{An outline of the proof} 
To conclude this section we outline the proof of the main theorem \ref{RSthm};
a more elementary discussion is in \S \ref{classical} and \S \ref{classical2}.  Consider two (generic) automorphic representations $\pi_{1}$, $\pi_{2}$; for simplicity we assume that  both are cuspidal\footnote{However, we devote some time and effort to handling
the necessary regularizations in the general case; to our surprise, these modifications are not ugly but rather beautiful.}  and that $\pi_{2}$ is fixed.
We aim for a subconvex bound for the Rankin/Selberg $L$-function central value of the form
$$\frac{ L(\pi_{1}\otimes\pi_{2},1/2)}{C(\pi_{1}\otimes\pi_{2})^{1/4}}\ll_{\pi_{2}}C(\pi_{1})^{-\delta},\ \delta>0.$$
By Rankin/Selberg theory we may realize the left-hand side above as a triple product period $\int_{\PGL_{2}(\Qq)\bash\PGL_{2}(\Aa)}\varphi_{1}\varphi_{2}E(g)dg=\peter{\varphi_{1},\ov{\varphi_{2}E}}$ for suitable  $\varphi_{i} \in \pi_i \ (i=1,2)$ and
$E$ belonging to the Eisenstein series of type ``$1 \boxplus \chi$,'' where $\chi$ is the inverse
of the product of the central characters of $\pi_1$ and $ \pi_2$. 
Now $|\peter{\varphi_{1},\ov{\varphi_{2}E}}|^2$
is bounded by
\be\label{commute}\peter{\varphi_{2}E,\varphi_{2}E}=\peter{\varphi_{2}\ov{\varphi_{2}},E\ov E},\ee
and it suffices to show this is at most $C(\pi_1)^{-\delta}$. 
Now, $E\ov E$ is not square
integrable and one needs a regularized version of the above inner product;
 this is described in \S \ref{regularsection} and \S \ref{triple}.   
This being done, we obtain by spectral expansion (see\S \ref{regularIP}) 
\be\label{introexpand}
\peter{\varphi_{2}\ov{\varphi_{2}},E\ov E}_{reg}=
\peter{\varphi_{2}\ov{\varphi_{2}},\varPhi}_{reg}+\int_{\pi }\sum_{\mcB(\pi)}\peter{\varphi_{2}\ov{\varphi_{2}},\varphi}_{reg}
\peter{\varphi,E\ov E}_{reg}d\mu_{P}(\pi)
\ee
where the subscript $reg$ denotes regularized inner products (one has of course $\peter{\varphi_{2}\ov{\varphi_{2}},E\ov E}=\peter{\varphi_{2}\ov{\varphi_{2}},E\ov E}_{reg}$) and $\varPhi$ is a certain non-unitary Eisenstein series, which is an artifact of the integral regularization; 
on the other hand, $\pi$ varies over the automorphic dual of $\PGL_{2}$, $d\mu_P$ is a ``Plancherel measure,'' and  $\varphi$
 varies over an orthonormal basis, $\mcB(\pi)$, of factorizable vectors in $\pi$.  Now:
   \begin{enumerate} 
 \item  The term $\peter{\varphi_{2}\ov{\varphi_{2}},\varPhi}_{reg}$ is handled via the amplification method (\S \ref{degeneratebound}).
\item  
The terms $\peter{\varphi,E\ov E}_{reg}$ are bounded by some negative power of $C(\pi_{1})$  (\S \ref{genericbound}).  Rankin/Selberg theory implies that $\peter{\varphi,E\ov E}_{reg}$ factor into a product of local integrals to which bounds for matrix coefficient can be applied (see below) {\em times} the central value $L(\pi,1/2) L(\pi \times \chi, 1/2)$; we eventually need only bounds for
$L(\pi \times \chi, 1/2)$, where $\pi$ is essentially fixed\footnote{
 As it turns out, the vector $\varphi_{2}\in\pi_{2}$ depends only on $\pi_{2}$ (up to archimedean components, this is just the new vector), as the latter is fixed, the quantities $\peter{\varphi_{2}\ov{\varphi_{2}},\varphi}_{reg}$ decay very rapidly as the eigenvalue or level of $\varphi$ increases. So, we may regard $\varphi$ as essentially fixed.}
and $\chi$ varying.

\end{enumerate}

The phenomenon of reduction to another case of subconvexity was noted by the first-named author in \cite{M04}.  We establish this case in Theorem \ref{thmcharacter} (again, in all aspects);
the proof generalizes \cite{Ve}, and
 we refer to the introduction of \cite{Ve}
and to \S \ref{classical} for intuition about it.

 \subsubsection{Local computations}
Let us be more precise about the local integrals that occur in, e.g., the Rankin-Selberg method. 
For our purpose they must be examined carefully; we are particularly interested in their {\em analytic}
 properties, i.e., how large or small they can be. This is the purpose of Part III of the paper; we deal also with the local Hecke integrals.  
  
For example, the local integral occuring in the Rankin-Selberg method can be interpreted
 as a linear functional $\ell$ on the tensor product $\pi_1 \otimes \pi_2 \otimes \pi_3$ 
 of three representations $\pi_j$ of $\GL_2(k)$. 
The space of such functionals is at most one dimensional; $|\ell|^2$ is therefore proportional
to the Hermitian form
 $$x_{1}\otimes x_{2}\otimes x_{3}\mapsto \int_{\PGL_{2}}\peter{x_{1},g.x_{1}}\peter{x_{2},g.x_{2}}\peter{x_{3},g.x_{3}}dg,$$
and this is of tremendous utility for the analytic theory.  Similarly, when studying the Hecke integral, it is most convenient to study 
 the Hermitian form defined by $x \mapsto  \int_{A(k)}\peter{a.x,x}d a$ where $A$ the group of diagonal matrices (modulo the center).  
 
 It is a wonderful observation of Waldspurger (see also \cite{IchIk}) that global period formulas
 become very simple when expressed using such canonically normalized local functionals.

\subsubsection{Spectral identities of $L$-functions} The  identity \eqref{commute}, although evident, is the keystone of our argument. Its usage in the ``period'' form presented above seems to have been noticed independently by the present authors in their attempt
to geometrize \cite{M04},  and by Bernstein and Reznikov in their work on the subconvexity problem for triple product $L$-functions \cite{BRlater}.

However \eqref{commute} also manifests itself at the level of $L$-functions, and looks
rather striking in this guise: indeed \refs{introexpand} may be recognized as an identity between (weighted) sums of central values of triple product $L$-functions
and (weighted) sums of ``canonical'' square roots of similar central values:
\begin{multline*} \int_{\pi_{1}}w(\pi_{1})L(\pi_{1}\otimes\pi_{2}\otimes\pi_{3},1/2)d\mu_{P}(\pi_{1})\\ =\int_{\pi}\tilde w(\pi)\sqrt{L(\pi_{2}\otimes\tilde\pi_{2}\otimes\pi,1/2)}.\sqrt{L(\pi_{3}\otimes\tilde\pi_{3}\otimes\pi,1/2)}d\mu_{P}(\pi)\end{multline*}
In this form, this identity was discovered already by N. Kuznetsov \cite{Kuz} (with $\pi_{2}$, $\pi_{3}$ Eisenstein series see also \cite{KuzMot}) and an interesting application of this result was made by M. Jutila \cite{Jut}.  This period identity has -- implicitly or explicitly --  played an important role in the analytic theory of $\GL_2$ forms.
  We refer to \S \ref{fulllevel} for some more discussion of how to convert
\eqref{commute} to an identity of $L$-functions.

In fact, it was shown by A. Reznikov (in the archimedean setting at least), that such phenomena is not isolated, and may be systematically described through the formalism of {\em strong Gelfand configurations}: that is commutative diagrams of algebraic groups
$$
\begin{diagram}[thick]
&\mcG&\\
\mcH_{1}\ruInt(1,1)&&\luInt(1,1)\mcH_{2}\\
&\luInt(1,1)\mcF\ruInt(1,1)&\\
\end{diagram}
$$

in which the pairs $(\mcF,\mcH_{i})$ and $(\mcH_{i},\mcG)$ are strong Gelfand pairs. The present paper
corresponds to the configuration
$\mcG=\GL_{2}\times\GL_{2}\times\GL_{2}\times\GL_{2}$, $\mcH_{i}=\Delta\GL_{2}\times\Delta\GL_{2}$ with two different diagonal embeddings and $\mcF=\Delta\GL_{2}$ (diagonally embedded). As is explained in \cite{AndreJAMS}, such  strong Gelfand configurations yield naturally to spectral identities between periods and then between $L$-functions.
We refer to {\em loc. cit.} for more interesting examples of that sort.

\subsubsection{Structure of the paper; reading suggestions}
The paper splits into five parts; the first four parts are largely independent of each other,
and contain various results of independent interest. The fifth part brings together these results
to prove the main theorem.

The reader may wish to skip directly to Parts IV and V of the paper; the results of Parts II and III
are largely technical and not particularly surprising. 

\begin{itemize}
\item[-] In the remainder of Part I,   \S \ref{classical}, \S \ref{classical2} we
consider two corollaries to our main theorem -- which can both be phrased {\em without} $L$-functions -- and we explain how the proofs work in these instances. Indeed, the general proof
is obtained by adelizing and combining these two particular cases.  The corollaries we consider are:

\begin{enumerate}
\item In \S \ref{classical}, we discuss the ``Burgess bound,'' which relates
to the issue of the smallest quadratic non-residue modulo a prime $q$.
\item In \S \ref{classical2}, we discuss a problem in analysis on a negatively curved surface,
viz.: how large can the Fourier coefficients of an eigenfunction along a closed geodesic be ?
\end{enumerate}

\item[-] Part II (viz.  \S \ref{plancherel} to \S \ref{sobolevnormsproofs}) is of more general nature:
 we discuss a system of Sobolev norms on adelic quotients, inspired largely by work of Bernstein and Reznikov.  This section exists to give a suitable language for talking about adelic equidistribution, and the norms
are a somewhat cleaner version of those appearing in \cite{Ve}. Some of the remarks
here are of independent interest, although they are of technical nature. 

\item[-] Part III discusses some of the analytic theory of torus-invariant functionals
on a representation of $\GL_2(k)$ (where $k$ is a local field), and of trilinear functionals
on representations of $\GL_2(k)$. 

\item[-] Part IV discusses the global theory of torus periods on $\GL_2$ and
the diagonal period on $\GL_2 \times \GL_2 \times \GL_2$. 

\item [-]Part V gives the proof of Theorem \ref{RSthm} along with the important intermediary result
Theorem \ref{thmcharacter} (a subconvex bound of $L$-function of character twist uniform in the character aspect).

\end{itemize}

\begin{acknowledgement} The present work started during a visit of the first author at the Courant Institute (New York) and ended basically during the workshop ``Analytic Theory of GL(3) Automorphic Forms and Applications'' at the AIM (Palo Alto); parts of it were written during visits at the RIMS, at the IHES and at Caltech for the first author and visits at the IAS and the IHES for the second.
 Both authors would like to thank these institutions for their support and the very nice accommodations provided.
 
 We would like to express our thanks to J. Bernstein, E. Lapid  for many helpful discussions
 related to regularization, and to Y. Sakellaridis for discussions related to normalization
 of functionals. In addition, we have learned
 greatly from the work of J. Bernstein and A. Reznikov. We thank them both
 for sharing their knowledge generously; the reader will find many techniques
 inspired by their work in this paper.  We would also like to aknow\-ledge the influence of the series of papers of Duke, Friedlander, Iwaniec; besides the use of the, by now standard, amplification method, the reader will recognize, in a disguised form, many of the techniques invented in their work.
 Finally, we would like to thank D. Ramakrishnan and W. T. Gan for useful discussions related to this work, F. Brumley, E. Fouvry and N. Templier for their comments and corrections and we wish to address  special thanks to Peter Sarnak for his continuous encouragements and advices. 
    
This project has been stretched over several years and we would like to apologize for the long delay between the announcement (e.g. in  \cite{MV-ICM}) of some of the results of this paper and its appearance. By way of excuse,
 the results here are now substantially stronger. 
Perhaps more importantly, the present proofs are far more streamlined than our original proofs; although
the general outline is as in \cite{MV-ICM}, the use of regularization and canonical inner product
formula made many of the details far more manageable. 
\end{acknowledgement}
\subsection{The Burgess bound and the geometry of $\mathrm{SL}_2(\mathbb{Z})\backslash \mathbb{H}$}
 \label{classical}

In this section and also in \S \ref{classical2} we present some of the ideas of the general proof
in the most down-to-earth setting as we could manage.   In both these sections, we have
by and large eschewed mention of $L$-functions.

\subsubsection{The Burgess bound}
Let $\chi$ be a Dirichlet character to the modulus $q$.  It is well-known that subconvexity
for the Dirichlet $L$-function $L( \chi,1/2)$, in the $q$-aspect, is substantively equivalent
to a bound of the nature:
\begin{equation} \label{Burgess-subconvex} \left| \sum_{i=1}^{M} \chi(m)  \right| \leq M q^{-\delta},\end{equation} 
with $M \approx q^{1/2}$ and some absolute constant $\delta>0$. 

 In this question, $q^{1/2}$ is a ``threshold'': indeed it is rather easy to establish \eqref{Burgess-subconvex} when $M = q^{0.51}$ (the Polya-Vinogradov inequality). The bound \refs{Burgess-subconvex}, for $M$ in the range $M \approx q^{1/2}$ was proven (in a wider range) by Burgess \cite{Burgess}. 

In the present section (see also \cite{Ve-Pisa}) -- which, we hope, will make sense to the reader without any knowledge of $L$-functions -- we explain how \eqref{Burgess-subconvex} is related to an equidistribution equation
on the space of lattices,  and then discuss how to prove the uniform distribution statement. 
A key part of the paper -- \S \ref{proof} -- will implement the discussion of this section in a more general context. 

{\em To simplify that discussion, we assume for the rest of this section  that $q$ is prime and that $\chi$ is the Legendre symbol.}

\subsubsection{The space of lattices}
Put $X = \SL_2(\Z) \backslash \SL_2(\R)$, the space of unimodular lattices in $\R^2$. 

We say a sequence of finite subsets $S_i \subset X$ is becoming uniformly distributed if, for any $f \in C_c(X)$, we have $$\frac{1}{|S_i|} \sum_{S_i} f \rightarrow \int_{X} f,$$ the latter integral being
taken with respect to the unique $\SL_2(\R)$-invariant probability measure on $X$.

\subsubsection{Burgess bound and lattices} 

For $x \in \R$ consider the lattice
$$\Lambda_x = \frac{1}{\sqrt{q}}\bigl(\Zz.(1,x)+\Zz. (0,q) \bigr) \in X.$$ 

As $x$ varies, $\Lambda_x$ moves on a horocycle in $X$ -- an orbit
of the group of upper triangular, unipotent matrices.   This horocycle is in fact {\em closed},
since $\Lambda_{x+q} = \Lambda_{x}$. 

Given $0<\eta<1$, let $F: X \rightarrow \mathbb{R}$ be defined by $F(L) =  | L\cap [0, \eta]^2|-1$, i.e. $F$ counts the number of non-trivial lattice points in a small square box. 
A simple computation shows that for $x\in\Zz$, $F(\Lambda_x)$ equals the number of nonzero solutions
$(\alpha, \beta) \in [0, \eta \sqrt{q}]^2 \cap \mathbb{Z}^2$ to the equation
$$\beta  \equiv \alpha x \hbox{ mod } q.$$
It follows that
$$
\frac{1}{q}\sum_{x \modu q} F(\Lambda_x) \chi(x) = \frac{1}{q} \big|\sum_{\beta \in [1, \eta \sqrt{q}]} \chi(\beta)\big |^2$$ 

Now, the Burgess bound would follow if we knew
that the right-hand side was small. The above identity therefore relates the Burgess bound to a type of equidistribution statement:
 we must show that the sets \begin{equation} \label{ud-1}\{\Lambda_x: x \in \Z/q\Z \mbox { is a quadratic residue}\},\end{equation}
and the similar set for quadratic nonresidues, are uniformly distributed\footnote{Strictly speaking, the function $F$ is not of compact support on $X$; in fact, it grows at the cusps. We shall ignore this technical detail for the purpose of explanation.} on $X$. 

\begin{Rem}
This connection between the bound \eqref{Burgess-subconvex}, and a uniform distribution statement on the space of lattices, is not an accident: it is a special case of the connection between $L$-functions and automorphic forms.  Indeed, the uniform distribution \eqref{ud-1} encodes much more than \eqref{Burgess-subconvex}: it encodes, at once, subconvex bounds for twists $L(\frac{1}{2}, f \times \chi)$ where $f$ is a fixed $\SL_2(\Zz)$-modular form (of any weight); the latter specializes to the former when $f$ is an Eisenstein series. Our point above, however, is that the connection between \eqref{Burgess-subconvex}
and the space of lattices can be made in an elementary way. 
\end{Rem}

It is possible to visualize the desired uniform distribution statement by projecting
from the space of lattices to $Y := \mathrm{SL}_2(\Z)\bash\mathbb{H}$ (at the price of losing the group actions). 
The lattice $\Lambda_x$ projects to the class of $z_x = \frac{i}{q} + x$.

\begin{figure}
\includegraphics[scale=0.5]{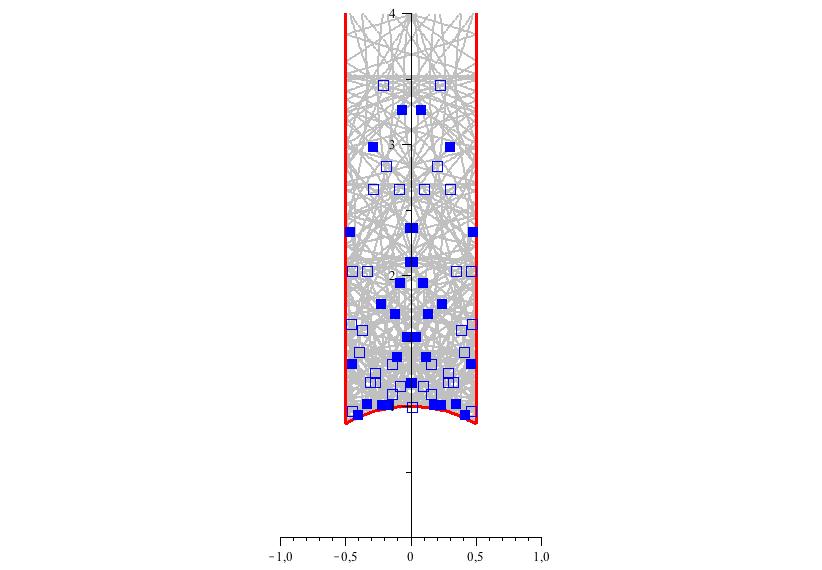}
\caption{
The horocycle $\frac{x}{173}+\frac{i}{173}$. Gray dots: $x\in\Rr$, plain (resp. empty) square  $x \in \Z$ a quadratic (resp. non-quadratic) residue mod $173$.
}
\end{figure}

\subsubsection{Equidistribution statements}\label{classicaled}
We shall try to establish \eqref{ud-1} by first proving uniform distribution of a ``bigger'' set,
and then refining that statement. 
 Consider, then, the following three equidistribution statements, 
as $q \rightarrow \infty$:
\begin{enumerate}
\item[(A)]  The closed horocycle $\{\Lambda_x: x \in [0,q]\}$ becomes uniformly distributed on $X$; 
\item[(B)]  $\{ \Lambda_x: x \in \Z \cap [0,q] \}$ becomes uniformly distributed on 
$X$; 
\item[(C)] $\{\Lambda_x:  x \in \Z \cap [0,q] \}$ becomes u.d. on $X$, {\em when each point $\Lambda_x$ is weighted by $\chi(x)$.}
\end{enumerate}

For $\chi$ the quadratic character modulo $q$, we might rewrite (C) as:
\begin{enumerate}
\item[(C2)]$\{\Lambda_x:  x \in \Z \cap [0,q] , x \mbox{ a quadratic residue mod $q$}\}$ becomes u.d. on $X$. 
\end{enumerate} 

We have already discussed informally, and it is true in a very precise sense,
that (C) and (C2) have substantively the same
content as the subconvexity result that we are aiming for.  Note that (C) and (C2) seem ``unnatural'' at first; it seems like $z_x$ is parameterized by an {\em additive} structure, i.e. $\Zz/q\Z$;
thus it is odd to restrict attention to a ``multiplicatively'' defined set. 
But in truth the examples of (C) 
and (C2) have -- as we shall see --  an underlying ``multiplicative'' symmetry; the fact that it appears additive
is a reflection of the degeneration of a torus in $\GL_2$ to a unipotent group.

Observe that (A), (B), (C2) are asserting the equidistribution of smaller and smaller sets.
So what we need, besides a proof of (A) -- which happens to be an old result of Peter Sarnak --  is a method to pass from the equidistribution of a large set,
to the equidistribution of a smaller subset. It is provided by the following easy principle (cf. \cite{Ve-Pisa}
for a further discussion of its applicability in this type of context):

\begin{princ}\label{Prinzip} Suppose a group $H$ acts ergodically on a probability space $(X,\mu)$, and $\nu$
is an $H$-invariant measure which is dominated by some  positive multiple of $\mu$.
Then $\nu$ is a scalar multiple of $\mu$. 
\end{princ}
By ``$\nu$ is dominated by some  positive multiple of $\mu$'' (written $\nu \ll \mu$) we mean that there is a constant $A$
such that for any measurable set $S\subset X$, $|\nu|(S)\leq A\mu(S)$. Indeed, 
if this is so, $\nu$ is absolutely continuous with respect to $\mu$, and thus may be expressed as $f \mu$ for some $f \in L^1(\mu)$; then $f$ must be a $H$-invariant function, necessarily constant by ergodicity. 
The principle is therefore trivial; on the other hand, its consequences in the number-theoretic
context are surprising.

Using the principle, we can pass from the equidistribution of $\mu$ to the equidistribution of $\nu$. 
(Of course, one needs a more quantitative form of this principle; we enunciate such a form
in the context we need in \S \ref{Prinzip-quant}.) 

By applying the ergodic principle to the group $\{n(t): t\in \mathbb{Z}\}$, we establish
the implication $(A) \implies (B)$. To show that $(B) \implies (C)$ is a little more subtle,
because the choice of $H$ is not clear; we will need to pass to a covering to uncover it!
 We now discuss this in more detail.

\subsubsection{The entry of the ad{\`e}le group and the implication $(B) \implies (C)$} \label{Adele}

It is well-known that the inverse limit (over the principal congruence subgroups)
$$\widetilde{X} := \ilim_{q} \Gamma(q) \backslash \SL_2(\R)$$
carries not only an action of $\SL_2(\R)$ but of the much larger group $\SL_2(\Aa_{\Q})$,
where $\Aa_{\Q}$ is the adele ring of $\Q$; there is a natural projection $\pi: \widetilde{X} \rightarrow \SL_2(\Z) \backslash \SL_{2}(\Rr)$. 

Let $\nu$ be the measure implicit in (C), that is to say,
$\nu = \sum_{x=1}^{q} \left( \frac{x}{q} \right) \delta_{\Lambda_x}$. 
As we have commented, $\nu$ has no apparent invariance.  

{\em However, there is a closed subgroup $H^{(1)} \subset \SL_2(\adele_{\Q})$, and
a $H^{(1)}$-equiva\-riant measure $\tilde{\nu}$ on $\tilde{X}$, which projects to $\nu$. }
In other terms, the measure of (C) {\em acquires} invariance after lifting
to the adeles. (The group $H^{(1)}$ admits a surjection onto $(\Z/q\Z)^{\times}$, and this surjection
is compatible with the natural action of   $(\Z/q\Z)^{\times}$ on the sets in (B), (C)). 
Now, by a suitable application of the ergodic principle on $\widetilde{X}$, rather than $X$, 
we deduce that $(B) \implies (C)$. 

In this way, the role of adeles in our proof is not merely to provide
a convenient language, but also the group actions that we use simply do not exist
at the level of $\SL_{2}(\Zz)\bash\SL_{2}(\Rr)$ -- or rather, only their shadows, the Hecke operators, are visible.  

\begin{Rem} 
(The adeles for dynamicists). A somewhat more intuitive way of constructing this is as follows:
fix a prime $p$, and consider the $p$th Hecke operator $T_p$ on $X$. 
It is a multi-valued function. 

We can {\em formally} turn it into an invertible single-valued function  by considering the space of sequences:
$$\tilde{X}_p = (\dots, x_{-2}, x_{-1}, x_0, x_1, \dots ) \in X^{\Z}, x_{i+1} \in T_p x_i \mbox{ for all } i.$$
Then the ``shift'' operation $S: \tilde{X}_p \rightarrow \tilde{X}_p$ can be considered
as a version of $T_p$ that has been forced to be invertible. 

It can be verified that $\tilde{X}_p$ is isomorphic to the quotient:
$$\SL_2(\Z[p^{-1}]) \backslash \SL_2(\R) \times \SL_2(\Q_p) / M,$$
where $M$ is the subgroup of diagonal matrices in $\SL_2(\Q_p)$ whose entries
belong to $\Z_p^{\times}$. Moreover, $S$ is identified with the right action of a diagonal matrix
in $\SL_2(\Q_p)$. 

If we imitate this procedure for {\em all primes simultaneously}, 
one is naturally led to the space $\tilde{X}$. 
\end{Rem}

\subsection{Geodesic restriction problems} \label{classical2}

In this section, we present another corollary to our main results and discuss the idea of its proof,
again, largely without mention of $L$-functions. This section is phrased in the language of analysis on a Riemannian manifold. Our discussion can be considered a variation on the sketch of proof
that was already presented in our ICM article \cite{MV-ICM}. 

\subsubsection{Geodesic restriction problems: the results of \cite{BGT}}

Let $M$ be a Riemannian surface of finite volume, with Laplacian $\Delta_M$, and let
$\varphi_{\lambda}$ be an eigenfunction of $\Delta_M$ with eigenvalue $-\lambda^2$.
Let $\geod$ be a closed geodesic of length $L$ on $M$; we fix a parameterization $t \mapsto \gamma(t)$
of $\geod$ by arc length, so that $\gamma(t+L) = \gamma(t)$.

  In this section, we shall discuss the restriction of $\varphi_{\lambda}$ to $\geod$. A theorem
   of Burq, Gerard, and Tzvetkov \cite{BGT}, generalizing a result of Reznikov \cite{Rez}, asserts the following general bound \
   $$ \frac{ \|\varphi_{\lambda}\|_{L^2(\geod)} } { \|\varphi_{\lambda}\|_{L^2(M)} } \leq C(\geod, M)   \ \lambda^{1/4} ,$$
   the constant $C(\geod, M)$ depending only on $M$ and the geodesic $\gamma$.

  This is in fact the ``worst possible behavior'' as the following basic example shows: let $M= S^2$, embedded as the unit sphere in $\mathbb{R}^3$, with the induced metric. Let $\varphi_n$ be the restriction to $M$ of $(x,y,z) \mapsto (x+iy)^n$. Then $\varphi_n$ is a Laplacian eigenfunction of $L^2$-norm $\asymp n^{-1/4}$; on the other hand, when restricted to the equatorial geodesic $z=0$
  it corresponds to the function:  $$\gamma(t) \mapsto e^{i n t}.$$
  In particular, $\varphi_n| \gamma$ is a {\em single Fourier mode}, and
  $\|\varphi_n\|_{L^2(\geod)} \asymp n^{1/4} \|\varphi_n\|_{L^2(M)}$.

Returning to the general case, let $\omega$ be an integral multiple of $\frac{2 \pi}{L}$. Consider the ``$\omega$-th'' Fourier coefficient along $\geod$, namely,  $$a(\varphi_{\lambda}, \omega) =
\int_0^{L} \varphi_{\lambda}(\gamma(t)) e^{i \omega t}  dt .$$ It measures
the correlation of $\varphi_{\lambda}|\geod$ with a single Fourier mode. By Cauchy-Schwarz:
$$ \frac{ |a(\varphi_{\lambda}, \omega)|}{ \|\varphi_{\lambda}\|_{L^2(M)} } \leq L^{1/2}  \leq  C(\geod, M) \lambda^{1/4}.$$
As the above example shows, this bound is indeed sharp.

\subsubsection{Geodesic restriction problems in the arithmetic case}
 As was explained to us by A. Reznikov, a consequence to our main result, Theorem \ref{RSthm}, is that this behavior {\em never} occurs
on an surface which is of {\em arithmetic type} and when $\varphi_{\lambda}$ varies amongst a suitable orthogonal basis of Laplacian eigenfunctions. We give definitions of these concepts below.

 In such a situation, we obtain a much stronger result: {\em let $\{ \varphi_{\lambda}\}$ be a basis of Hecke-Laplace eigenfunctions. Then there is an absolute constant $\delta>0$ such that}
\begin{equation} \label{bettergeodesicbound} \frac{ | a(\varphi_{\lambda}, \omega) |}{ \| \varphi_{\lambda} \|_{L^2(M)}}  \leq C(M, \geod) \lambda^{-\delta}. \end{equation}
The bound on the right hand side is independent of $\omega$, and includes the ``difficult'' case
when $\omega$ and $\lambda$ are close. 
Thus the Fourier coefficients of an Laplace/Hecke eigenform along a fixed geodesic decay {\em uniformly}: such uniformity is a direct consequence of the {\em hybrid} nature of the subconvex bound proven in Theorem \ref{RSthm}.

  If $M$ is an arithmetic {\em hyperbolic} surface, then it is expected that \refs{bettergeodesicbound} holds for {\em any} orthonormal basis of Laplace eigenforms. This comes from the fact that, in the hyperbolic case, the multiplicities of Laplace eigenvalues are expected to be small, so that any Laplace eigenfunction could be expressed as a short linear combination of Laplace/Hecke eigenfunctions.

\subsubsection{The definition of arithmetic hyperbolic manifold}\label{arithHdefn} 

By an {\em arithmetic hyperbolic manifold} we shall mean  the quotient of the upper half-plane by a lattice that arises from a quaternion algebra over $\mathbb{Q}$.    However, for 
simplicity of exposition, we restrict ourselves to a slight subclass of $(M , \geod)$; we describe
this subclass here, and also translate some of our data into automorphic language. We strongly suggest skipping this section at a first reading.

Let $D$ be a quaternion algebra over $\mathbb{Q}$, split at $\infty$.   
Let $\mathbf{G}$ be the algebraic group $\GL_1(D)/\mathbb{G}_m$, 
so that the $\Q$-points of $\mathbf{G}$ are $D(\Qq)^{\times}/\mathbb{Q}^{\times}$. 
 Let $\adele$ resp. $\adele_f$ be the ring of adeles resp. finite adeles of $\mathbb{Q}$. 
Let $K_f$ be an open compact subgroup of $\mathbf{G}(\adele_f)$ with the property
that $\mathbf{G}(\Q) \cdot \mathbf{G}(\R) \cdot K_f = \mathbf{G}(\adele_f)$.  
Writing $\Gamma $ for $\G(\Q) \cap K_f$, we have a natural homeomorphism
from $\Gamma \backslash \mathbf{G}(\R)$ to
 $\mathbf{G}(\Q) \backslash \mathbf{G}(\adele) / K_f$. 
Since $\mathbf{G}(\R) \cong \PGL_2(\mathbb{R})$, it acts on $\mathbb{H}^2$ (although
only the connected component preserves orientation). We refer to the quotient
$\Gamma \backslash \mathbb{H}^2 $ as an arithmetic hyperbolic manifold. 

We therefore have a projection: 
 $$ \G(\Q) \backslash \G(\adele) = \adele^{\times} D(\Qq)^{\times} \backslash D(\adele)^{\times} \longrightarrow \Gamma \backslash \mathbb{H}^2,$$

For simplicity of exposition, we shall also restrict our discussion to geodesics $\geod$
that arise as a projection of  of a full adelic orbit $\left( \mathbf{T}(\Q) \backslash \mathbf{T}(\adele) \right). g$ 
to $\Gamma \backslash \mathbb{H}^2$, where 
 $\mathbf{T} \subset \mathbf{G}$ is a maximal torus and $g \in \mathbf{G}(\adele)$. 
(In general,   such a projection is the union of $\geod$ with finitely many closed geodesics, the number of such geodesics being the class number of a suitable quadratic order; this restriction, therefore, amounts to the requirement that this class number is $1$;  \refs{bettergeodesicbound} however remains true without such requirement). 

We now associate automorphic data to our eigenfunctions and frequencies:
\begin{itemize}
 \item[-] Let $\pi_{\lambda}$ be the automorphic representation of $\mathbf{G}$
generated by the pull-back of $\varphi_{\lambda}$.  

\item[-]  The torus $\mathbf{T}(\Q)$ 
is of the form $E^{\times}/\Q^{\times}$, where $E$ is 
a real quadratic field extension of $\Q$.  
\item[-] Associated to the character $\gamma(t) \mapsto e^{i \omega t}$ is a character of $\mathbf{T}(\Q) \backslash \mathbf{T}(\adele)$, and, in particular, a character $\omega_E$
of $\adele_{E}^{\times}/E^{\times}$. 
Let $\pi_{\omega}$
be the automorphic representation of $\GL_2$ over $\mathbb{Q}$ obtained
by automorphic induction from $(E, \omega_E)$; thus, the $L$-function of $\pi_{\omega}$
coincides with the $L$-function $L(s, E, \omega_E)$. 
\end{itemize} 
 
By a result of Waldspurger,   $|a(\varphi_{\lambda}, \omega)|^2$
is proportional to the central value of the completed $L$-function $\Lambda(\pi \times \pi_{\omega},\frac{1}{2})$.

\subsubsection{Sketch of the proof of \refs{bettergeodesicbound}}

As mentioned above, \refs{bettergeodesicbound} is a consequence of Theorem \ref{RSthm}. We will now outline a proof of this corollary, for arithmetic hyperbolic surfaces, in purely geometric terms.  

The function $t\mapsto \varphi_{\lambda}(\gamma(t))$ oscillates over a length scale of size $\lambda^{-1}$. 
It therefore stands to reason that the most interesting case of \eqref{bettergeodesicbound}
is when $\omega \sim \lambda$. This is indeed so, and the key step of the proof 
of Theorem \ref{RSthm} corresponds -- in this present language -- 
to the use of certain identities to switch from the ``difficult'' range when $\omega \sim \lambda$
to the ``easier'' range when $\lambda$ is small and $\omega$ is large. 
{\em These identities are extremely specific to the arithmetic cases under consideration;}
we do not know how to prove anything like \eqref{bettergeodesicbound} for a general
hyperbolic surface $M$.

Notation as previous.  One may 
construct (for ``deep'' number-theoretic reasons) the following auxiliary data:
\begin{itemize}
\item[-] another  {arithmetic hyperbolic surface} $\widetilde M$, depending only on $M$; 
\item[-] a Laplace/Hecke eigenfunction $\tilde\varphi_{\lambda}$ on $\widetilde M$ of eigenvalue $-\lambda^2$; 
\item[-] For each admissible $\omega$, i.e each integral multiple of $2 \pi/L$, we associate a Laplace eigenfunction $\theta_{\omega}$ on $\widetilde M$ with eigenvalue $-\frac{1}4-\omega^2$.  
\end{itemize}
moreover, whenever $\omega_{1},\omega_{2}, \frac{\omega_1 + \omega_2}{2}$ are admissible, one has
\begin{equation}\label{basic}\left| a(\varphi_{\lambda}, \omega_1) a(\varphi_{\lambda}, \omega_2) \right|^2 
\sim \left|  \int_{\widetilde M} \tilde\varphi_{\lambda} \theta_{\frac{\omega_1 + \omega_2}{2}} \theta_{\frac{\omega_1 - \omega_2}{2}}\right|^2.\end{equation}

The $\sim$ here does not indicate approximate equality, but rather, equality up to a constant
that is precisely computable; it is essentially a ratio of $\Gamma$-functions. 

The identity \eqref{basic} has the remarkable feature that the left-hand side has a ``quadrilinear'' nature,
whereas the right hand side has a ``bilinear'' nature.  Thus the map $\varphi_{\lambda} \mapsto 
\tilde\varphi_{\lambda}$ is in no natural sense linear; rather, it is defined element-by-element
over a special basis of Laplace/Hecke eigenfunctions.%

\subsubsection{Number-theoretic explanation} Supposing $\omega_{1},\omega_{2}, \frac{\omega_1 + \omega_2}{2}$ admissible,  
$$\Lambda(\pi \times \pi_{\omega_{1}},\frac{1}{2}) 
\Lambda(\pi \times \pi_{\omega_{2}},\frac{1}{2}) = \Lambda(\pi  \times\pi_{\omega_{+}}
\times \pi_{\omega_{-}},\frac{1}{2}),\ \omega_{\pm}:=\frac{\omega_1 \pm \omega_2}{2}.
$$
  \eqref{basic} now follows using the main result of \cite{Ich}: 
The manifold $\widetilde M$
is a quotient of the unique quaternion algebra $D'$ which is nonsplit
at those places where $\epsilon_v(\pi \times \pi_{\omega_{+}}
\times \pi_{\omega_{-}}) = -1$ -- if no $D'$ exists, the left-hand side of \eqref{basic} vanishes -- and $\widetilde{\varphi}_{\lambda}, \theta_{\omega_+}, \theta_{\omega_-}$
belong to the Jacquet-Langlands transfer to $D'$ of $\pi, \pi_{\omega_+}, \pi_{\omega_-}$.

\subsubsection{The switch of range from $\omega \sim \lambda$ to $\lambda=O(1)$}
Take $\omega_1 = \omega_2 = \omega$, and apply Cauchy-Schwarz to the right-hand side of \eqref{basic} to obtain:
$$|a(\varphi_{\lambda}, \omega)|^4  \lesssim \|\tilde{\varphi}_{\lambda}\|_{L^2(\widetilde M)}^2 \langle \theta_{\omega}^2 , \theta_0^2 \rangle_{L^2(\widetilde{M})}. $$

We shall analyze this by expanding both $\theta_{\omega}^2$ and $\theta_0^2$
into constituents. Since $\theta_0$ is itself a Laplacian eigenfunction with eigenvalue $1/4$, all $\Delta$-eigenfunctions that occur in the spectral expansion of $\theta_0^2$
will have ``small'' eigenvalue. Carrying this out\footnote{In the cases we encounter,
there will be difficulties with convergence, and regularization of the following
expression will be needed.}
\begin{equation} \label{basic2} \langle \theta_{\omega}^2 , \theta_0^2 \rangle =   \sum_\stacksum{\Delta \psi_{\mu} = \mu^2}{\mu\ll 1} \langle \theta_{\omega}^2, \psi_{\mu} \rangle \langle \psi_{\mu}, \theta_0^2 \rangle,\end{equation}
where the $\psi_{\mu}$-sum ranges over a basis of Hecke-Laplace eigenfunctions
on $\widetilde M$.

It is, in fact, possible to now apply \eqref{basic} once more to understand the term $\langle \theta_{\omega}^2, \psi_{\mu} \rangle$; however, we apply it ``in the reverse direction.''
This shows that there exists a manifold $\check M$ and $\check\psi_{\mu}$ so that 
$$| \langle \theta_{\omega}^2, \psi_{\mu} \rangle|^2 \sim |a( \check\psi_{\mu},2\omega)a( \check\psi_{\mu},0)|^2.$$

We have achieved our objective and {\em switched the range}:
starting with the analysis of $a(\varphi_{\lambda}, \omega)$ with $\lambda,\omega$ essentially arbitrary,
we have reduced it to the analysis of $a(\check\psi_{\mu}, 2 \omega)$ where $\mu$ may be assumed small relative to $\omega$.
 Although we are not done, this allows us to out-flank the most tricky
case of the question: when $|\lambda - \omega| = O(1)$. 
\subsubsection{The range when $\omega$ is large}
We now discuss bounding the Fourier coefficient $a(\varphi_{\lambda}, \omega)$ 
when $\omega$ is very large compared to $\lambda$.   
 For simplicity, let us assume in the present section that we are dealing with a {\em fixed} eigenfunction $\varphi = \varphi_{\lambda}$, and analyze the question of bounding $a(\varphi_{\lambda}, \omega)$ as $\omega \rightarrow \infty$.

 We need to be more precise about what is necessary to prove. In this context, it is evident (by real-analyticity) that $a(\varphi_{\lambda}, \omega)$ decays exponentially with $\omega$; however,
 this is not enough. For our previous argument we require:\footnote{The exponential factors arise from, in essence, the $\Gamma$-functions that were suppressed in \eqref{basic}.} 
\begin{equation} \label{needy} e^{\pi |\omega|/2} \left| a(\varphi_{\lambda}, \omega) \right| \leq C(\varphi_{\lambda}) |\omega|^{-\delta},\end{equation} 
for some $\delta > 0$.

Now, $a(\varphi_{\lambda}, \omega)$ is the integral of $\varphi_{\lambda}(\gamma(t))$
against $e^{i \omega t}$. To eliminate the exponential factors, we {\em deform the path}: we replace $t \mapsto \gamma(t)$
by a path $t \mapsto \gamma'(t)$ so that, first of all, $\gamma'(t)$ approximates -- at least locally -- a horocycle; secondly, the deformed integral still determines the original integral,
but is larger than it by a factor of size $e^{\pi |\omega|/2}$. Thus, to prove \eqref{needy},
we need only prove polynomial decay for the deformed integral.
The deformed integral is analyzed using dynamical properties of the horocycle flow,
especially mixing; it is related to the analysis in \cite[1.3.4]{Ve}.

\section{Norms on adelic quotients}

The classical Sobolev norms on $\R^n$, or on a real manifold, measure
the $L^p$- norms of a function {\em together with its derivatives}.  For example,
let $\|f\|_{2,k}$ be defined as the sum of the $L^2$-norms of the first $k$ derivatives
of $f \in C^{\infty}(\R/\Z)$, and let $S_{2,k}$ be the completion with respect to this norm. Then:
\begin{enumerate}
\item[Sa.] Sobolev inequality:  The Sobolev norms control point-evaluation, e.g. $$|f(0)| \ll \|f\|_{2,1};$$ 
\item[Sb.] Distortion: If $h : \R/\Z \rightarrow \R/\Z$ is a diffeomorphism, then $$\|f\circ h\|_{2,1} \leq
(\sup_{x} |h'(x)| ) \|f\|_{2,1}.$$ 
\item [Sc.]  Sobolev embedding: The Sobolev norms are compact with respect to each other:
the inclusion of $S_{2,k}$ into $S_{2,k'}$ is compact for $k < k'$. 
\item [Sc.*] Sobolev embedding in quantitative form: if $k' \geq k +2$,
then the trace of $S_{2,k}$ with respect to $S_{2,k'}$ is finite. This means
that, if $W_k$ is the completion of $C^{\infty}(\R/\Z)$ with respect to $S_{2,k}$, 
then the induced homomorphism $W_k \rightarrow W_{k+2}$ is trace-class.

\item[Sd.] Fourier analysis: $$|S_{2,k}(f)|^2 \ll \int (1+|\lambda|)^{2k} |\hat{f}(\lambda)|^2,$$
where $\hat{f}(\lambda)$ is the Fourier transform.
\end{enumerate}

It is very convenient to have a system of norms on adelic quotients with corresponding properties. 
We shall present them in terms of a list of axiomatic properties they satisfy, before giving the definition (\S \ref{sobolev}). 
These properties are intended to be analogous to (Sa) -- (Sd) above.  Prior to doing this,
we need to first recall $L^2$-spectral decomposition (\S \ref{plancherel}). 

We strongly recommend that the reader {\em ignore the definition} of the Sobolev norms
and rather {\em work with its properties.} 

\subsection{Notation}

\subsubsection{On implicit constants} \label{implicit} 

We use throughout the notation $A \ll B$ of Vinogradov to mean: there exists
a constant $c$ so that $A \leq c B$. If we write $A \ll_{\delta} c B$, it means
that the constant $c$ is permitted to depend on $\delta$, and so on. 
We shall also use a modification of this notation: $A \ll B^{\star}$ means
that there exist constants $c_1, c_2$ so that $A \leq c_1 B^{c_2}$. 

Although the notation $A \ll B$ is generally understood to mean that the implicit constant is {\em absolute}, 
it is 
extremely convenient in our context to allow it to depend on certain predetermined parameters (e.g., the number field over which we work) without explicit comment. 
We gather together at this point references to where these conventions
are introduced, for the convenience of the reader. 
To wit: in Part II, constants may depend on the isomorphism class of $(\mathbf{G}, \rho)$ over the number field $F$; 
in Part III, the constants may depend on the discriminant of the local 
field (\S \ref{implicitconstants}); in Parts IV and V, they may depend on the isomorphism class
of the ground field $F$.

  For Parts IV and V, we shall in fact make a more stringent use of the notation
where we require certain implicit constants to be polynomial; see \S \ref{ipnot}.

Later in the text we shall use indexed families of norms -- the Sobolev norms $\Sob_d$. 
They will depend on an indexing parameter $d$, as the notation suggests. 
In a similar fashion to the Vinogradov convention, we allow ourselves to write inequalities omitting the parameter $d$; see \S \ref{norm-intro} and \S \ref{SSmod-1} for a further discussion of this point.

\subsubsectionind{} 
Let $F$ be a number field. 
We denote by $\xi_F$ the complete $\zeta$-function of $F$. 
It has a simple pole at $1$; the residue is denoted by $\xi_F^*(1)$.

\subsubsectionind{} 
Let $\G$ a reductive algebraic $F$-group. Choose a faithful representation $\G\hookrightarrow \SL(F^{r})$ for some $r\geq 1$;
we shall suppose that it contains a copy of the adjoint representation. Henceforth
we shall feel free to identify $\G$ with a matrix group by means of this embedding. 
In the case of $\GL_2$, we shall fix the faithful representation to be
$$\rho: g \mapsto \left( \begin{array}{cc} g & 0\\ 0 & (g^{t})^{-1} \end{array}\right) \in \SL_4, $$
where, as usual, $g^t$ denotes the transpose.

Let $\mfg$ be the Lie algebra of $\G$; if $v$ is a place 
of $F$, we write $\mfg_v =\mfg\otimes_F F_v$.
  We also fix a basis for $\mathfrak{g}$. In that way, we regard
  the adjoint embedding as a map $\Ad: \mathbf{G} \rightarrow \GL(\dim \mathfrak{g})$. 

  Let $$\bfX := \G(F) \backslash \G(\adele).$$    We denote by $C^{\infty}(\bfX)$ the space of smooth functions on $\bfX$; the space $\bfX$ 
can be understood as an inverse limit of quotients of real Lie groups by discrete subgroups,
and a smooth function simply means one that factors through a smooth function on one of these quotients. 

\subsubsectionind{}  \label{basicnotn:TB}
We fix a maximally $F$-split torus $\mathbf{T} \subset \G$
and, correspondingly, a minimal parabolic $F$-subgroup $\mathbf{B}$ containing $\mathbf{T}$. 

\subsubsectionind{} 
 For $v$ non-archimedean and $m\geq 0$ an integer, we denote by
$K_{v}[m]$ the open-compact (principal congruence) subgroup
$$K_{v}[m]:=\G(F_{v})\cap \bigl\{g\in \GL_r(\mcO_{v}),\ g\equiv\Id_{r}\ (\varpi_{v}^m)\bigr\},$$
where $\varpi_v$ is a uniformizer in $F_v$. 

Choose, for each place $v$, a maximal compact
subgroup $K_v$ so that :
\begin{enumerate}
\item $ K_v \supset K_v[0]$ when $v$ is nonarchimedean (this implies that
$K_v = K_v[0]$ for almost all $v$);
\item For $v$ nonarchimedean, $K_v$ is {\em special}, i.e.
it is the stabilizer of a special vertex in the building of $\G(F_v)$. 
\end{enumerate}
This entails, in particular, that if $P_  v$ is the set of $F_v$-points of any parabolic subgroup,
then $P_v K_v= \G(F_v)$. 

There exists a constant $A$, depending on $\G$ and the chosen faithful representation, 
so that 
\be\label{volumebound} (1+A/q_v)^{-1} \leq  \frac{ [K_v: K_v[m_v]] }{q_v^{m_v \dim(\G)}} \leq (1+A/q_v) .\ee

We denote by $\underline{m}:v\mapsto m_{v}$ any function on the set of places of $F$ 
to the non-negative integers, which is zero for almost all $v$. 
Write $\|\underline{m}\| = \prod_{v} q_v^{m_v}$ (we take $q_v = e = 2.718\dots$ for archimedean places); we note that $$| \{\underline{m}: \|\underline{m}\| \leq N \}|= O(N).$$
For such $\underline{m}$, we set $$K[\underline{m}] := \prod_{v  \ \mathrm{ finite}} K_v[m_v].$$
We also put $K = \prod_{v} K_v$.

\subsubsectionind{} \label{measurenormalization} 
We now fix a normalization of left-invariant Haar measures on the various groups:
the $F_v$ and $\adele$-points of $\G$, as well as any parabolic subgroup $\mathbf{P} \subset \mathbf{G}$, as well as of any Levi factor $\mathbf{M} \subset \mathbf{P}$; and finally, on $K_v$, for every $v$.
These measures should have the following properties: 

\begin{enumerate}
\item The measures on adelic points should be the product of local measures;
\item The projection $\mathbf{P}(F) \backslash \mathbf{P}(\adele) \rightarrow \mathbf{M}(F) \backslash \mathbf{M}(\adele)$ corresponding to the decomposition $\mathbf{P} = \mathbf{M} \mathbf{N}$ is measure-preserving (for the {\em left} Haar measures!) 
\item For all $v$, the measure on $K_v$ has mass $1$. For nonarchimedean $v$, it is the restriction
of the measure from $\G(F_v)$. 
\item 
The map $\mathbf{P}(F) \backslash \mathbf{P}(\adele) \times K \rightarrow \G(F) \backslash \G(\adele)$ should be measure-preserving (i.e., the preimage of any set has the same measure as the set). 
\end{enumerate}

It is not difficult to construct such measures; we shall make an explicit choice, in the case of $\G = \GL_2$,  in \S \ref{measures-local}.  It is worth observing that, since we are only concerned
with upper bounds in this paper, and not exact formulas, precise choices are never of importance,
so long as they remain consistent. 
With any such choice, it follows from  \eqref{volumebound} that $\|\underline{m}\|^{-\eps} \ll \vol(K[\underline{m}])  \|\underline{m}\|^{\dim(\G)} \ll \|\underline{m}\|^{\eps}$; in fact, one can replace the upper and lower bounds by constants in the case that $\G$ is semisimple. 
 
Finally, put on  $\bfX$  the corresponding quotient measure. 

If $H$ is any locally compact group, we define the {\em modular character}
$\delta_{H}: H \rightarrow \mathbb{R}^{\times}$ via the rule
$\mu(S h)= \delta_{H}(h)^{-1} \mu(S)$, where $\mu$ is a left Haar measure on $H$, and
$S $ is any set with $\mu(S)  > 0$.  In other words, if $d_l h$ is a left Haar measure, then
$d_l(h h') = \delta_H(h')^{-1} d_l h$. %

\subsubsectionind{} \label{redtheory}

For $g \in \G(F_v)$, we define $\|g\| = \sup_{ij} |\rho_{ij}(g_v)|$; in the adelic case, we take the product over all places. We set $\| \Ad(g)\|$ to be
defined as $\|g\|$, but with $\rho$ replaced by the adjoint embedding. 
(The ``functional'' difference between these two norms lies
in the fact that $\|\Ad(g)\|$ is invariant under the center, whereas $\|g\|$ is not.)

For $g \in \G(\Aa)$, define the {\em height} by
$$\height(g)^{-1} := \inf_{x \in F^r - \{0\}} \prod_{v} \sup_{i=1\dots r} |(\Ad(g).x)_i|_v.$$
This descends to a function on $\bfX$, and -- if the center of $\G$ is anisotropic -- the map $\height: \bfX \rightarrow \R_{\geq 1}$ is proper.   

\begin{Lemma*}
Fix $x_0 \in \bfX$. 
For any $x \in \bfX$, there exists $g \in \G(\adele)$ with
$x_0 g= x$ and $\|\Ad(g)\| \leq \height(x)^{\star}$. 

\end{Lemma*}
\proof This is well-known; see e.g. \cite{ELMV1}, footnote 15 for a proof in the case of real groups, from
which the stated result is easily deduced. 
\qed

\subsubsectionind{} 
Throughout this paper the phrase ``$\pi$ is a unitary representation of the group $G$''
will be understood to mean that the underlying space of $\pi$ is a Hilbert space
and $G$ acts by isometries on that space.\footnote{In some contexts, unitary representation
is used to mean ``unitarizable'', i.e., such that there exists {\em some} inner product.
We shall always understand it to mean that we have fixed a {\em specific} inner product.}

Let $V$ be a unitary representation of $\G(F_v)$. 
The space of smooth vectors $V^{\infty}$ is defined, in the case when $v$ is nonarchimedean, as the subspace of $V$
comprising vectors whose stabilizer is open; in the case when $v$ is archimedean,
it is that subspace for which the map $g \mapsto g.v$ defines a smooth map
from $\G(F_v)$ to $V$. {\em It is always dense in $V$.}

Let $V$ be a unitary representation of $\G(\adele)$. It factorizes as a tensor
product of unitary representations of $\G(F_v)$, and we define the smooth subspace
$V^{\infty}$ as the (image in $V$ of the) tensor product of the local smooth subspaces.

\subsection{Structure of adelic quotients and the Plancherel formula}
\label{plancherel}

 In the present section we are going to recall the ``Plancherel formula'' for $L^2(\bfX)$, that is to say, its decomposition into irreducible $\G(\adele)$-repre\-sentations.

\subsubsection{Eisenstein series} \label{Eisenstein}
There is a standard parameterization of the automorphic spectrum 
via the theory of ``Eisenstein series'' that we shall now recall. See also \cite{Arthur}
for a r{\'e}sum{\'e} of the theory, and \cite{MWDS} for a detailed treatment.

Let $\data$ denote the set of pairs $(\mathbf{M}, \sigma)$, where 
$  \mathbf{M}$ is a $F$-Levi subgroup of a $F$-parabolic subgroup, containing $\mathbf{T}$,  and let $\sigma$
be an irreducible subrepresentation of the space of functions on $\mathbf{M}(F) \backslash
\mathbf{M}(\adele)$, which  is ``discrete series'' in the following sense: all $f \in \sigma$ are square-integrable with respect
to the inner product $$f\mapsto \|f\|^2_{\sigma}=\int_{Z_{\bfM}(\Aa)\bfM(F)\bash \mathbf{M}(\Aa)} |f|^2,$$ where $Z_{\bfM}(\Aa)$ denote the center of $\mathbf{M}(\Aa)$.
 
We can equip $\data$ with a measure in the following way:  We write $$\data = \bigsqcup_{\mathbf{M}} \data_{\mathbf{M}},$$ indexed by Levis containing $\mathbf{T}$. We require that
for any continuous assignment of $\chi \in \data_{\mathbf{M}}$ to $f_{\chi}$ in the underlying space of $\chi$,
$$\int_{\bfM(F) \backslash \bfM(\adele)} \left|  \int_{\chi} f_{\chi}  d\chi  \right|^2 = \int_{\chi} \|f_{\chi}\|_{\sigma}^2 d\chi.$$ 
This uniquely specifies a measure $d\chi$ on $\data_{\mathbf{M}}$, and so also on $\data$.  

There exists a natural equivalence relation $\sim $ on $\data$:
declare $(\mathbf{M}, \sigma)$ and $(\mathbf{M}', \sigma')$ to be equivalent
if there exists $w$ in the normalizer of $\mathbf{T}$
with $\Ad(w) \mathbf{M} = \mathbf{M}'$ and $\Ad(w) \sigma =\sigma'$.
There is a natural quotient measure on $\data/\sim$.

For $\chi = (\mathbf{M}, \sigma) \in \data$, we denote by $\mcI(\chi)$
the unitarily induced representation $\mathrm{Ind}_{\mathbf{P}(\adele)}^{\mathbf{G}(\adele)} \sigma$, 
where $\mathbf{P}$ is any parabolic subgroup containing $\mathbf{M}$.  (Its isomorphism class is  -- not obviously -- independent of the choice of $\mathbf{P}$.) 
It consists of functions $\G(\adele) \rightarrow V_{\sigma}$
(where $V_{\sigma}$ is a vector space realizing the representation $\sigma$)
satisfying the transformation property $$f(pg) = \delta^{1/2}(p) \sigma(m_p) f(g),$$ where
$\delta$ is the modular character (cf. \S \ref{measurenormalization}) and $p \mapsto m_p$ the projection $\mathbf{P} \rightarrow \mathbf{M}$. 
We define a norm on $\mcI(\chi)$ by
$$\|f\|_{Eis}^2 := \int_{K} \|f(k)\|_{V_{\sigma}}^2dk,$$
where $K$ is equipped with the Haar probability measure. 

There exists a natural intertwiner
(the ``unitary Eisenstein series'', obtained by averaging over $\mathbf{P}(F) \backslash \mathbf{G}(F)$
and analytic continuation):
$$\mcI(\chi) 
\stackrel{\Eis}{\rightarrow} C^{\infty}(\bfX).$$

The map $\Eis$ is an isomorphism away from a set of parameters $\chi$ of measure zero; we call the latter the set of {\em singular} parameters.
For instance in the case of $\GL_{2}$ the parameters are pairs of unitary characters $(\chi^+,\chi^-)$ for $\bfM=\GL_{1}\times\GL_{1}\hookrightarrow\bfB$ (the standard Borel) and the singular ones are the ones for which $\chi^+=\chi^-$; for this reason, we shall define for this case a variant $\Eisreg$ of $\Eis$ which is non-zero (see \S \ref{singulareisenstein}.) 

In any case, to almost every $\chi \in \data$ is associated an automorphic representation -- 
the image of $\mcI(\chi)$ -- 
denoted $\Eis(\chi)$. The resulting automorphic representation depends only on the class of $\chi$
in $\data/\sim$.   If $\mathrm{Eis}$ is an isomorphism -- as is so
away from a set of $\chi$ of measures $0$ -- we equip $\Eis(\chi)$ with a norm by requiring $\mathrm{Eis}$ to be an isometry. 
 Whenever defined, we call this norm the {\em Eisenstein norm} on the space $\Eis(\chi)$.  \footnote{The terminology may be slightly misleading; if $\pi \subset L^2(\bfX)$ is, e.g., a cuspidal representation, then $(\mathbf{G}, \pi) \in \data$,
and the Eisenstein norm on $\pi$ is simply the restriction of the $L^2$-norm.}

\begin{Theorem*} (Langlands). 
The map
$\int_{\chi \in \data} \mcI(\chi) \rightarrow L^2(\bfX)$, 
defined by integrating the map $\Eis$, extends to an isometric isomorphism
of the Hilbert space $\int_{\data} \mcI(\chi) d\chi$ and $L^2(\bfX)$. 
\end{Theorem*}
For a slightly more precise formulation, we refer to \cite{Arthur}.

We shall call any automorphic representation that occurs as an $\Eis(\chi), \chi \in \data$,
a {\em standard} automorphic representation. Not every automorphic representation is standard. 
For example, every standard automorphic representation is abstractly unitarizable. 
The standard automorphic representations are precisely those needed for unitary decomposition. 

The set of standard automorphic representations will be denoted by $\automorphicdual$
and the push-forward of the measure $d\chi$ on $\data$ to $\hat{G}_{Aut}$ will 
be denoted by $d\muP$. 

\begin{Remark}\label{pointwise}The Plancherel decomposition of a function in, e.g., $C^{\infty}_c(\bfX)$ is pointwise defined;
explicitly, for $\varphi \in C^{\infty}_c(\bfX)$, we have the equality of continuous functions,
$$ \varphi = \int_{\chi \in \data} \sum_{f \in  \mathcal{B}(\chi)} \langle \varphi, \Eis(f) \rangle \Eis(f)\ d\chi,$$
where $\mathcal{B}(\chi)$ is an orthogonal basis for $I(\chi)$, and the right-hand
side is absolutely convergent. This is not a triviality; it follows, for example, from the results of W. M{\"u}ller
\cite{muller}; it may be that there is a more elementary proof also. 
\end{Remark}

\subsubsection{The canonical norm for $\GL_2$} \label{canGL2}
In the case of $\GL_n$ for {\em generic standard} representations -- although not necessarily cuspidal --  it is possible to give a simple description of a canonical norm
on the space of any standard automorphic representation. 
We explain for $\GL_2$; the general case is obtained by replacing the role
of $a(y)$ below with $g\in\GL_{n-1}\subset \GL_{n}$ embedded as usual.

 Suppose $\pi$ is generic.  Let $\Whit_{\pi}=\otimes_{v}\mcW_{\pi,v}$ be the Whittaker model of
$\pi$. 
There are two natural inner products that one can equip $\mcW_{\pi,v}$ with,
namely -- see \S \ref{measures-local} for the measure normalizations -- 

\begin{eqnarray} \label{localinner}
\peter{W_v, W'_v}  &&= \int_{F_v^\times} W_{v}(a(y))\ov{W'_{v}}(a(y)) d^{\times}y, \peter{W_{v},W'_{v}}_{reg}\\
&&=\frac{\int_{F_v^\times} W_{v}(a(y))\ov{W'_{v}}(a(y)) d^{\times}y }{\zeta_{v}(1)L_{v}(\pi,\Ad,1)/\zeta_{v}(2)}. \nonumber
\end{eqnarray}
The latter inner product has certain good normalization properties, eg. $$ \peter{W_{v},W_{v}}_{reg}\\ = |W_v(1)|^2$$
for almost all $v$.

 We define an inner product on $\mcW_{\pi}$, by its value on pure tensors $W=\prod_{v}W_{v}$.
\be\label{globalinner}\|W\|^2_{reg} := \Lambda^*(\pi, \Ad,1) \times\prod_{v}
\peter{W_{v},W_{v}}_{reg};\ee
where 
$$\Lambda^*(\pi, \Ad,1)=\lim_{s \rightarrow 1} \frac{\Lambda(\pi, \Ad,s)}{(s-1)^r},$$  with $\Lambda(\pi, \Ad,s)=\prod_{v}L(\pi_{v},\Ad,s)$ denotes the completed $L$-function  and $r$ is taken to be the largest non-negative integer for which the limit is nonzero. 
The regularized value $ L^*(\pi, \Ad,1)$ satisfy (\cite{GHL})
\be\label{GHL}
L^*(\pi,\Ad,1)=C(\pi)^{o(1)},\ \mathrm{as}\ C(\pi)\ra\infty. \ee

Finally, we define the {\em canonical norm} on the space of $\pi$ by the rule
\begin{equation} \label{canreg} \|\varphi\|_{can}^2 =  
 \whitintertwinerconstant \|W_\varphi\|_{reg}^2,  \end{equation}
   where $\varphi \mapsto W_\varphi$ is the usual intertwiner \eqref{AW}. 
The terminology is justified by

\begin{Lemmat}\label{lem:cannorm}
Suppose $\pi$ generic and standard. Then, for $\varphi \in \pi$, $ \|\varphi\|_{can}^2 =   \|\varphi\|_{L^2(\bfX)}^2$
if $\pi$ is cuspidal; and $\|\varphi\|_{can}^2 = 2 \xi_F(2)\|\varphi\|_{Eis}^2$
if $\pi$ is Eisenstein and nonsingular. 
\end{Lemmat}

The verification of this equality for cusp forms is a consequence of the Rankin-Selberg method (cf. \S \ref{Subsec:rs}); the 
Eisenstein case is detailed in \S \ref{Eis-fourier}.

\begin{Remark}   \label{normdef}
One unfortunate consequence (perhaps unavoidable) is that the canonical
norm does not always behave continuously in families.  It is possible
to have a family of automorphic forms $\phi(s)$ belonging to standard generic automorphic
representations so that $\varphi(s) \rightarrow \varphi(0)$ pointwise on $\mathbf{X}$,
but the canonical norms do not converge. This happens when the order of pole
of the adjoint $L$-function jumps, e.g. at singular parameters. 
\end{Remark}

\subsection{Norms on adelic quotients}\label{sobolevnorms}\label{sobolev}

\subsubsectionind{} \label{norm-intro}

We shall set up two families of Hilbert norms, valid for any $d  \in \mathbb{R}$:
\begin{itemize}
\item[-]  A notion of ``Sobolev norm'' $\Sob^V_d$ on any unitary $\G(\adele)$-representation $V$. 
\item[-] A finer notion of ``Sobolev norm'' $\Sob^{\bfX}_d$ for functions on $\bfX$.
\end{itemize}

We shall follow the following convention: If $\ell \in V^*$ is a functional and we write $|\ell(f)| \leq \Sob^V(f)$, 
without a subscript $d$, it means {\em there exists a constant $d$, depending
on only the isomorphism class of $\mathbf{G}$ over $F$,  so that}
$|\ell(f)| \leq \Sob^{V}_d(f)$. In particular, $\ell$ is continuous in the topology defined by the family
of norms $\Sob^V_d$.

\begin{enumerate} 
\item We would like to warn the reader that the constructions are not totally formal. Namely,
some of the subtler features of the norms rely on Bernstein's uniform admissibility theorem as
well as M{\"u}ller's theorem \cite{bernstein,muller}. Indeed, one of the properties of the norms is established
only for the group $\G = \GL_n$. In fact, for the purpose of the present paper,
none of these deep results are important.

\item We observe that all our norms can take the value $\infty$. (In precise terms,
we understand a norm $N$ on a vector space $V$ to be a function $N: V \rightarrow [0, \infty]$ 
that satisfies the usual axioms; equivalently, we could regard $N$ to be
a pair consisting of a subspace $W \subset V$, and a (usual, finite-valued) norm on $W$.)
\end{enumerate}

$\Sob^V_d$ will always take finite values on $V^{\infty}$; 
$\Sob^{\bfX}_d$ will always take finite values on the space of compactly supported smooth functions
$C^{\infty}_c(\bfX)$. We shall sometimes refer to the completion of $V^{\infty}$ in the norm induced by $\Sob^V_d$ 
as the {\em Hilbert space associated to} $\Sob^V_d$. Similarly for $\Sob^{\bfX}_d$.

\subsubsection{The Sobolev norms on a unitary representation} \label{SSur}
Let $V$ be a unitary  admissible representation of either $\G(F_v)$ (some place $v$) or $\G(\adele)$. 
We shall define, in both contexts, a {\em generalized Laplacian operator}
$\Delta: V^{\infty} \rightarrow V^{\infty}$. 
This being so, we define the $d$th Sobolev norm
via \begin{equation} \label{Sobnormdef} \Sob^{V}_d(f) :=  \|\Delta^d f\|_V.
\end{equation}
The Laplacian will have the property that $\Delta$ is invertible, and a suitable
power of $\Delta^{-1}$ is trace class from $V$ to itself.

\subsubsection{Local}\label{sec:laplacelocal}
Let $v$ be a place and $V$ an unitary admissible representation of $\G(F_v)$. We shall
make a certain (Hilbert) orthogonal decomposition 
$V = \oplus_{m \geq 0} V[m]$; roughly speaking, vectors in $V[m]$ have ``higher frequency''
as $m$ grows. 

\begin{itemize}
\item[-] If $v$ is a finite place, we define $V[m]$ to be the orthogonal complement
of the $K_v[m-1]$-invariants vectors inside the $K_v[m]$-invariants vectors. 

\item[-]If $v$ is archimedean, fix a basis $\{X_i\}$ 
for the Lie algebra $\mfg_{v} := \mfg \otimes_F F_v$ and let $\mathcal{C} := \sum_{i}  (1- X_i^2)$;
let $V[m]$ be the direct sum of all $\mathcal{C}$-eigenspaces
with eigenvalue in $[e^m, e^{m+1})$ (here $e = \exp(1) = 2.718\dots$). 
\end{itemize}

The space $V[m]$ is finite dimensional. Write $e_v[m]$ for the projector onto $V[m]$, and put:
$$\Delta_{v}=\sum_{m\geq 0}q_{v}^me_{v}[m].$$
Note that $\sum_m e_v[m]$ is the identity.

\subsubsection{Global}

In the global setting, we use, as before,  the notation $\underline{m}$ for a function
$v \mapsto m_v$ from places to non-negative integers, and set 
$$e[\underline{m}] := \prod_{v} e_v[m_v],
\Delta_{\Aa} = \sum \|\underline{m}\| e[\underline{m}]$$
Note that one has $\sum_{\underline{m}} e[\underline{m}] = 1.$

\subsubsection{Sobolev norms on  $C^{\infty}(\quot)$}

 On $C^{\infty}(\quot)$ we shall introduce an increasing system of Sobolev norms $\Sob_d^{\quot}$
which will take into account the noncompactness of the space $\quot$. Let $H$
be the operation of ``multiplication by $1+ \height(x)$.'' Put:
$$\Sob_d^{\quot}(f) := \| H^d \Delta_{\Aa}^d f\|_2^2.$$

\subsection{Properties of the Sobolev norms.}
This section enunciates the properties of the Sobolev norms, defined in the prior section.

\subsubsection{Properties of the unitary Sobolev norms}

Write $G$ for either $\G(\adele)$ or $\G(F_v)$.  
\begin{enumerate}
\item[S1a.] \label{sobineq} (Sobolev inequality) For $f \in C(G)$ and 
a fixed smooth function $\omega: G \rightarrow \C$ of compact support with $\omega(1) = 1$, 
$$|f(1)| \ll_\omega \Sob^{L^2(G)}( f \cdot \omega). $$
\item[S1b.] \label{sobdist}(Distortion property) There
is a constant $\adconstant$, depending only\footnote{If the representation $\rho$ contains 
the adjoint representation, one may take $\adconstant=1$.} on $(\G, \rho)$, so that, for $V$  a unitary representation of $G$, $$\Sob_d^V(g f) \ll \|\Ad(g) \|^{\adconstant d} \Sob_d^V(f).$$
\item [S1c.](Embedding) For each $d$, there exists $d' > d$ so that the trace of $\Sob^V_d$ w.r.t.
$\Sob^V_{d'}$ is finite.  (See \cite{BR}*{Appendix} for definitions; this means
that the inclusion from the Hilbert space associated to $\Sob^V_d$, to
the Hilbert space associated to $\Sob^{V}_{d'}$, is trace-class. 
\item[S1d.] \label{soblocglob} (Linear functionals can be bounded place-by-place.) 
Let $\pi = \otimes \pi_v$ be a unitary representation of
$\G(\adele)$; let $\ell = \prod_v \ell_v, \ell_v  \in \pi_v^*$ be a factorizable functional
with the property that $|\ell(x_v)|\leq 1$ when $x_v \in \pi_v$ is spherical and of norm one.\footnote{If this condition is satisfied not for all $v$, but for all $v \notin T$, where $T$ is some fixed
finite set of places, then the constant has to be replaced by $A^{1+|T|}$.} 
Then: 
\begin{equation*}
 |\ell_v| \leq A \Sob^{\pi_v}_{d} \mbox{ for all $v \implies$} |\ell | \leq A' \Sob^{\pi}_{d'}\end{equation*}
where $d'$ depends on $d$, and $A'$ depends on $A,d$. 

\end{enumerate}

\subsubsection{ Properties of the $\quot$-Sobolev norms }
\begin{enumerate}
 \item[S2a.] \label{sobinfty}(Sobolev inequality) There exists $d_0$ so that $\Sob_{d_0}^{\quot}$ majorizes $L^{\infty}$-norms.  
\item[S2b.] \label{sobdist2} (Distortion property) There exists a constant $\adconstantX$,
depending only on $(\G, \rho)$, so that
 $$\Sob_d^\quot(g f) \ll \|\Ad(g)\|^{\adconstantX d} \Sob^\quot_d(f), f \in C^{\infty}(\quot), g \in \G(\adele).$$
\item[S2c.]  \label{sobtrace} (Embedding) 
For each $d$, there exists $d' > d$ so that the trace of $\Sob^{\quot}_d$ w.r.t. $\Sob^{\quot}_{d'}$
is finite.
\end{enumerate}

\subsubsection{  Relationship between unitary Sobolev norms on $L^2(\quot)$  and $\quot$-Sobolev norms }
\begin{enumerate}
\item[S3a.]
We have a majorization $$\Sob_d^{L^2(\quot)} \ll \Sob_d^{\bfX}.$$
\item[S3b.] 
$$\Sob^{\bfX}_d(f) \ll \Sob^{L^2(\bfX)}_{d'}(f), \ \  f \in L^2_{cusp}$$  where $d'$ depends on $d$,  
and $L^2_{cusp}$ is the cuspidal subspace of $L^2(\bfX)$. 
\footnote{A stronger and more natural statement is: the truncation operator $\wedge^T$ of Arthur is continuous from the $\Sob^{L^2(\bfX)}$-topology
to the $\Sob^{\bfX}$ topology.
This statement is actually a quantitative form of \cite[Lemma 1.4]{Arthur-truncation},
which proves the same result but with the finite level fixed.} 

\item[S3c.]
 Let $\ell$ a linear form on $ C^{\infty}(\bfX)$. Suppose, for each standard automorphic representation
$\pi$, we have $\left| \ell|_{\pi^{\infty}} \right| \leq \Sob_{d}^{\pi}$,
where the unitary structure on $\pi$ is the Eisenstein norm. 
Then $$|\ell(v) | \ll \Sob^{\quot}_{d'}(v),$$ whenever both sides are defined; $d'$ depends only on $d$. 
\end{enumerate}

\subsubsection{Other properties} \label{interpolation}

Suppose $\ell$ is a linear functional on either $C^{\infty}(\bfX)$ or a unitary $\G(\Aa)$-representation,
and let $A(d)$ be the operator norm of $\ell$ with respect to $\Sob_d$. 
\begin{enumerate}
\item [S4a]
$\log A(d)$ is convex with respect to $ d$. \footnote{This will be helpful in bounding operator norms
when $d$ is not integral, for, like the case of $\R^n$, the Sobolev norms are most accessible when $d \in \mathbb{Z}_{\geq 0}$. }
\item [S4b.]
Given $d$ and $j \geq 1$, there exists $d' > d$ and common orthogonal bases
$e_1, e_2, \dots$ for $\Sob_{d}, \Sob_{d'}$ with the property that 
$\frac{\Sob_{d}(e_i)}{\Sob_{d'}(e_i)} \leq (1+|i|)^{-j}$. 
\item [S4c.]  $\Sob^V_{d}$ and $\Sob^{V}_{-d}$ are self-dual;
\item [S4d.]
If $\ell: V \rightarrow W$ is any linear functional from $V$ to another normed vector space:
\begin{eqnarray}\label{rt1}\ell(v) \leq A \| \Delta^ d v \|_V&&( v \mbox{ a $\Delta$-eigenfunction)} 
\\ &&\implies\  |\ell(v)| \leq A' \Sob^V_{d'}(v),  \ \ \mbox{ all $v \in V$, some $d' > d$.}\nonumber \end{eqnarray}
\item[S4e.] $\Sob^{\bfX}_d(f g) \ll_d \Sob^{\bfX}_d(f) \Sob^{\bfX}_d(g)$.  
\end{enumerate}

(S4d) is a consequence of the fact that a suitable power of $\Delta^{-1}$ is trace class.

\subsection{Examples.}

We shall now discuss a number of examples of using the Sobolev norms to quantify
mixing or uniform distribution. The essence of all our examples is well-known. We simply
want to make the point that the axioms of Sobolev norms make it simple to derive the results in great generality. 

In the sequel we will use the following notation: for $\Sob_{?}^\star$ one of the families of Sobolev norms discussed previously, we will write
$A\leq \Sob^{\star}(f)B$ to mean that there is a constant $d>0$ so that $A\leq \Sob_{d}^{\star}(f)B$.

\subsubsection{Bounds for matrix coefficients and canonically normalized functionals: local bounds} \label{mc-local}
For any place $v$, we denote by $\Xi_v(g)$ the Harish-Chandra spherical function on 
$G_v := \G(F_v)$, i.e. $\Xi_v(g) := \langle g v, v \rangle$ where
$v$ is the $K_v$-invariant function in the representation of $G_v$ unitarily induced
from the trivial character of a  minimal parabolic subgroup. 

It follows from \cite{CHH} (as extended to reductive groups in \cite{Oh}) that, if $\pi_v$ is a tempered representation of $\G(F_v)$
 one has, for any $x_1, x_2 \in \pi_{v}, g\in G_v$
 and some constant $d\geq 0$ depending on $\G$ only
$$| \peter{g x_1, x_2}|  \leq A_v
\Sob^{\pi_{v}}_d(x_1) \Sob^{\pi_{v}}_d(x_2) \Xi_v(g),$$
where $A_v = [K_v: K_v[0]]$ is equal to $1$ for almost all $v$. \footnote{Indeed, if $x_1, x_2$
are both stabilized by a subgroup of $K_v[0]$ of index $\heartsuit$, the bound may be taken to be
$\heartsuit \|x_1\| \|x_2\| \Xi_v(g)$; we shall use this later. }

Indeed, \cite{CHH}
implies, in fact, that $| \langle g x_1, x_2 \rangle | \leq [K_v: K_v[m]] \Xi_v(g)$
for $x_1, x_2 \in V[m]$.   Now \eqref{volumebound} and \eqref{rt1}
establish the desired statement.

In the case of non-tempered representations, one has a corresponding, but weaker, bound. 
In this paper we will only use the special case of $\bfG=\GL_{2}$; and moreover,
 $\pi_{v}$ will always be a local constituent at $v$ 
of a generic automorphic representation of $\bfG(\Aa)$. 
In that case, one has
\be\label{localglobalbounds}| \peter{g x_1, x_2}|  \leq A_v
\Sob_d^{\pi_{v}}(x_1) \Sob_d^{\pi_{v}}(x_2) \Xi_v(g)^{1-2\theta},\ee
for some absolute constant $\theta<1/2$. This is due to Selberg for $F=\Q$,
and \cite{GJ} in general.  As of now, it follows from the work of
Kim and Shahidi \cite{KSh}, that \refs{localglobalbounds} holds for $\theta \geq  3/26$. 

It will be convenient, for the purpose of this document, to allow $\theta$
to be any number in $]3/26, 1/4[$. More precisely, any statement involving $\theta$
will be valid for any choice of $\theta$ in this interval. This notational convention
will suppress $\epsilon$s at a later point. The reader may safely substitute
$3/26+\epsilon$ every time he/she sees the symbol $\theta$.

\subsubsection{Bounds for matrix coefficients and canonically normalized functionals: global bounds} 
\label{MC0}
Let $\iota:\tilde{\G}\mapsto \G$ be the simply connected covering of $\G$. 
For $f \in C^{\infty}(\bfX)$, set 
\be\label{locconst}\scrP{f}(x) = \int_{g \in \tilde{\G}(F)\backslash \tilde{\G}(\adele)}
f(\iota(g) x)dg\ee for $x \in \bfX$; here $dg$ is the invariant probability measure
on $\tilde{\G}(F) \backslash \tilde{\G}(\adele)$.  The endomorphism $\scrP$ realizes the orthogonal projection onto locally constant functions on $\bfX$, and is in particular $L^2$- and $L^{\infty}$-bounded; note also that it commutes with the $\bfG(\Aa)$-action
(cf. \cite{Ve}).

Take $g \in \G(\adele)$ and $f_1, f_2 \in C^{\infty}(\quot)$. 
Then there exists  $\beta > 0$ and $d$, both depending only on the isomorphism class of $\mathbf{G}$ over $F$, so that
\begin{equation} \label{MC} \left| \langle g. f_1, f_2 \rangle - \langle g. \scrP{f_1}, \scrP{f_2}  \rangle \right| \ll \|\Ad(g)\|^{-\beta} \Sob^{\bfX}_d(f_1) \Sob^{\bfX}_d(f_2).\end{equation}

This is a consequence of property $(\tau)$, which has been established through the work of many people; the proof was completed (together with the most difficult case) by L. Clozel \cite{clozel}.  Property $(\tau)$
was put in a quantitative form in \cite{GMO}, and the quoted statement\footnote{It should be noted that we could not follow the quantitative arguments of \cite{GMO} in the case which relies on Clozel's work;
however, their argument certainly furnishes a bound.}   follows from 
Lemma 3.3 and Theorem 3.10 of {\em loc. cit.} 
Here we will need \refs{MC} only for the case of $\G= \GL_2$: it follows from \refs{localglobalbounds} that the exponent
 $\beta = 1/2-\theta$ is admissible, any $\theta$ as above.

\begin{Rem}  
Let us explain one of the reasons that \eqref{MC} is so useful in the context of the present paper.  It has been observed J.-L. Waldspurger \cite{Wald} and greatly extended by Ichino-Ikeda \cite{Ich,IchIk}
that a wide variety of period functionals can be expressed by integrating matrix coefficients;
roughly speaking, there exists a variety of pairs $(\mathbf{H} \subset \mathbf{G})$ so that, with suitable choice of measure {\em and suitable regularization}, 
\begin{equation} \label{Canonical} \left| \int_{\mathbf{H}(F) \backslash \mathbf{H}(\adele)} \varphi\  dh\right|^2 = \int_{\mathbf{H}(\adele)} \langle h. \varphi, \varphi \rangle dh.\end{equation} 
The right-hand side is usually divergent, and is interpreted by a suitable regularization;
in many situations, almost every local factor is equal to the local factor of a suitable $L$-function,
which suggests a regularization involving a special value of that $L$-function. 

In any case, \eqref{Canonical} 
was  used in \cite{CU} and \cite{ELMV3} in order to reduce equidistribution results
for $\mathbf{H}(F) \backslash \mathbf{H}(\adele)$ -- or translates thereof -- to bounds for matrix coefficients. 
We shall use the same technique in \S \ref{sec:hjl} of the present paper. 

We note, however, that this technique is not universally applicable, and in many
instances the correct formulation of a result analogous to \eqref{Canonical} remains mysterious. 
\end{Rem}

\subsubsection{A quantitative form of the Ergodic principle II} \label{Prinzip-quant}

Next, let us give a quantitative form of the Ergodic principle \ref{Prinzip} from \S \ref{classical}.  

\begin{Lemma*}
Suppose $H \subset \G(\adele)$ is noncompact, and $\chi: H \rightarrow \mathbb{C}^{\times}$
a unitary character. Suppose that $\bfX$ has finite measure, i.e., that the center of $\G$ is anisotropic.\footnote{The lemma could be easily adapted to the general case by introducing a character; but we will only use it for $\G = \PGL_2$.}

 Let $\nu$ be a (possibly signed) $\chi$-equivariant measure on $\quot$, i.e.
 $\nu^h = \chi(h) \nu$ for $h \in H$. 
 Let 
$\mu$ be the $\G(\adele)$-invariant (Haar) probability measure, and suppose that, 
for some $d \geq 0$,  we have the majorization:
$$ |\nu|(f)  \ll  \mu(f)+  \epsilon \Sob_d^{\bfX}(f) \ \ (f \geq 0).$$

Let $\sigma$ be any probability measure on $H$. Then, for any $f$ with $\mathscr{P} f= 0$, 
$$ |\nu(f) - \delta_{\chi=1} \mu(f) |^2 \ll \left( \epsilon \|\sigma\|_{d'}^2 + \|\sigma \star \check{\sigma}\|_{ -\beta}  \right) \Sob^{\bfX}_{d'}(f)^2.$$ 
Here $\delta_{\chi=1}$ is $1$ if $\chi$ is trivial and zero otherwise,  $$\|\sigma\|_d := \int_{h} \|\Ad(h)\|^d d\sigma(h),$$ and similarly for $\|\sigma \star \check{\sigma}\|$, 
 $\check{\sigma}$ denote the pullback of $\sigma$ by $g\mapsto g^{-1}$,  $\beta$ is as in  \eqref{MC}, and $d'$ depends only on $d$. 
\end{Lemma*}

This is {\em indeed} a quantitative form of the Ergodic principle \ref{Prinzip}: if, e.g., $\epsilon = 0$, 
we see at once that that the left-hand side may be made arbitrarily small
to choosing $\sigma$ to have large support.

\proof We may assume that $\mu(f)=0$.
For $\sigma$  a probability measure on $H$ we set\ $$ f \star_{\chi} \sigma := \int_H  \chi(h) (h \cdot f) d\sigma(h).$$ We see:
$$|\nu(f)|^2= |\nu(f \star_{\chi} \sigma)|^2 \leq |\nu|(|f \star_{\chi} \sigma|^2) \leq \epsilon \Sob_d^\bfX(|f \star_{\chi} \sigma|)^2
+ \| f \star_{\chi} \sigma\|_2^2,$$
The required bound now follows from (S2b), (S4e) and from the bounds for matrix coefficients of \S \ref{MC}. 
\qed

\begin{Rem}In the context of the subconvexity problem for $L$-functions, this amounts to the {\em amplification method} 
(\cite{Iw}) and the measure $\sigma$ play the role of the {\em amplifier}.  
\end{Rem}

\subsubsection{Uniform distribution of horocycles} \label{horoUD}
As a final example, we discuss the uniform equidistribution of horospheres. 
A result of Peter Sarnak proves  that {\em a closed horocycle of length $L$} on a hyperbolic surface
becomes uniformly distributed as $L$ approaches $\infty$. What we present below
is an adelized version of the result for general groups, but only in the easiest version: where ``closed horocycle''
is replaced by (in more explicit language than what follows) ``closed orbit
of a maximal horospherical subgroup.''

Let us fix $\mathbf{N}$ to be the unipotent radical of a minimal proper parabolic subgroup $\mathbf{B}$ of $\mathbf{G}$; 
for $f \in C^{\infty}(\bfX)$ let $$f_N(g) := \int_{\mathbf{N}(F) \backslash \mathbf{N}(\adele)} f(ng) dn,$$ 
denote the constant term. Then for $f \in L_0^2(\bfX)\cap C^{\infty}(\bfX),$ and $b \in \mathbf{B}(\adele)$
\be\label{horcyclebound}
|f_N(b)| \ll \delta_{\bfB}(b)^{1/2} \Sob^\bfX(f),
\ee
here by $L_0^2(\bfX)$, we mean the orthogonal complement of locally constant functions (i.e. the kernel of \refs{locconst}), 
and $\delta_{\bfB}: \mathbf{B}(\adele) \rightarrow \mathbb{R}_{>0}$ is the modular character (cf. \S \ref{measurenormalization}).

 \proof 
By Property (3c) of Sobolev norms,  it is enough to verify \refs{horcyclebound} when $f$ belongs to an 
standard automorphic representation $\pi \subset L_0^2(\bfX)$. The constant term
is zero unless (notation of \S \ref{plancherel}) $\pi = \mathrm{Eis}({\mathbf{M}, \sigma})$, where
$\mathbf{M}$ is a minimal Levi subgroup.  Let $\mathbf{B}$ be a parabolic subgroup with unipotent radical $\mathbf{M}$. 
 {\em As an abstract unitary representation}, $\pi$
is isomorphic to the space $V_{\pi}$ of functions from $\G(\adele)$ to the representation space of $\sigma$, satisfying:
$$f(bg) = \sigma(b) \delta_{\bfB}(b) f(g), \ \  b \in \bfB(\adele),$$
equipped with the inner product $\int_K \|f\|^2$. (Here $K$ is the fixed maximal compact subgroup of $\G(\adele)$.) 

Let $\omega$ be a smooth compactly supported bi-$K$-invariant function on $\bfG(\adele)$ 
so that $\omega(1)=1$. 
Consider $\mathrm{res}: V_{\pi} \rightarrow L^2(\bfG(\adele))$
defined by $f \mapsto \omega f$. 
Apply (S1a) to see  $|f(1)|  \ll \Sob^{L^2(\bfG(\adele))} (\omega f)$. The map
$\mathrm{res}$ is $L^2$-bounded, and commutes with the action of the finite part of $K$;
moreover, if $\mathcal{D}$ is any element of the universal enveloping
algebra $\mathfrak{U}$ of $\mathfrak{g}$, then $\mathcal{D} \omega = \sum_{i} \omega_i \mathcal{D}_i$
for various $\mathcal{D}_i \in \mathfrak{U}$ and $\omega_i$ smooth, compactly supported on $\G(\adele)$. 
It follows
$\Sob^{L^2(\bfG(\adele))}(\omega f) \ll  \Sob^{V_{\pi}}(f).$
\qed 

\subsection{Proofs concerning Sobolev norms.}
\label{sobolevnormsproofs}

In this section, which can be safely skipped in the course of reading the paper,
we give the proofs that the Sobolev norms defined in \S \ref{sobolevnorms} 
have the good properties indicated in that section.   

The proof of the distortion properties, (S1b), (S2b) are elementary, and we omit them. 
Similarly we omit the proof of properties (S4a) -- (S4d). 

\subsubsectionind{} \label{1apar}
We will establish (S1a), from which (S2a) follows easily. 
We will establish (S1a) in the adelic setting; the local setting is easier. 

Let $\omega$ be as in (S1a). 
Take $f \in L^2(G) [\underline{m}]$, where, as usual, $\underline{m}$
is a function from places of $F$ to non-negative integers.

The function $f$ factors through the  $\dim(\G) [F:\Q]$-dimensional manifold\\
$\G(\adele)/K[\underline{m}]$; this has around $\vol(K[\underline{m}])^{-1}$
connected components. We can apply the usual Sobolev inequality (for real manifolds)
to the connected component of the identity.  It shows that $f(1)$ is bounded
by the $L^2$-norm, on the connected component of the identity, of
(e.g.) the first $\dim(\G)$ derivatives of $f\omega$. 

From this it follows that:
$$|f(1)| \ll  \vol(K[\underline{m}])^{-1} \prod_{v |\infty} \|m_v\|^{[F:\Q] \dim \G}\|f \omega\|_{L^2} \ll
\|m\|^{d_0} \|f\|_{L^2}
$$
where we take any $d_0 \geq  2 [F:\Q] \dim(\G)$.  This, in view of \eqref{rt1}, implies (S1a). 

To deduce (S2a), take $f \in L^2(\bfX)$. Fix $x_0 \in \bfX$. By pull-back of functions,
we immediately deduce from (S1a) that $|f(x_0)| \ll \Sob_{\linftyexp}^{L^2(\bfX)}(f)$. 
Now, by the Lemma of \S \ref{redtheory}, for any $x \in \bfX$, there exists $g \in \G(\adele)$ with
$x_0 g= x$ and $\|g\| \leq \height(x)^{\star}$. Therefore, the distortion properties imply that
\begin{equation} \label{S2a-inter} |f(x)| = |f^g(x_0)| \ll  \Sob_{\linftyexp}^{L^2(\bfX)}(f^g) \ll \height(x)^{\star} \Sob_{\linftyexp}^{L^2(\bfX)}(f).\end{equation}
Thus (S2a). 

\subsubsectionind{}
We prove (S3b).  It suffices to check that for any $m \geq 1$ there exists $d$ so that
\begin{equation}
 |f(x)| \ll_m (\height x)^{-m} \Sob^{L^2(\bfX)}_d(f), \ \ f \in L^2_{cusp}. \end{equation}

It suffices, by (S4d), to establish this estimate for $f \in L^2_{cusp}[\underline{m}]$.

Let $x \in \bfX$.  Let $\mathbf{T}$ be as in \S \ref{basicnotn:TB}. 
It is a consequence of reduction theory that there exists
a compact subset $\Omega \subset \G(\adele)$ and $R > 0$  
and $a \in \mathbf{T}(\adele)$ so that
$x = \G(F) . a \omega$, where $\omega \in \Omega$, and $a \in \mathbf{T}(\adele)$
is so that there exists at least one simple root $\alpha$ with 
 $|\alpha(a)|_{\Aa} \geq \height(x)^{\star}$. 
  Let $\mathbf{U}$ be the unipotent
radical of the maximal proper parabolic subgroup associated to $\alpha$.

We choose a sequence $\mathbf{Z}_0 = \mathbf{Z} \subset \mathbf{Z}_1 \subset \mathbf{Z}_j = \mathbf{U}$ where $\mathbf{Z}_{j+1}/\mathbf{Z}_{j}$
is central in $\mathbf{U}/\mathbf{Z}_j$ for $j \geq -1$ (we interpret $\mathbf{Z}_{-1}$ to be trivial). 
Now expand $f$ in a Fourier series along $\mathbf{Z}_0$, i.e.
$$f(x) = \sum_{\psi} f_{\psi}(x), 
f_{\psi}(x) := \int_{\mathbf{Z}(F) \backslash \mathbf{Z}(\Aa)} f(ux) \psi(u)du;$$
 the measure is the invariant probability measure, whereas the $\psi$
 sum extends on characters of $\mathbf{Z}(\Aa)$ trivial on $\mathbf{Z}(F)$. 

Only those characters $\psi$ which are trivial on $\mathbf{Z}(\Aa) \cap K[\underline{m}]$
contribute to this summation, and the number of such characters is bounded by
$\ll \|\underline{m}\|^{\dim \mathbf{Z}}$. 
Let $\mathfrak{u}$ be the Lie algebra of $\mathbf{Z}(F \otimes_{\Q} \R)$.    To any
$\psi$ as above, we may associate its differential $d\psi \in \mathfrak{u}^*$.
Fix a norm on $\mathfrak{u}^*$. Integrating by parts, we conclude that for any $n \geq 1$, 
$$\sum_{\psi \neq 1} f_{\psi}(x) \ll_n |\alpha(a)|_{\Aa}^{-n} \|m\|^{\star} \sup_{u \in \mathbf{U}(\adele), D \in \mathcal{D}}
|\mathcal{D} f(ux)|,$$
where $\mathcal{D}$ is a finite set of $\G(\Aa)$-invariant differential operators.  
Applying \eqref{S2a-inter}, we deduce that the sum is $\ll_n \height(x)^{-n} \|m\|^{a(n)} \Sob^{L^2(\bfX)}_{d(n)}(f)$.

It remains to estimate the term $f_{\psi}$ for $\psi=1$. For this, we expand it in a Fourier series
along $\mathbf{Z}_1/\mathbf{Z}_0$.  Proceeding in this way, we arrive inductively at the conclusion. 

\qed

\subsubsectionind{} \label{1bpar}
We prove (S1d). First of all, 

\begin{Lemma*} Let $\pi_v$ be an irreducible admissible representation of $\G(F_v)$.
Then the image of $e_v[m]$ has dimension bounded above by $A_1 q_v^{A_2 m}$, where $A_2$ depends only on the isomorphism class of $\G$. 
\end{Lemma*}

This follows (for $v$ nonarchimedean) from Bernstein's proof of uniform admissibility \cite{bernstein} and (for $v$ archimedean)
since any irreducible representation of a maximal compact subgroup occurs in $\pi_v$
with a multiplicity bounded by its dimension (a consequence of e.g. the subrepresentation theorem;
recall we interpret $q_v = e = 2.71828\dots$ for $v$ archimedean). 
We give, for illustration, Bernstein's proof for $\PGL_2$; this, and the mild variant
of $\GL_2$ with a fixed central character, is the only case in which we shall use
the result later in this paper. 
\proof (in the case of $\G = \PGL_2$).  Let $\mathcal{C}$  the characteristic function of $K[\varpi_v^m] \left( \begin{array}{cc} \varpi_v & 0 \\ 0 & 1 \end{array} \right) K[\varpi_v^m]\subset \GL_{2}(F_{v}).$ Then $\mathcal{C}$ acts on $\pi^{K[\varpi_v^m]}$,
the $K[\varpi_v^m]$-fixed vectors in $\pi$, 
as does the finite group $\mcG = \bar{K}_v / \bar{K}[\varpi_v^m]$. 
Let $\C[\mcG]$ be the algebra of endomorphisms of $\pi^{K[\varpi_v^m]}$ generated
by $\mcG$.  Let $\langle \mathcal{C} \rangle$ be the algebra of endomorphisms
of $\pi^{K[\varpi_v^m]}$ generated by $\mathcal{C}$. Then $\C[\mcG] \cdot
\langle \mathcal{C} \rangle \cdot \C[\mcG]$ must be the whole
algebra of endomorphisms of $\pi^{K[\varpi_v^m]}$, because a corresponding fact holds for the Hecke algebra
 of $K[\varpi_v^m]$-bi-invariant functions on $\GL_2(F_v)$, which acts irreducibly on $\pi^{K[\varpi_v^m]}$. 
 On the other hand,
$\dim \langle \mathcal{C} \rangle \leq \dim \pi^{K[\varpi_v^m]}$. Therefore,
the dimension of 
$\C[\mcG] \cdot
\langle \mathcal{C} \rangle \cdot \C[\mcG]$ is at most $|\mcG|^2 .  \dim \pi^{K[\varpi_v^m]}$;
since it acts irreducibly, its dimension is $=( \dim \pi^{K[\varpi_v^m]})^2$, yielding the result. 
\qed

Now, let $\pi$ be as in (S1d) -- i.e. a representation of $\G(\adele)$ --  and let $V$ be the underlying vector space. 
The Lemma shows that any $X \in V[\underline{m}]$ is a sum
$X=\sum_i  x_i$ where
where $x_i$ is a pure tensor, 
and the number of $x_i$s is at most $\prod_{m_v \neq 0} (A_1 q_v^{A_2 m_v}) 
\leq  \|m\|^{C_2}$, where $C_2 = A_2 + \frac{\log A_1}{\log 2}$. 
For any such pure tensor $x_i =\prod_{v} x_{i,v}$, we have (notation as in the statement
of [S1d]): 
$$| \ell(x_i)| =|  \prod_{m_v \neq 0}  \ell(x_{i,v}) \prod_{m_v = 0} \ell(x_{i,v}) |
\leq A' \|m\|^{d+\epsilon} \|x\|.$$

Here we used the fact that $\prod_{m_v \neq 0} A \leq A' \|m\|^{\eps}$, where $A'$
can be taken to a function of $A$ and $\epsilon$.  Thus an application of Cauchy-Schwarz yields:
$$|\ell(X)| \leq \|m\|^{C_2/2} A' \|m\|^{d+\epsilon}  \left( \sum_{i} \|x_i\|^2 \right)^{1/2}
\leq A'(\epsilon) \|\Delta_{\adele}^{d+C_2/2+\epsilon} X  \| $$ 
(S1d) now follows from \eqref{rt1}. 
\qed
\begin{Remark} \label{s1dremark}
It will be useful to apply (S1d) to multilinear functionals. Suppose 
that we are given (e.g.) a bilinear functional $B: \pi_1 \otimes \pi_2 \rightarrow \mathbb{C}$,
which factorizes as $B = \prod_{v} B_v$, and, for every $v$,
$B(x_{1,v}, x_{2,v}) \leq A \Sob^d(x_{1,v}) \Sob^d(x_{2,v})$;
also $|B_v| \leq 1$ when both $x_{1,v}, x_{2,v}$ are spherical of norm $1$. 

Then $B(x_1, x_2) \ll A' \Sob^{d'}(x_1) \Sob^{d'}(x_2)$ (similar notation to (S1d)). 
This follows formally: apply (S1d) first to the linear functional
$x \mapsto B(x,y)$, for $y$ fixed; we see it is bounded in absolute value by $A' \Sob^{d'}(x) \Sob^d(y)$. 
Now apply (S1d) to the functional $y \mapsto B(x,y)$, for $x$ fixed; note that
the condition on spherical inputs is satisfied not for all $v$, but only for those
$v$ where $x_v$ is spherical, this being handled by the footnote to (S1d). 
\end{Remark}

\subsubsectionind{} We prove (S2c). 
This proof is modeled on ideas from \cite{BR}.
For $x \in \bfX$, consider the linear form on $C^{\infty}_c(\bfX)$ given by:
$\ell_x:  f \mapsto \height(x)^{r} \Delta_{\adele}^r f$.  By (S2a), 
$\ell_x$ is bounded with respect
to $\Sob_{d_0}^{L^2(\bfX)}$, for some $d_0 > 0$ and uniformly in $x$ (at least after increasing $d_0$ as necessary).  It follows that,
for suitable $d > d_0$, the trace of $|\ell_x|^2$ with respect to the square of $\Sob^{\bfX}_{d}$
is finite and independent of $x$. Integrate over $x \in \bfX$ to obtain the desired conclusion.

\subsubsectionind{} \label{S3cproof} We prove (S3c); however, some preliminaries are necessary.  
We shall need to make use of an assertion substantively equivalent to the 
 well-known theorem of W. M{\"u}ller \cite{muller} that discrete automorphic spectrum is trace-class.  That fact, in turn is closely related to the fact that the Eisenstein series on a general group have finite order.   
 We believe that to handle the general case should be routine
using the techniques of M{\"u}ller's paper \cite{muller}; but we have not verified the detail, and
this is why the proof is presently only for $\GL_n$, where it is possible to give
a direct proof of M{\"u}ller's trace-class theorem -- using the work of Moeglin and Waldspurger. 

 For any automorphic representation, set $$\Sobcond(\pi) = \inf_{f \in \pi, \|f\|_2 = 1} \|\Delta_{\adele} f\|.$$
 It might be referred to as the ``Sobolev-conductor'' of the representation $\pi$; as we shall see momentarily, it is bounded above and below in terms of the analytic conductor in the case of $\GL_n$. 

We begin with:
\begin{equation}\label{cusp-bound}\mbox{The number of cuspidal $\pi$ with $\Sobcond(\pi) \leq X$
is at most polynomial in $X$}\end{equation}
\proof
There exists (by (S3a) and (S2c)) some $d$
so that the inclusion of $L^2(\bfX)$ into the Hilbert space
associated to $\Sob^{\bfX}_d$ is trace class. 
There exists $d'$ so that $\Sob^{d}_{\bfX}$ is majorized by $\Sob^{L^2(\bfX)}_{d'}$
for cuspidal functions, by (S3b). 
Choose $e_{\pi}$
in each cuspidal representation\footnote{It is not difficult to check, although it is in fact not essential, 
that the infimum which defines $\Sobcond(\pi)$ is realized.} so that $\Sob_{d'}(e_{\pi}) = \Sobcond(\pi)^{d'} \|e_{\pi}\|_2$. 
Then $\Sob_d^{\bfX}(e_{\pi}) \ll  \Sobcond(\pi)^{d'} \|e_{\pi}\|_2$. 
The trace-class feature forces  $\sum_{\pi} \Sobcond(\pi)^{-d'} < \infty$.  Thus \eqref{cusp-bound}.
\qed

With notation as in \S \ref{plancherel}, we now claim:
\begin{equation}\label{muller} \mbox{There is $d_1 \geq 1$ so that}
\int_{\chi \in \data} \Sobcond(\mcI(\chi))^{-d_1} d\chi < \infty. \end{equation}
For $\G =\GL_n$, the classification of discrete spectrum (\S \ref{Eisenstein}) 
together with the classification of discrete series (\cite{MWDS}) reduces this to \eqref{cusp-bound}.

Finally, 
 \eqref{muller} $\implies$ (S3c). Let $\ell$ be a linear form as in (S3c). We express an arbitrary $f \in L^2(\bfX)$
 as an integral $f = \int f_{\chi} d\chi$, where $f_{\chi} \in \mcI(\chi)$ (notation of \S \ref{plancherel}; 
 we have already remarked in \S \ref{pointwise} that the Plancherel decomposition is pointwise defined). 
 Using \eqref{muller} and Cauchy-Schwarz, 
 
\begin{eqnarray} |\ell(f) |^2  =  \left| \int \ell(f_{\chi}) d\chi \right|^2   && \ll  \int \Sobcond(\mcI(\chi))^d \left| \ell(f_{\chi}) \right|^2d\chi \\
 \nonumber && \leq \int \|\Delta_{\adele}^{d'} f_{\chi} \|^2 d\chi \ll \Sob_{d'}^{L^2(\quot)}(f).\end{eqnarray}

Let us note the following corollary, which can also be proved in other ways:

\begin{Corollary*}
 The number of cuspidal representations
on $\GL_n(F)$ with fixed central character and conductor $\leq X$ is bounded by a power of $X$.
\end{Corollary*}
It follows from the prior Lemma and the following Lemma due to W. T. Gan \footnote{This statement was derived from a result sketched by Gan to one the authors (A.V.), namely, the conductor of any representation of $\GL_n(k)$ -- for $k$ a nonarchimedean local field -- admitting a vector fixed under the $r$-th principal congruence subgroup is at most $r n^2$} 

\begin{Lemmat}  As $\pi$ varies through
automorphic cuspidal representations of $\GL_{n}(\Aa)$,
\begin{equation} \label{gan}\frac{ \log C(\pi)} {\log \Sobcond(\pi) } \mbox{ is bounded above and below} \end{equation} 
\end{Lemmat}

It is of interest (from the point of view of the ``Selberg-class'') to determine what the correct growth rate in this corollary is; 
among other things, it is related to the question of the number of twists
needed in converse theorems; in that connection:
\begin{Question*}
Let $N(X)$ be the number of cuspidal representations on $\GL_n(\Q)$ with analytic conductor
$\leq X$. Determine the asymptotic behavior of $N(X)$; for example, is it true that $\frac{\log N(X)}{\log X}$
approaches $n+1$?
\end{Question*}

\section{Integral representations of $L$-functions: local computations.}

\subsection{Notations}
This part, comprising \S \ref{sec:whitmodels} -- \S \ref{endpartIII},  will be essentially of local nature; we collate here some notations that will be used. 

\subsubsection{Local fields}
Throughout this part, we shall work over a local field $k$. We shall always suppose
$k$ to be given as the completion of $F$ at a place $v$; thus, we regard the place $v$
and an isomorphism $F_v \cong k$ to be given along with $k$.\footnote{The reason for this
is that it will be helpful in handling implicit constants.}

We denote by $|.|=|.|_{v}$ the absolute value on $k$ normalized so that $|x|_{\Rr}=\max\{x,-x\}$, $|z|_{\Cc}=z.\ov z$ and for $v$ non-archimedean, 
$|\varpi|=q^{-1}$ for $\varpi$ an uniformizer of $k$ with $q$ the size of the residue field; in that later case
we denote by $\mfo$ the maximal order of $k$. If $k$ is archimedean, we denote by $\deg(k)$ its degree over $\Rr$.

\subsubsection{Additive characters} We equip once for all each local field $k=F_{v}$ with an additive character
 $\psi$ obtained as follows: we let $e_{\Qq}:\Aa_{\Qq}\ra\Cc^\times$ be the unique additive character which coincides with
  $x\ra\exp(2\pi i x)$ on $\Rr$ and set $e_{F}=e_{\Qq}(\tr_{F/\Qq})$. We denote by $\psi$ the restriction of $e_{F}$ to $F_{v}=k$.
  In the sequel any local considerations involving an additive character of $k$ will refer to that character $\psi$ or its complex conjugate. If $k$ is nonarchimedean, we denote its conductor by $d_{\psi}$. We set $d_{\psi}=0$ if $k$ is archimedean.

\subsubsection{Implicit constants} \label{implicitconstants}
Unless otherwise specified, given some parameter $\eps$ say, the implicit constant in $\ll_{\eps}$ will depend (of course) on the parameter $\eps$ but may also depend on the quantity $q^{d_{\psi}}$ {\em  as $k$ varies through completions of $F$.} Note that $q^{d_{\psi}}=1$ for almost every $k$.

\subsubsection{Subgroups}\label{notationsubgroups}We denote by $\rmZ$, $\rmN$, $\rmB$ the usual upper triangular unipotent (resp. Borel) subgroup of $\GL_2$ or $\PGL_2$ and by $\rmA$ the diagonal subgroup with lower diagonal entry equal to $1$. We denote $G = \GL_2(k)$, $\bar{G}$ the quotient of $G$ by its center and $Z,N,B,A$ the group of $k$ points of $\rmZ,\rmN,\rmB,\rmA$. For $t,x,y$ in any ring $R$, we set $$z(t)=\left(\begin{array}{cc}t & 0 \\0 & t\end{array}\right),\ n(x)=\left(\begin{array}{cc}1 & x \\0 & 1\end{array}\right),\ a(y)=\left(\begin{array}{cc}y	 & 0 \\0 & 1\end{array}\right).$$ 
We let $K$ be the standard maximal compact subgroup of $G$, i.e.
the group of integral matrices if $k$ is nonarchimedean, and the stabilizer
of the ``standard'' orthogonal form $(x_1, x_2) \mapsto x_1^2 +x_2^2$ or the ``standard''
Hermitian form $(z_1, z_2) \mapsto z_1\ov{z_{1}} + z_2\ov{z_{2}}$ if $v$ is real or complex. 

For $v$ non-archimedean and $m\geq 0$ we denote by $K[m]\subset K_{0}[m]\subset K$ the subgroups of matrices 
$\left(\begin{array}{cc}a & b \\c & d\end{array}\right) \in \GL_2(\mfo)$ such that
$a\equiv d\equiv 1,\ b\equiv c\equiv 0\ (\varpi^m)$, or $c\equiv 0\ (\varpi^m)$, or with no constraint on $a,b,c,d$,  respectively.

\subsubsection{Measure normalizations} \label{measures-local}
Fix an additive haar measure on $k$ which is self-dual with respect to the
character $\psi$.  (Note that for $k \cong \mathbb{R}$, this gives Lebesgue measure;
for $k \cong \mathbb{C}$, this gives the measure $i dz \wedge d\bar{z} = 2 dx \wedge dy$;
and for $k$ nonarchimedean, it assigns mass $q^{-d_{\psi}/2}$ to the maximal compact subring of $k$.) 

Transport it to the haar measure $dn$ on $N$ via $x \mapsto n(x)$;
we equip $A = \{a(y): y \in k^{\times}\}$ and $Z=\{z(u): u \in k^{\times}\}$ with their Haar measures $\zeta_k(1) \frac{dy}{|y|}$, $\zeta_k(1) \frac{du}{ |u|}$ respectively. 

The product measure $\zeta_k(1) \frac{du}{|u|}dx \zeta_k(1)  \frac{dy}{|y|}$ defines the right Haar measure $d_{R}b$ on $B= \{z(u) n(x) a(y): u, y \in k^{\times},  x\in k\}$ while the left Haar measure $d_{L}b$ is $|y|^{-1}d_{R}b$. The left and right Haar measures on the Borel of $\ov G$ are defined explicitely through the same parameterisations (with $\zeta_{k}(1){du}/{|u|}$ removed). We denote by $dk$ the Haar probability measure on $K$.
The Haar measures on $G,\ \ov G$ are then given by $d_{L}bdk$.

The measure on $G$ determines a Plancherel measure $d\mu_P$ on the unitary dual of $G$. 
This is normalized by the formula:
$$f(1) = \int_{\pi \in \hat{G}} \langle f, \chi_{\pi} \rangle_{L^2(G)} d\mu_P, f \in C_c(G)$$
where $\chi_{\pi}$ is the character of $\pi \in \hat{G}_{Aut}$. 

For $k$  archimedean, we fix a norm, $\|.\|_{\mfg}$ say,  on $\mfg$ the Lie algebra of $G$.

\subsubsection{A measure computation} \label{meascomp1} It will be of use, as a check, to later compare $dg$ to the measure $d'g := d_L b dn $
pushed forward from the map $B \times N \rightarrow G$, $(n,b) \mapsto n w b$ ( here $w=\left(\begin{array}{cc}0 & -1 \\1 & 0\end{array}\right)$ denote the Weyl element; recall that the image of that map is dense in $G$ by the Bruhat decomposition). The later
coincides with the measure pushed forward from $dn \times d_R b$ on  $N \times B$ via $(n,b) \mapsto  nw b$,
for both are invariant and have the same behavior near the element $w$.  In fact,
\begin{equation} \label{NtoK1} d'g = \knconst dg, \ \ \knconst = q^{-d_{\psi}/2} \frac{ \zeta_k(1)}{\zeta_k(2)} .\end{equation}
Indeed, let $f \in C_c(G)$, and let $\bar{f}(g) = \int_{b \in B} f(bg) d_L(b)$. 
Then $$\int f(g) dg = \int \bar{f}(k) dk,\hbox{ while }\int f(g) d'g = \int \bar{f}(w n(x))dx.$$
The function $\bar{f}$ has the transformation property  $\bar{f}(bg) = |b|^{-1}\bar{ f}(g)$ where $|b|=|n.a(y)|:=|y|$; 
so we are reduced to computing the ratio $\frac{ \int \bar{f}(k)}{\int \bar{f}(w n(x)) dx}$
for a {\em single} $\bar{f}$ with this property.  We choose the square of the ``height'' function
$\bar{f}(g) = \frac{|\det g|}{ \| (0,1 ) g\|^{2}}$, where $\|\cdot\|$ is a $K$-invariant norm on $k^2$, 
which we normalize so that $\|(0,1)\| = 1$. We are reduced, thus, to evaluating $ \knconst = \int_{x \in k} 
\|(1, x)\|^{-2} dx$. 

\begin{enumerate}
\item 
If $k$ is nonarchimedean, with
residue field of size $q$, 
$$\knconst = q^{-d_{\psi}/2} (1 + (q-1) (q^{-2} + q^{-3} + \dots) )  = q^{-d_{\psi}/2}
\frac{\zeta_k(1)}{\zeta_k(2)}.$$
\item If $k = \mathbb{R}$, then $\knconst  = \int \frac{dx}{1+x^2}  =  \pi = \frac{\Gamma_{\R}(1)}{\Gamma_{\R}(2)}$. 
\item If $k = \mathbb{C}$, then $\knconst  = 2 \int \frac{dx \wedge dy}{(1+x^2+y^2)^2} =
4 \pi \int_{r =0}^{\infty} \frac{r \ dr}{(1+r^2)^2} = 2 \pi = \frac{\Gamma_{\C}(1)}{\Gamma_{\C}(2)}$. 
\end{enumerate}

\subsubsection{Representations: convention}\label{convention}
 In the sequel, the irreducible admissible representations of $G$ that will occur are assumed to be 
local constituents of unitary automorphic representations of $\GL_{2}(\Aa)$. In particular, by the work of Kim and Shahidi \cite{KSh}, such representations will never be complementary series with parameter $\geq\theta$.  (cf. \refs{localglobalbounds}).

\subsubsection{The local analytic conductor} \label{sec-localanalyticconductor}Given $d=1,2$ and $\pi$ an irreducible admissible representation of $\GL_{d}(k)$ acting on a vector space $V$ say we define the (local) analytic conductor of $\pi$, $C(\pi)$,  to be
\begin{itemize}
\item[-] $k$ non-archimedean: $C(\pi)=q^{f(\pi)}$ where $f(\pi)$ is the conductor of $\pi$.
\item[-] $k$ archimedean: the $L$-factor of $\pi$ has the form
$$L(\pi,s)=\prod_{i=1}^{d}\Gamma_{k}(s+\mu_{\pi,i})$$
with $\Gamma_{\Rr}(s)=\pi^{-s/2}\Gamma(s/2)$, $\Gamma_{\Cc}(s)=2(2\pi)^{-s}\Gamma(s)$ and  $\mu_{\pi,i}\in\Cc$; the analytic conductor is given by
$$C(\pi)=\prod_{i=1}^d(2+|\mu_{\pi,i}|)^{\deg(k)}.$$
\end{itemize}
In the sequel, the usual notation of Vinogradov $A\ll_{\pi} B$ will have the following slightly
 more precise meaning: there are constants $C,d\geq 0$ (independent of $\pi$) such that $|A|\leq C.C(\pi)^d.|B|$; we will then say that $A$ is bounded by $B$, polynomially in $\pi$.

\subsubsectionind{} \label{SSmod-1}
We recall that we have defined previously (\S \ref{SSur}) a system of Sobolev norms $\Sob^V_d$
on any unitary representation $V$ of $G$; these norms have the property, among others,
that $\Sob^V_d(v) = \|\Delta^d v\|$ for a certain positive self-adjoint operator $\Delta$ on $V$. 
As a complement to our earlier convention \ref{implicitconstants}, the indices of the Sobolev norms
 occurring in inequalities such as $A\leq \Sob^{\star}(f)B$ will always be independent of the field $k$. 

\subsubsection{Principal series representations}\label{ps-local} For $(\omega,\omega')$ a pair of characters of $k^\times$, we denote by $\mcI(\omega,\omega')$ or
 $\omega\boxplus\omega'$ the corresponding {\em principal series representation} of $G$ which is unitarily induced from the corresponding representation of $B$: the $L^2$-space of functions $f$ on $G$ such that
$$f(\left(\begin{array}{cc}a & b \\0 & d\end{array}\right)g)=|a/d|^{1/2}\omega(a)\omega'(d)f(g),$$
with respect to the inner product $$\peter{f,f}=\int_{K} |f(k)|^2 dk.$$

If $\omega,\omega'$ are unitary, the resulting inner product is  $G$-invariant
and $\mcI(\omega, \omega')$ is thus a unitary representation. Such representation will be called a {\em unitary principal series}. In this case,  for $f \in \mcI$, 
\begin{equation} \label{ntok} \|f\|^2 = \int_{k \in K} |f(k)|^2 =\knconst^{-1} \int_{x \in k} |f(w n(x))|^2 dx,\end{equation}
where $\knconst = \knconstant$ is as in \S \ref{meascomp1}.

Even if $\omega, \omega'$ are not unitary, 
a principal series representation $\mcI(\omega,\omega')$ may be {\em unitarizable}
 although not, in general, with respect to the above inner product. In particular, in such cases, $\ov\omega^{-1}={\omega'}$.%

\subsubsection{Deformation} \label{sss:deform} 
It will be of utility to deform induced representations. To that end we introduce the following notation:
Suppose that $\pi$ is induced from\footnote{The choice of $\omega,\omega$ is unique up to order.} two unitary characters $\omega, \omega'$, i.e.
$\pi = \mcI(\omega, \omega')$ and
and $s$ is a complex number. 
We set $$\pi_{s}:=\omega|.|^s\boxplus\omega'|.|^{-s};$$ it is realized in the induced space
$\mcI_s$. For $f\in\mcI$, we define $f_{s}$ to be the vector in $\mcI_{s}$  whose restriction to the maximal compact $K$ coincide with that of $f$.

The map $f \mapsto f_s$ is $K$-equivariant from $\pi$ to $\pi_s$. It is not $G$-equivariant in general. 

\subsubsection{$\gamma$-factors} \label{localgamma}
If $\pi$ is a unitary representation of $\GL_r(k)$, we define the $\gamma$-factor, as usual,
as the ratio 
 $$\gamma(\pi, \psi,s) =
\eps(\pi,\psi, s) \frac{ L(\tilde\pi ,1-s)}{L(\pi, s)}.$$

We recall that, if $\psi$ is unramified and $k$ non-archimedean, 
then
 \be\label{epsilonfactorbound}\eps(\pi, \psi, s) = q^{-f_{\pi} (s-1/2)},\ee
  where $f_{\pi}$ is the local conductor
of $\pi$ (in the sense of \cite{JPSS}).   
When $r=1$ and $\pi$ is a character $\chi$ of $k^{\times}$, we obtain the $\gamma$-factor of Hecke and Tate. Recall that, for $\Phi$ a Schwarz function on $k$
and $$Z(\Phi, \chi,s) = \int \Phi(x) \chi(x)|x|^s d^{\times}x$$
 we have
  \be\label{Tate}Z(\widehat{\Phi},\chi^{-1},1-s) = 
\gamma(\chi, \psi,s) Z(\Phi, \chi,s),\ee where $\Phi \mapsto \widehat{\Phi}$ is the Fourier transform;
 also $\gamma(\chi, \psi,s)$ can be expressed formally as $\int_{y \in k} \psi(x) \chi(x)^{-1}|x|^{-s} dx$,
and it has poles only when $\Re(\chi\alpha^s) = 1$.

We shall use the following estimtes for the $\gamma$-factor: given $0<\epsilon<1/100$, one has for $GL_{1}$,
\begin{equation} \label{gammaboundchi}  | \gamma(\chi,  \psi, \sigma)|  \asymp_{\epsilon, q^{d_{\psi}}} \Cond(\chi)^{-(\sigma-1/2)}, 
\ \ (\sigma \in [\epsilon, 1/2-\epsilon] ),  \end{equation}
and for $\GL_2$,
\begin{equation} \label{gammabound}  | \gamma(\pi,  \psi, \sigma)|  \asymp_{\epsilon,q^{d_{\psi}}} \Cond(\pi)^{-(\sigma-1/2)}, 
\ \ (\sigma \in [3/26+\epsilon, 1/2-\epsilon] )  \end{equation}
This follows from the prior remark, the dependence of $\eps$-factors on $\psi$,
and the fact that $L(\pi, s)$ and $L(\tilde \pi, 1-s)$ are uniformly bounded above and below in the stated 
region, if $k$ is nonarchimedean; and by a direct computation (see \cite[Chap. 1, \S 5 \& 6]{JL}) when $k$ is archimedean. 

For $\GL_1$, the same bound remains valid (with the same proof), so long
as $\pi$ is unitary; we may even replace $(3/26, 1/2)$ by any compact subinterval of $[0,1]$.

\subsubsectionind{}\label{gausssumcomputation}

Let $f \in C^{\infty}_c(\mathbb{R}_+)$, and consider
the integral $$G(\lambda, t) := \int_0^{\infty} f(x) |x|^{i \lambda} e^{i t x} dx.$$
The method of stationary phase shows that
the dominant contribution to this integral comes from stationary
points of the phase function $\lambda \log(x) + t x$, i.e., when $x = \lambda/t$.
If $f$ is fixed, it therefore follows that $G(\lambda, t)$ is small unless
$\lambda$ and $t$ are of the same order of magnitude. Moreover, 
stationary phase predicts that, when $\lambda \sim t$, the magnitude
$|G(\lambda, t)| \sim \lambda^{-1/2}$.

A similar phenomenon occurs when we replace $\R$ by a nonarchimedean local field
with residue field of size $q$; 
in that case, the pertinent integrals are closely related to Gauss sums. 
If $\chi$ is a character of conductor $C=q^{f_{\chi} } >1$, and $D=q^{d_{\psi}}$ the discriminant of $\psi$,  then: for $t\in k^\times$,
\be\label{gausscomp}\int_{|x| = |u|} \chi(x) \psi(tx)  d^{\times} x =
 \begin{cases}0&, |u| \neq CD/|t| \\
 \eta C^{-1/2} |u|^{\Re\chi}&, |u| = CD/|t|
 \end{cases}
 \ee
 with $|\eta|=1$.  
 
The following general version will be helpful:

\begin{Lemmat} \label{almostgauss} Let $A < B \in \mathbb{R}$.
Let $k$ be a local field, $\phi \in C^{\infty}_c(k^{\times})$ a smooth function, and $\chi$
a character of $k$ of conductor $C$ such that, writing $|\chi(x)|=|x|^{\sigma}$,
one has $A \leq \sigma \leq B$.
 
Set, for $t \in k^\times$, 
$$G_{\phi}(\chi,t) = \int_{k} \phi(x)\psi(tx) \chi(x)  dx.$$
Then:
\begin{enumerate}
\item For every $N \geq 0, \epsilon > 0$, we have:
$$ | G_{\phi}(\chi,t)| \ll_{\phi,N, A, B, \epsilon } C^{-1/2+\epsilon}  \max\left(\frac{1+|t|}{C}, \frac{C}{|t|}\right)^{N};$$
\item Fix $\eps > 0$. Then there exists $t$, satisfying $\left| \frac{\log|t|}{\log (C+1)} - 1 \right| \leq \eps$, 
so that $|G_{\phi}(\chi,t)| \gg_{\varphi,\eps} C^{-1/2-2\eps}$. 
\end{enumerate}
\end{Lemmat} 
\proof 
 The statement (2) follows from statement (1) and Plancherel formula:
$$\int_{k}|G_{\phi}(\chi,t)|^2dt=\int_{k}|\phi(x)\chi(x)|^2dx,$$ noting
that the volume of the set $|t| \leq C$ is $\asymp C$.
It suffices thereby to prove (1).

We shall give the proof of the bound  $| G_{\phi}(\chi,t)| \ll_{\phi,N, A,B}   ((1+|t|)/C)^{N}$.
Let us suppose, first of all, that $k$ is archimedean, and fix
a basis $\{D_i\}$ of $k^{\times}$-invariant differential operators on $k^{\times}$
of degree $[k:\R]$.  Then $\chi$ is an eigenfunction of each $D_i$,
and there exists $i$ so that the eigenvalue $\lambda_i$ satisfies $|\lambda_i| \gg C$. 
Now, 
$$ | G_{\phi}(\chi,t)| \ll C^{-N} \int D_i^N(\phi(x) \psi(tx))  \chi(x) \ll_{\phi,N,A,B} C^{-N} (1+|t|)^N$$
If $k$ is nonarchimedean, there exists a constant $c > 0$, depending only on $\phi$
and the discriminant of $\psi$, 
so that $\phi(x) \psi(tx)$ is invariant under $x \mapsto xy$ whenever $|t(y-1)| \leq c$. 
In particular, $G_{\phi}(\chi,t)$ vanishes whenever $|t| C^{-1} < c$. 

The bound $| G_{\phi}(\chi,t)| \ll_{\phi,N, A,B}   (|t|/C)^{-N}$ is similar; we replace
the role of multiplicative translations and multiplicatively invariant differential operators
by additive translations and additively invariant differential operators. We also reverse the roles of $\chi$ and $\psi$ in the prior argument.

As for the bound $G_{\phi}(\chi,t)\ll_{\phi,N,A,B} C^{-1/2+\eps},$
we may assume that $\chi$ is unitary, after replacing $\phi$
by $\phi |x|^{\sigma}$ suitable $\sigma$. 
Fix $\eps>0$ small. By \refs{Tate}, applied to $x\mapsto\phi(x)|x|^{1-\eps}\psi(tx)$, we may 
write $G_{\phi}(\chi,t)$ as
$$\frac{ \int_{k}\widehat{(\phi\alpha^{1-\eps})}(x+t)(\chi\alpha^\eps)^{-1}(x) dx}{\gamma(\chi,\psi,\eps)}=
\frac{\int_{k}\widehat{(\phi\alpha^{1-\eps})}(x)\chi^{-1}(x-t)  |x-t|^{-\eps} dx}{\gamma(\chi,\psi,\eps)}.$$
Our assertion \label{Desid10} follows from
the rapid decay of $\widehat{(\phi\alpha^{1-\eps})}(x)$ and \refs{gammaboundchi}.
 \qed

\subsection{Whittaker models.} \label{sec:whitmodels}
Let $\psi$ be the non-trivial additive character of $k$ described previously. Given
 $\pi$ a generic representation of $G$, $\pi$
admits a unique realization in a Whittaker model, which we shall denote by $\Whit(\pi,\psi)$, and thereby in a Kirillov model, which we denote by $\kirill(\pi,\psi)$. 
We equip the Kirillov model
with that inner product given by the $L^2$-inner product on $k^{\times}$: 
\be\label{whitinner}\peter{W,W}=\int_{k^\times} |W(y)|^2\dti y.
\ee
Unless otherwise specified, the inner product \eqref{whitinner} is {\em always} the one that
we will put on any Kirillov or Whittaker model. 

Let us suppose, moreover, that $\pi$ is unitary, i.e. {\em equipped with a $G$-invariant inner product.}
Then we may choose an intertwiner $\pi \mapsto \kirill(\pi, \psi)$ which 
preserves inner products. Such an intertwiner is unique up to a scalar of absolute value $1$,
and we will regard it as having been fixed; thus we will identify $\pi$ with $\kirill(\pi, \psi)$. 
 Throughout this paper, since we are concerned only with the size
of various quantities, the scalar of absolute value $1$ will not prove a problem. 
Therefore -- in the setting where $\pi$ is equipped with an inner product -- we
will regard the following statement, for example, as meaningful:
\begin{equation} \mbox{ Choose $v \in \pi$ corresponding to the element $W \in \Whit(\pi, \psi)$ \dots } \end{equation}

\subsubsection{Intertwiners from the induced model to the Whittaker model}
It will be useful later to record these in somewhat explicit form. 
Notation as in \S \ref{ps-local}, let $\pi = \mathcal{I}(\omega, \omega')$.
An explicit intertwiner $\pi \rightarrow \mathcal{K}(\pi, \psi)$ is given by:
\begin{equation}
\label{principaltokirillov} f \mapsto W_f,  \ \
W_{f}(g)= \whitconst^{-1/2} \int_{k}f(wn(x)g)\psi(x)dx, \whitconst = \zeta_k(1) \knconst 
\end{equation}
 If $\omega, \omega'$ are unitary, then the map \refs{principaltokirillov} is isometric by \eqref{ntok}; this was the reason for the inclusion of $\whitconst$. 

Taking $g=a(y)$, one has
$$W_f(a(y)) = \whitconst^{-1/2} |y|^{1/2} \omega'(y) \int_{x \in k} f( w n(x) )  \psi(xy) dx.$$
In general, the right-hand integral is not absolutely convergent, but it is conditionally convergent --
interpreted as a limit of integrals over increasing
compacta $|x| \leq R$  -- 
if $f$ is a smooth vector. 

It follows from the above expressions, using the shorthand $\alpha (y) := |y|$, that:
\begin{equation} \label{wfmellin} \whitconst^{1/2} \int W_f(y) \chi(y) d^{\times}y =  \left(  \int f(w n(x)) \chi'(x)^{-1} dx  \right) \gamma(\chi'^{-1},\psi,1) \zeta_k(1), \end{equation} 
where $\chi' := \chi \omega' \alpha^{1/2}$.

\begin{Remark} \label{convergence-remark} 
It is useful to discuss convergence. Suppose $\omega, \omega'$ unitary.  The left hand side of \refs{wfmellin}
is holomorphic for $\Re(\chi) > -1/2$. On the other hand, the right-hand side integral
is absolutely convergent and holomorphic in the range $-1/2 < \Re(\chi) < 1/2$. 

We observe that the left-hand side of \eqref{wfmellin}
has poles precisely when  $\chi \omega = \alpha^{-1/2}$ or $\chi \omega' =\alpha^{-1/2}$.
These poles arise from the asymptotics of $W(y)$ as $|y| \rightarrow 0$.  On the right-hand side,
the integral involving $f$ has poles at $\chi \omega' = \alpha^{1/2}$, arising from the 
behavior as $|x| \rightarrow 0$; this pole is cancelled by the zero
of the $\gamma$-factor $\gamma(\chi^{-1} \omega'^{-1}, \psi,1/2)$ at this point. 
Moreover, one has poles at $\chi \omega =\alpha^{-1/2}$ (arising
from the asymptotics of $f$ as $|x| \rightarrow \infty$) and $\chi\omega' =\alpha^{-1/2}$ (arising from the pole of $\gamma$). 
\end{Remark}

\subsubsection{The Jacquet-Langlands functional equation} 
Let $\pi$ be generic, not necessarily unitary. 
For $\chi$ a character of $F^\times$, let $\Re \chi$ denote the real part of $\chi$ (
$|\chi(.)|=|.|^{\Re \chi}$) and set for $W \in \kirill(\pi, \psi)$, 
\begin{equation} \label{lchidef} \ell^{\chi}(W):={\int_{k^\times}W(u)\chi(u)d^\times u};\end{equation}
This integral is absolutely convergent for $\Re(\chi)=\sigma > 0$.
By the theory of Hecke and Tate, $\ell^{\chi}(W)$ admits a meromorphic continuation to $\widehat{k^\times}$; more precisely
 the ratio $\frac{\ell^{\chi}(W)}{L(\pi\otimes\chi,1/2)}$ is holomorphic; it satisfies the local functional equation
\begin{equation}\label{localfcteqn}\gamma(\pi\otimes\chi,\psi,1/2)\ell^{\chi}(W)=\ell^{\omega^{-1} \chi^{-1}}(\widetilde W)
\end{equation}
where $\widetilde W=\pi(w).W$, and $\omega$ is the central character of $\pi$.

Note we regard $\widetilde W$ as belonging to the Whittaker model of $\pi$; 
it would also be possible to incorporate the $\omega$-twist into it and regard
it as being in the Whittaker model of $\tilde{\pi}$, bearing in mind that
$\tilde{\pi} \cong \omega(\det)^{-1} \otimes \pi$. In any case, the Sobolev norms
of $\widetilde W$ and $W$ are the same, because the element $w$ belongs to $K$.

\begin{Prop}\label{zeroboundprop} (Pointwise bounds for the Whittaker function.) For any $N\geq 1$, any $\eps>0$ and all $W \in \mathcal{K}$, there exists $d=d(N)$ so that:
$$| W(y)|\ll_{N} \Sob_d^{\kirill(\pi,\psi)}(W) \begin{cases} |y|^{1/2-\theta }, &  |y| \leq 1, \\  |y|^{-N}, & |y| \geq 1 \end{cases}.$$
\end{Prop}

\proof
The assertion for $|y| \geq 1$ is elementary; it is a quantification
of the well-known fact that Whittaker functions decay at $\infty$.  For example, in the real case,
one may use the fact that there exists an element of the Lie algebra of $G$ which
acts on the Kirillov model by $W(y) \mapsto y W(y)$. We leave it to the reader (see also \cite[(29)]{BH2}). 
As a consequence for the assertion with $|y| > 1$, 
we deduce that 
\be\label{Re>0}|\ell^{\chi}(W)| \ll S^{\kirill(\pi,\psi)}(W)\ee whenever $\Re(\chi) > 0$.

By  inverse Mellin transform\footnote{The measure on the space of $\chi$ is obtained as follows:
the $\chi$-space is a principal homogeneous space for the dual group to $k^{\times}$, and
we transport the dual measure to that space.} and \refs{localfcteqn} , we have
\begin{align*} W(y)&=\int_{\Re \chi=\sigma\gg 1} \chi^{-1}(y) d\chi\int_{k^\times}W(u)\chi(u)d^\times u\\\
&=\int_{\Re \chi=\sigma\gg 1}\chi^{-1}(y)  
\frac{\ell^{\chi^{-1} \omega^{-1}}(\widetilde W)}{\gamma(\pi\otimes\chi,\psi,1/2)}.d\chi
\\ &= \int_{\Re \chi=\sigma'}\chi^{-1}(y)\frac{\ell^{\chi^{-1} \omega^{-1}}(\widetilde W)}{\gamma(\pi\otimes\chi,\psi,1/2)}.d\chi.
\end{align*}
At the last stage, we have shifted the contour to $\Re\chi=\sigma'$ for $\sigma'\in (\thetaval-1/2,0)$: under our convention, the contour shift does not cross any pole of the local $L$-factor $L(\pi \times \chi,1/2$), and therefore the integrand remains holomorphic.

We observe that in the nonarchimedean case $$\int_{k^\times}\widetilde W(u)\chi^{-1}(u) \omega^{-1}(u) d^\times u=0$$ as long as the conductor of $\chi \omega$ is strictly greater than any $r$ for which $\widetilde W$ is $e[r]$-invariant. 
Applying then (S4d),  \refs{gammabound}, the previous bound \refs{Re>0} in the range $\Re(\chi^{-1})  >0$ and the previous observations, we conclude that, for any $m \geq 1$,

$$|W(y)| \ll_{m, \sigma'} |y|^{-\sigma'} \Sob^{\kirill(\pi,\psi)}(W) \int_{\chi} C(\chi \otimes \omega)^{-m} C(\pi \otimes \chi)^{\sigma'}.$$
The inner integral may be expressed as $\int_{\chi} C(\chi)^{-m} C(\tilde{\pi} \otimes \chi)^{\sigma'}$. 
Since $\sigma'<0$ the inner integral is bounded independently of $\pi$ by taking $m$ large enough. 
(This is only for convenience; an implicit dependence on $C(\pi)$ could be absorbed
into the Sobolev norm using the results of \S \ref{S3cproof}.) In the archimedean case, we reason similarly, replacing the
statement ``$\dots = 0$ as long as the conductor of $\chi$ is strictly greater than ... '' by integration by parts\footnote{
 in the archimedean case, there is an element $Z\in\mfg$  (more precisely in the Lie algebra of $a(k^\times)$) of bounded norm so that
$\int_{k^\times}\widetilde W(u)\chi^{-1}(u)d^\times u\ll C(\chi)^{-1}\int_{k^\times}Z.\widetilde W(u)\chi^{-1}(u)d^\times u.$}.

\qed 

We shall need later the following variant also.   Its role in our proof will be the following:
we shall need, at certain points, to utilize non-unitarizable representations, and
this Proposition that follows allows us to maintain control on their Whittaker functions
so long as they are ``close enough'' to being unitary. 

\begin{Prop} \label{zeroboundprop-def} Notation as in \S \ref{sss:deform}, so that $\pi$ is a unitary principal series, $f=f_0 \in \pi=\pi_0$, and $f_s \in \pi_s$ a deformation as in \S \ref{sss:deform}. Let
$W_s$ be the image of $f_s$ under the intertwiner \eqref{principaltokirillov}. Then 
there exists an absolute $d \geq 0$ so that, for $0\leq  \delta \leq 1/20$, and $|\Re(s)|\leq \delta$,
\begin{equation} \label{desideratum1} \int_{k^\times} | W_s (a(y)) |^2  \max(|y|, |y|^{-1})^{\delta} d^{\times}y \ll \Sob_{d \delta}^{\pi} (f)^2 , \end{equation} 
\end{Prop}

\proof
The Mellin transform $\int_{k^\times} W_{0}(x) \chi(x) d^{\times} x$
is  absolutely convergent for $ \Re(\chi)>-1/2$, cf. Remark \ref{convergence-remark}.  

Fix any $\eta \in (0,1/2)$; we will use it to parameterize a compact subinterval of $[0,1]$. 
Take
any $\beta \in (-1+\eta, -\eta)$; by Mellin inversion and \eqref{wfmellin}, 

\begin{equation} \label{mellininv-deriv} f_0(w n(x)) =  \frac{ \whitconst^{1/2}}{\zeta_k(1)} \int_{\Re(\chi) = \beta}  \chi(x) \frac{\ell^{\chi  \omega'^{-1} \alpha^{1/2}}(W_{0})}{ \gamma(\chi^{-1},\psi,0)} \ d\chi,\end{equation} 
with the measure $d\chi$ is as before.
We now bound $f_{0}(w n(x))$: 
we apply  Proposition \ref{zeroboundprop}  and \eqref{gammabound} to equation \eqref{wfmellin}, 
and use similar reasoning to Proposition \ref{zeroboundprop} to pass from
a pointwise bound for the integrand in \eqref{mellininv-deriv}, to the integral.

We arrive at the existence of $d$ so that for $x \in k$ and $ \beta \in [-1 + \eta, - \eta]$: 
\begin{equation} \label{f0bound} |f_0(w n(x))| \ll_{\eta} \Sob^{\pi}_{d}(f_0) |x|^{\beta}.\end{equation}

 Observe that:
$f_s( w n(x)) = f_0( w n(x)) q(x)^{s},$
where $q(x) = \|(1,x)\|^{-2}$ (notations of \S \ref{meascomp1}). 
Combining this with \eqref{f0bound} and substituting into \eqref{wfmellin}, we see that for $\eta\in(0,1/2)$ there exists $d=d(\eta) \geq 0$ so that, whenever 
$|\Re(\chi)| + |2\Re(s)| \leq 1/2-\eta$,
$$ \left| \int W_s(y) \chi(y) d^{\times}y  \right| \ll \Sob_d^{\pi}(f_0)  C(\chi)^d.$$

Now fix $\beta=1/20$.  Estimate $\int |W_s(y)|^2|y|^{\pm\beta} d^{\times}y$ by Plancherel,
yet again controlling the $\chi$-integration as in Proposition \ref{zeroboundprop}.  Thus, there exists $d$ so \be\label{betabound} \int |W_s(y)|^2 \max(|y|, |y|^{-1})^{\beta}  d^{\times}y \ll \Sob_{d}^{\pi}(f_0)^2 \ \
(|\Re(s)| \leq \beta) \ee

This is almost what we want; however, the ``number of derivatives'' in the Sobolev norm
in \eqref{desideratum1} is ``small,'' whereas \eqref{betabound} involves an unspecified number of derivatives (i.e., we have no control of $d$). To get around this, we interpolate. 

Define for every complex $s$ with $\Re(s) \in [0,\beta]$ and $\tau\in[-1,1]$, a function $F_s \in L^2(k^{\times})$ via the rule
$$F_s(y) = W_{\tau s}(\Delta^{-\frac{d}{\beta}s} f_0) \cdot \max(|y|, |y|^{-1})^{s/2}$$
Here $\Delta$ denote the local laplacian discussed in \S \ref{sec:laplacelocal} and $W_{\tau s}(\Delta^{-\frac{d}{\beta}s} f_{0})$ denotes the result of applying the intertwiner
\eqref{principaltokirillov} to the representation $\pi_{\tau s}$ and the vector
whose restriction to $K$ is that of $\Delta^{-\frac{d}{\beta}s} f_0$.  Also, $d$
is as in \eqref{betabound}. 

Then $s \mapsto F_s$ is a holomorphic function from $\{s \in \C: |\Re(s)| \leq \beta\}$
to $L^2(k^{\times})$. Moreover, independently of $\tau$, $\|F_s\|^2 \ll \|f_{0}\|^2$ for $\Re(s) = 0$, while for $\Re(s) = \beta$, we have by  \refs{betabound}
$$\|F_s\|^2 \ll \Sob_{d}^\pi(\Delta^{-\frac{s}{\beta}d}f_{0})^2=\|\Delta^{(1-\frac{s}{\beta})d}f_{0}\|^2 = \|f_{0}\|^2;$$ 
hence we have $\|F_{s_{0}}\|^2 \ll \|f_0\|^2$ in the region $\Re(s_{0}) \in [0,\beta]$ (independently of $\tau$) : 
for every $s_0$, apply the usual maximum principle to the holomorphic function $s \mapsto \langle F_s, F_{s_0} \rangle$. We have shown that for all $|\Re(s)| = \delta$, all $\tau \in [-1,1]$, 
$$\int |W_{\tau s}(f)|^2 \max(|y|,|y|^{-1})^{\delta}  \ll \Sob^{d \delta/\beta}(f_0)^2,$$
which implies the desired assertion. 
\qed

\subsubsection{Computations in the Kirillov model}\label{seckirillov}

Suppose that $k$ is archi\-medean; the action of the Lie algebra of $G$, $\mfg$, on $\kirill(\pi,\psi)$ is well known \cite[Chap. 1 \S 5,6]{JL}.
For instance, it follows from the explicit computations of \cite[\S 8.1.1]{Ve} that
given $Z_1, \dots, Z_m \in \mathfrak{g}$, a smooth compactly supported function 
$W\in C^\infty_{c}(k^\times)\subset\kirill(\pi,\psi)$ one has
 $$\|Z_1 \cdot Z_2 \cdot \dots Z_m \cdot W\|\ll C(\pi)^{m/\deg(k)}\Sob_{2m}^{k^\times}(W),$$
 where $$\Sob_{2m}^{k^{\times}}(W) := \int_{k^{\times}} (|y|+|y|^{-1})^{2m}  |\sum_{i} \mathcal{D}_i W |^2 d^{\times}y,$$
and $\mathcal{D}_i$ ranges over a basis of $k^{\times}$-invariant differential operators
of degree $\leq 2 m$. 

\begin{Rem} This result is stated in \cite[Lemma 8.4]{Ve}, except
the exponent of $C(\pi)$ is given as $2m/\deg(k)$ rather than
$m/\deg(k)$. 
however, the improved statement follows from a simple computation: 
if $\lambda_{\pi}$ denote the Casimir
eigenvalue, $\nu_{\pi}$ the eigenvalue for the action of the Lie algebra of diagonal matrices and $D := y \frac{d}{dy}$, then there are generators of the Lie algebra (independent of $\pi$) which act respectively as $ \nu_{\pi} +D, y , y^{-1} (\lambda_{\pi} +  ( \nu_{\pi} + D) - (\nu_{\pi} + D)^2)$. Taking
into account that both $\nu_{\pi}$ and $\lambda_{\pi} - \nu_{\pi}^2$ are bounded by $C(\pi)$, we are done. 
\end{Rem}
By the mean-value Theorem, we have, for $\delta \in \Rr$, $0<\delta< 1$:
\begin{eqnarray} \nonumber
 \| \exp(\delta Z).W - W\|_{\pi} &&\ll  \delta\sup_{\delta' \in [0,\delta]} \| \exp(\delta' Z) \cdot Z \cdot W \| = \delta \|Z.W\|\\
&&  \ll \delta C(\pi)^{1/\deg(k)}\Sob_{2}^{k^\times}(W).\label{meanvalue}
\end{eqnarray}

Moreover, all the implicit constants in the above inequalities may be taken to depend continuously on
$Z_i$ or on $Z$.
In particular, they may be taken to be uniform when $Z$ is restricted
to a fixed compact set. 

\subsubsection{Some norms related to Whittaker models} \label{Whitmodel-0}

Let $\pi_2, \pi_3$ be two generic irreducible representations of $G$ with $\pi_{3}$ isomorphic to a unitary principal series $\mathcal{I}_{3}$ (see \S \ref{ps-local}); let $\chi_{i}$, $i=2,3$
denote their central characters. Let $L^2(N \backslash G;\chi_{2}\chi_{3})$ be the $L^2$-space of
$N$-equivariant functions on $G$ which transform by $\chi_{2}\chi_{3}$ under multiplication by the center, with respect to the inner product
$$\peter{\varphi,\varphi}_{N\bash G}=\int_{N\bash\bar{G}}|\varphi(g)|^2dg$$

\begin{Lemmat}\label{jisometry}
The $\GL_2(k)$-equivariant map
$$J = J_{\pi_2, \pi_3}: (W, f) \in \Whit(\pi_2) \times \mathcal{I}_3 \mapsto Wf.$$
 defines an isometry from $\pi_2 \otimes \pi_3$ into $L^2(N \backslash G;\chi_{2}\chi_{3})$.
\end{Lemmat}
\proof It suffices
to check that it preserves inner products of pure tensors.  Take $(W,f), (W', f')$, and observe that
-- using the Iwasawa decomposition -- 
\begin{equation} \label{miracle-factor} \langle J(W,f), J(W', f') \rangle_{N \backslash G} = \int_{y \in k^{\times}} \int_{k \in K}
W(a(y) k) \overline{W'(a(y) k)} f(k) \overline{f'(k)} d^{\times}y \end{equation}
Since the integral $\int_{y \in k^{\times}} W(a(y) k) \overline{W'(a(y) k)}$ equals the inner product
$\langle W, W' \rangle$ {\em for any $k \in K$}, the right-hand side
factors as a product $\langle W, W' \rangle_{\pi_2} \langle f, f' \rangle_{\pi_3}$. We are done. 
\qed

\begin{Lemmat} \label{degboundlem} Suppose $k$ non-archimedean. Let $$\pi_2=\chi_{2}^+\boxplus\chi_{2}^-, \pi_3=\chi_{3}^+\boxplus\chi_{3}^-$$ be unramified principal series representations with $\pi_{3}$ tempered (i.e. the characters $\chi_{3}^\pm$ are unitary). We assume moreover that the additive character $\psi$
is unramified.

Let $W$ be a spherical vector, in the Whittaker model of $\pi_2$ and $f$ spherical in the induced model of $\pi_3$. 
Let $\pi_2(s), \pi_3(t)$ and $W_{s} \in \pi_2(s), f_{t} \in \pi_3(t)$ be deformations
of respectively $W, f$ parameterized by complex parameters $s$ and $t$. Here $W_{s}$ is obtained, as before,
by applying the intertwiner \eqref{principaltokirillov} to the $f_2(s)$, 
where $f_2 \in \pi_2$ is the preimage of $W$ under the intertwiner
\eqref{principaltokirillov}.

We denote the map $J_{\pi_{2},\pi_{3}}$ by $J$ and $J_{\pi_{2}(s),\pi_{3}(t)}$ by $J_{s,t}$. 
For $0\leq \delta<1/10$, $|s|, |t| \leq \delta/2$, and $u\in G$ one has
\begin{equation} \label{degbound}\frac{  \langle u.J(W \otimes f), J_{s,t}(W_{s} \otimes f_{t}) \rangle_{N\bash G}}{\|W\|^2 \|f\|^2} \ll
\|\Ad(u)\|^{\theta+ \delta -1}. \end{equation}
\end{Lemmat}

\proof We may assume that $W(1)=f(1)=1$.

We shall reduce the bound \eqref{degbound} to a bound where the representation $\pi_2,\pi_3$
are unitary, although $\pi_2$ is not tempered. 
In fact, we claim that:
\begin{equation} \label{boundwanted}
|\peter{u.J(W \otimes f), J_{s,t}(W_{s} \otimes f_{t})}_{N\bash G}| \leq 2 \peter{u.J(W^{0}_{\theta+\delta}\otimes f^{0}),J(W^{0}_{\theta+\delta}\otimes f^{0})}_{N\bash  G};\end{equation}
where $\pi_{2}$ replaced by the 
complementary series $|.|^{\theta + \delta}\boxplus|.|^{-\theta - \delta}$ with parameter $\theta + \delta$
and $\pi_3$ replaced by $1 \boxplus 1$; $W^{0}_{\theta+\delta}$, $f^{0}$ denote the spherical vectors in the corresponding models of the respective representations, normalized as above. Then by Lemma \ref{jisometry}, the latter integral factors as
\begin{multline*}\peter{u.J(W^{0}_{\theta+\delta}\otimes f^{0}),J(W^{0}_{\theta+\delta}\otimes f^{0})}_{N\bash G}=\peter{u.W^{0}_{\theta+\delta},W^{0}_{\theta+\delta}}\peter{u.f^{0},f^{0}}\\
\ll \|\Ad(u)\|^{\theta+ \delta -1}\|W^{0}_{\theta+\delta}\|^2\|f^{0}\|^2 
\ll \|\Ad(u)\|^{\theta+ \delta -1}\|W\|^2\|f\|^2
\end{multline*}
as follows from bounds for matrix coefficients and a computation of the norm in the Whittaker model.

As for the claim: we note that by the Cartan decomposition and $K$-invariance, we may assume that $u$ is of the form
$a(t_0)$ with $0<|t_0|\leq 1$; again by the Cartan decomposition
$$\peter{u.J(W \otimes f), J_{s,t}(W_{s} \otimes f_{t})}_{N\bash G}=\int_{k^\times}\int_{K}u.(Wf)(a(y)k)\ov{W_{s}f_{t}}(a(y))|y|^{-1}dk\dti y.$$
Fix a unitary character $\chi$, and 
consider the space of  functions $F$ on $G$ sa\-tisfying
$F(z(t)n(x)gk)=\chi(t)\psi(x)F(g),\ k\in K$. 
This space is equipped with a natural inner product, namely, that arising
from $L^2(N \bash \bar{G})$.  This inner product
is a Hermitian form in the values $(F(a(\varpi^\alpha)))_{\alpha\in\Zz}$. In particular,
the matrix coefficient $\langle u F_1, F_2 \rangle$ may be expressed:
\begin{multline*}
\langle u. F_1, F_2 \rangle = \int_{k^\times}\int_{K}u.F_{1}(a(y)k)\ov{F_{2}}(a(y))|y|^{-1}dk\dti y\\=\sum_{\alpha,\beta\in\Zz}c^{u}_{\alpha,\beta,\chi}F_{1}(a(\varpi^\alpha))\ov{F_{2}(a(\varpi^\beta))}.
\end{multline*}
Moreover the $c^{u}_{\alpha,\beta,\chi}$ satisfy
$|c^{u}_{\alpha,\beta,\chi}|\leq c^{u}_{\alpha,\beta,1}.$ It follows that if
$\bar{F}_i$ ($i=1,2$) as above, are invariant under the center and satisfy
$|F_i(a(\varphi^\alpha))| \leq \bar{F}_i(a(\varphi^\alpha))$ for every $\alpha$, then also
$|\peter{F_{1},F_{2}}_{N\bash G}|\leq \peter{\ov F_{1},\ov F_{2}}_{N\bash G}.$

We claim that  $|J_{s,t}(W_s \otimes f_t) a(\varpi^{\alpha})| $
and $|J(W \otimes f)(a (\varpi^{\alpha})|$ are both bounded by 
$2 J_{\theta,\delta}(W^0_{\theta+\delta} \otimes f^0 (a(\varpi^{\alpha})$,
the latter quantity being non-negative. This will complete the proof of the claim
\eqref{boundwanted}. 
 
By choice of data, all quantities vanish when $\alpha < 0$. 
When $\alpha \geq 0$, they may be explicitly evaluated. 
Our assertion amounts to the following:

\begin{equation} \label{olivIa} \Bigl( \sum_{\stackrel{ k\in[-n,n]}{ 2|(k-n)}} q^{k (\theta+ \delta/2)}  \Bigr) q^{n \delta/2}\leq 2 \Bigl( \sum_{\stackrel{ k\in[-n,n]}{ 2|(k-n)}} q^{k (\theta+\delta)}  \Bigr),  \mbox{$n\geq 0$, $q\geq 2$ and $\delta,\theta\geq 0$.}\end{equation}

\qed

\subsection{Bounds for the normalized Hecke functionals.} \label{ugly}
Let $\chi$ be a character of $k^{\times}$, and $\pi$ an irreducible unitary representation of $G$
with central character $\omega$. It is known that the space of functionals $\ell^{\chi}$
on $\pi$ satisfying
\begin{equation} \label{lchiprop} \ell^{\chi}\left( \left( \begin{array}{cc} a & 0 \\ 0 & d \end{array} \right) v \right) = \chi(a).\omega\chi^{-1}(d) \ \ell^{\chi}(v).\end{equation} 
 is at most one-dimensional.  Our goal is to normalize an element in this space (using the unitary structure on $\pi$) up to a scalar of absolute value $1$. Once this is done, we shall study
 the analytic properties of the resulting ``normalized'' functional.

\subsubsection{The normalized Hecke functional}

As before, {\em we fix once and for all an identification of $\pi$ with $\kirill(\pi, \psi)$
which carries the unitary structure on $\pi$ to the unitary structure on $\kirill(\pi, \psi)$
given by \eqref{whitinner}. } This being understood, define a functional 
$\ell^{\chi}$ on the Kirillov model $\kirill(\pi, \psi)$ --
and therefore on $\pi$ -- via the rule of \eqref{lchidef}  (interpreted by meromorphic continuation in general).   In the case $\chi$ is the trivial character,
we denote $\ell^{\chi}$ simply as $\ell$; thus
$\ell(W) = \int_{k^{\times}} W(y) d^{\times}y$.  The resulting functional satisfies \eqref{lchiprop} 

We have thus normalized a functional in the one-dimensional space of solutions
to \eqref{lchiprop},  up to a scalar of absolute value $1$;
this ambiguity arises
from the ambiguity in choosing the isometry between $\pi$ and its Kirillov model. 

Equivalently, the functional may be expressed in terms of matrix coefficients: one has
$$\ell^{\chi}(v)\ov{\ell^{\ov\chi^{-1}}(v)}  = \int_{ k^{\times}} \langle a(y)  v, v \rangle \chi(y)\dti y,$$
as follows from \eqref{whitinner}. 
In particular, if $\chi$ is unitary, the {\em associated Hermitian form} $w^{\chi} := |\ell^{\chi}|^2$
depends {\em only} on the unitary structure on $\pi$ -- i.e., there is not even an ambiguity of absolute value $1$ -- and can be expressed thus:
\be\label{matrixcoeff}h^{\chi}(v)  = \int_{ k^{\times}} \langle a(y)  v, v \rangle \chi(y)\dti y, \ \ \ v \in \pi. 
\ee

\begin{Lemmat}[\cite{Ve}, Lemma 11.7]\label{lem:hjllocalnonarch} 
Suppose $k$ non-archimedean; let $v$ the new vector,
and let $r$ be the conductor of $\chi$. Then:
$$ \ell^{\chi}(n(\varpi^{-r}) v)  =
\begin{cases}L(\pi\otimes\chi,1/2)&r=0\\
\eta&r>0\end{cases}$$
where $\eta$ has absolute value $q^{-r/2}(1-q^{-1})$. 
\end{Lemmat}

\begin{Lemmat} \label{lem:hjllocal}
Suppose $k$ archimedean. There is $v  \in \pi$, with the property that $\Sob_d^{\pi}(v) \ll_d C(\pi)^{2d}$ for 
every $d$,  and so that, for any $\eps>0$ there is $T \in k$
with  $\frac{ \log |T|}{\log C(\chi)} \in [1-\varepsilon, 1+\varepsilon]$ so that, 
$$\left| \ell^{\chi}(n(T)v)\right| \gg_{\eps} C(\chi)^{-1/2-\varepsilon}.$$
Moreover,  $\ell^{\chi_s}(n(T) v) \ll_{\epsilon} C(\chi)^{-1/2+\epsilon}$
for any $s$, the implicit constant uniform when $\Re(s)$ remains in a compact set. 
\end{Lemmat}

\proof
We take $v=W$ to
be (considered in the Kirillov model of $\pi$) a smooth bump function supported
in the ball of radius $1/10$ centered at $y=1$. (For concreteness, fix once and for all
such a function both in the case of $\R$ and $\C$. The implicit constants in the
result will depend on this choice.) 

 The bounds on Sobolev norms
are a consequence of the comments in \S \ref{seckirillov}. 
Noting that,    
$$\ell^{\chi}(n(T)v)=\int_{k} W(y) \psi(T y) \chi(y)d^{\times}y,$$ 
the assertions concerning upper and lower bounds for $\ell^{\chi}$ follow from Lemma \ref{almostgauss}. 
 \qed

\subsection{Normalized trilinear functionals.} \label{ugly:trilinear}
Let $\pi_{i}$, for $1 \leq i \leq 3$, be three unitary irreducible generic representations of $G$,
the product of whose central characters ($\chi_{i},\ i=1,2,3$) is $1$.
 In the sequel, we denote by $\omega$ the product $\chi_{1}\chi_{2}$.
 It is known (\cite{prasad}, \cite{Oksak}, \cite{Loke}) that the space of 
 trilinear invariant functionals
$$\pi_1 \otimes \pi_2 \otimes \pi_3 \rightarrow \C$$
is at most dimension $1$. Again, we wish to normalize one up to scalars of absolute value $1$,
and study the analytic properties of the resulting ``normalized'' functional.

\subsubsection{Integration of matrix coefficients} \label{subsec:wii}
We take the following normalization (cf. \cite{Wald}, \cite{Ich}, \cite{SaVe}).\be\label{RSmatrixlocal}|L_{W}(x_1, x_2, x_3)|^2 = \int_{g \in \bar{G}} \langle g x_1, x_1 \rangle \langle g x_2, x_2 \rangle
\langle g x_3, x_3 \rangle dg.
\ee
By \eqref{localglobalbounds} together with the convention \S \ref{convention}, the integral is convergent, and, indeed, 
bounded by a constant multiple of $\prod_{i=1}^3 \Sob^{\pi_i}(x_i)$. 
It is true, although not obvious, that the right-hand side is positive.
By direct computation (see \cite{Ich}), \begin{equation} \label{ich-form} \frac{|L_W|^2}{\prod_i \|x_i\|^2} = \zeta_k(2)^2 \frac{ L(\pi_1 \otimes \pi_2 \otimes \pi_3,\frac{1}{2} ) }{\prod_{i=1}^{3} L(\pi_i, \Ad,1)} \end{equation}
if $k$ is nonarchimedean and all vectors are unramified.

Another functional arises naturally from Rankin-Selberg integrals.
Suppose  that $\pi_{3}\simeq \chi_{3}^{+}\boxplus\chi_{3}^{-},$  is a principal series representation ($\chi_{3}^{+}\chi_{3}^{-}=\chi_{3}$ ) realized in the model $\mcI_{3}$,
and that $\pi_{i}$ $i=1,2$ are realized in the respective Whittaker models $\Whit(\pi_{1},\psi)$ and $\Whit(\pi_{2},\ov\psi)$.  (More precisely, 
we follow the convention of fixing once and for all an equivariant isometry between $\pi_i$ and
 their respective Whittaker models, equipped with the inner product \eqref{whitinner};
thus, for instance, we will regard a linear functional on the Whittaker model $\Whit(\pi_i)$
as a linear functional on $\pi_i$. )

Then,  for $W_i \in \Whit(\pi_i), f_3 \in \mcI_{3}$, we put
\begin{eqnarray} L_{RS}(W_1, W_2, f_3)&&= \zeta_k(1)^{1/2} \int_{N \backslash \bar{G}}  W_1{W_2} f_3, \label{RSlocal}\\ 
&&= \zeta_k(1)^{1/2} \int_{K}\int_{k^\times}W_{1}(a(y)k)W_{2}(a(y)k)f_{3}(a(y)k)|y|^{-1}\dti ydk
\nonumber
\end{eqnarray}
If $W_i, f_3$ are spherical, and $k$ nonarchimedean, then 
\begin{equation} \label{aers} \frac{L_{RS}(W_1, W_2, f_3)}{W_1(1) W_2(1) f_3(1) } = \zeta_{k}(1)^{1/2} 
\frac{L(\pi_1 \otimes \pi_2\otimes\chi_3^{+},{1/2})}{L(\chi_3^{+}/\chi_3^{-},1)}. \end{equation} If
we specialize to $\chi_3^{\pm}  = |\cdot|^{\pm 1/2}$, $W_{2}=\ov{W_{1}}$ then $L_{RS}=\zeta_{k}(1)^{1/2}\|W_{1}\|^2\int_{K}f(k)dk$, giving in particular for spherical vectors
$$\frac{ \|W_i\|^2}{|W_{i}(1)|^2} = \zeta_{k}(2)^{-1}  L(\pi_1 \times \overline{\pi_i}, 1)  = \frac{L( \pi_i , \Ad,1)}{\zeta_{k}(2)}  {\zeta_{k}(1)}.$$  We conclude that, in the unramified case, 
$$\frac{|L_{RS}|^2}{\prod_{i=1}^3 \|x_i\|^2} = \zeta_k(2)^2 \frac{ L(\pi_1 \otimes \pi_2 \otimes \pi_3,\frac{1}{2} ) }{\prod_{i=1}^{3} L(\pi_i, \Ad,1)}.$$

\begin{Lemmat}
Let $L = L_{RS}$ or $L_W$ be defined as in either \eqref{RSmatrixlocal} or \eqref{RSlocal}. 
Then for $x_1 \in \pi_1, x_3 \in \pi_3$:
\begin{equation} \label{planch-D}  \| x_1\|^2 \|x_3\|^2 =  \int_{\pi_2} d\mu_P(\pi_2)  \sum_{x_2 \in \mathcal{B}(\pi_2)} | L(x_1, x_2,  x_3)|^2,\end{equation} 
where $d\mu_P(\pi_2)$ is the Plancherel measure on the unitary dual of $G$, 
and $\mathcal{B}(\pi_2)$ is an orthonormal basis for $\pi_2$. 
In particular, $|L_{RS}|^2 = |L_W|^2$ 
whenever $\pi_2$ is tempered.  
\end{Lemmat}

The proof that follows is inspired by joint work of one of the authors (A.V.) with Y. Sakellaridis. 
 
\proof
The validity of \eqref{planch-D} for $L = L_{W}$ is an 
immediate consequence of the Plancherel formula.

The validity of \eqref{planch-D}, where $L= L_{RS}$ \eqref{RSlocal} follows
from the first Lemma of \S \ref{Whitmodel-0} together with the following:
it is known (the Whittaker-Plancherel theorem) that
for  $F \in L^2(N \backslash \bar{G})$,
\begin{equation} \label{Whit-planch} \|F\|_{L^2(N \backslash \bar{G}) }^2 =  \zeta_k(1) \int_{\pi} d \mu_P(\pi) \sum_{x \in \mathcal{B}(\pi)}
|\langle F, W_{x} \rangle|^2.\end{equation}
where, for each $\pi$, we choose an orthonormal basis $\mathcal{B}_{\pi}$
and a unitary intertwiner from $\pi$ to its Whittaker model $\Whit(\pi)$.
In fact, the key point in the validity of \eqref{Whit-planch} is the equality , for $W_i \in \Whit(\pi)$ between $\int_{x \in k} \psi(x) \langle
n(x) W_1, W_2 \rangle $ and $\zeta_k(1) W_1(1) \overline{W_2(1)}$. 
This follows from the usual Fourier inversion formula, taking into account that the Haar measure
on $k$ was chosen to be self-dual w.r.t. $k$.

 The second assertion (the coincidence of the two definitions of $|L|^2$ when $\pi_2$ is tempered)
 follows from uniqueness of unitary decomposition, applied to $\pi_1 \otimes \pi_3$,
 together with the fact  that the support of the Plancherel measure is the full tempered spectrum, 
and from a continuity argument. 
  \qed 

\begin{Rem} It may also be deduced that $|L_{RS}|^2$ and $|L_W|^2$
coincide even when not all $\pi_i$ are tempered.  To do this, write $\Pi = \otimes_{i=1}^{3} \pi_i$. 
 Both $|L_W|^2$
and $|L_{RS}|^2$ are Hermitian forms on $\Pi$ and can be polarized to 
maps $\Pi \otimes \tilde{\Pi} \rightarrow \mathbb{\C}$.   These polarized maps vary holomorphically in parameters, 
and therefore their coincidence may be extended from tempered representations by analytic continuation. 

\end{Rem}

\begin{Rem}
To give a trilinear functional on $\pi_1 \otimes \pi_2 \otimes \pi_3$ is equivalent to describing an equivariant intertwiner  $\pi_{1} \otimes \pi_{3} \rightarrow \widetilde{\pi_{2}}$. 
Dualizing, we denote by $I_{\pi_2}$ the functional
$I_{\pi_{2}}: \pi_{1} \otimes \pi_{3} \rightarrow \pi_2$
that arises from the normalized trilinear functional on $\pi_1 \otimes \widetilde{\pi_2} \otimes \pi_3$. 
In explicit terms, 
$$L(x_1 \otimes x_2 \otimes x_3) = \langle I_{\pi_2}(x_1 \otimes x_3), x_2 \rangle,$$
and one has
\be\label{normI}\|I_{\pi_{2}}(x_{1}\otimes x_{3})\|^2=\sum_{x_{2}\in\mcB(\pi_{2})}|L(x_1 \otimes x_2 \otimes x_3)|^2.
\ee
It is often convenient to think or phrase arguments in terms of $I_{\pi_2}$, instead of $L$. In the sequel, the notation $I_{\pi}$
will designate the intertwiner associated with $L_{RS}(\star,\pi,\star)$ or $L_{W}(\star,\pi,\star)$.

It should also be noted that Bernstein and Reznikov made extensive use of a local
study of such trilinear functionals in their beautiful papers. 
\end{Rem}

\subsection{Bounds for trilinear functionals, I: soft methods}\label{softsection} 

In this section we present some ``soft'' upper bounds on the normalized trilinear functionals.
In fact, the entire section is not strictly necessary for the present paper, since we reprove the bounds by brute force
in later sections. However, {\em under the Ramanujan conjecture}, 
one does not need the brute-force computations, so the reader may prefer
to focus on this section at first. 
 
Our main result is Lemma \ref{tempbound}. It bounds the Sobolev norms of $I_{\pi_3}(x_1 \otimes x_2)$
in terms of a certain norm $\|\cdot\|_{\Knorm}$ on the space of $\pi_1$ and $\pi_2$, assuming that all three representations
 are tempered. 

For the rest of this section, we will assume that the local field $k$ is {\em non-archimedean};
the proof in the archimedean case is similar.

\subsubsectionind{}\label{localnorms}

Suppose $V$ is any unitary $G$-representation. We define\footnote{If $K$ were a compact abelian group and $x \in L^2(K)$, this coincides
with the Gowers $U^4$-norm.} $$\|x\|_{\Knorm} = \left(  \int_{k \in K} | \langle k x, x \rangle|^2 \right)^{1/4},$$
where $K$ is endowed with Haar probability measure.  

\begin{Lemma*} $\|x\|_{\Knorm}$ is a norm.
\end{Lemma*}
\proof 
Let $v \mapsto v^K$ be the averaging operator over $K$. Then $\|x\|_{\Knorm}^2 = \| (x \otimes x)^K \|_{V \otimes V} $. Now,
Cauchy-Schwarz demonstrates that
\begin{equation} \label{cso} \|(x\otimes y)^K\|_{V \otimes V}^2 \leq  \| (x \otimes x)^K \|_{V \otimes V}  \| (y \otimes y)^K \|_{V \otimes V}.\end{equation} Since
$(x+y) \otimes (x+y) = x\otimes x+ x\otimes y + y \otimes x + y \otimes y,$ we conclude. 
\qed

\begin{Lemma*}There exists $d < d' < 0$ so that
$$\Sob^{V}_d \leq \|\cdot\|_{\Knorm} \leq \Sob^{V}_{d'}$$ on any unitary representation. 
\end{Lemma*}

We omit the proof, for this is never used in the present paper.

\begin{Lemma*} Suppose $g \in G$; then $\|g x\|_{\Knorm} \leq C(g) \|x\|_{\Knorm}$, where $C(g)$
is a continuous function of $g$; we may take
$C = [K g K: K]^{1/4}$. 

\end{Lemma*}
\proof
Write $U = g K g^{-1} \cap K$;
we equip it with the restriction of the measure from  $K$. 
Then $\int_{U} | \langle u g x , gx \rangle|^2 \leq \|x\|_{\Knorm}^4.$ 
Thus
\begin{equation}\label{uav}\ \|(gx \otimes g x)^U\|^2 \leq [K:U] \|x\|_{\Knorm}^4,\end{equation} where, once again,
$w \mapsto w^U$ denotes the projection onto $U$-invariants. 
Clearly, $\|(gx \otimes gx)^K\| \leq \|(gx \otimes g x)^U\|$ (averaging decreases norms!)
and we are done since $[K:U] = [K g K: K]$. 
\qed

 \begin{Lemmat} \label{ipibound}
For $x_1, x_2 \in \pi_1, \pi_2$ and  for any $d\geq0$ we have:
\begin{equation} \label{trisob} \Sob^{\pi_3}_d(I_{\pi_3}(x_1 \otimes x_2)) \ll_{d} \Sob_{d'}^{\pi_1}(x_1) \Sob_{d'}^{\pi_2}(x_2).\end{equation}
\end{Lemmat}
where $d'$ depend on $d$ only.
\proof 
We have seen after \eqref{RSmatrixlocal} that there exists an absolute $d_0$ so that
\begin{equation} \label{trivbond} \langle I_{\pi_3}(x_1 \otimes x_2), x_3 \rangle \ll   \prod_{i=1}^{3} \Sob^{\pi_i}_{d_{0}}(x_i),\end{equation}
which shows, in particular, that $\Sob^{\pi_3}_{-d_0}(I_{\pi_3})$ is
bounded by a quantity as in the right-hand side of \eqref{trisob}. 
However, we wish the same result where $-d_0$ is replaced by an arbitrary positive quantity. 
To establish this, observe first that, by equivariance, 
that $x_i \in \pi_i[m] \implies I_{\pi_3}(x_1 \otimes x_2) \in \pi_3[m]$. 
It follows, then, that a bound of the form \eqref{trisob} is true when $x_1 \in \pi_1[j], x_2 \in \pi_2[k]$
for any $j,k$. 
Now apply \eqref{rt1} (twice: first fix $x_1$, and then fix $x_2$) to conclude the same result
for arbitrary $x_1, x_2$. 
\qed

\begin{Lemmat} \label{tempbound}       
Suppose $\pi_i$ are all tempered.  For any $\eta > 0$, set $B_{\eta} = \{g \in G: |\Xi(g)| \geq \eta\}$. Then
there exists an absolute constant $d_0\geq 0$ (independent of $k$) so that:
$$\Sob_{-d_0}^{\pi_3} ( I_{\pi_{3}}(x_1 \otimes x_2))  \ll_{\eps} (\vol B_{\eta})^{\eps}  \|x_1\|_{\Knorm} \|x_2\|_{\Knorm} + 
\eta^{1-\eps}  \Sob^{\pi_1}(x_1)
\Sob^{\pi_2}(x_2)$$
\end{Lemmat}

The reader should ignore the second term at a first reading; the content is that $I_{\pi_{3}}(x_1 \otimes x_2)$
is bounded ``up to $\epsilon$s'' by $\|x_1\|_{\Knorm} |\|x_2\|_{\Knorm}$.

\proof 
 It is sufficient to show that there exists $n \geq 0$ so that,  whenever $x_3$ is invariant under an open compact subgroup $U \subset K$, 
$| L(x_1, x_2, x_3) |  \leq [K:U]^n \cdot   \mathrm{RHS},$ where $\mathrm{RHS}$ is the right-hand side
of the stated bound.

We endow both $U$ and $K$ with Haar probability measure. 
By a modification of \eqref{cso}, $\|(x_1 \otimes x_2 )^U \|^2  \leq [K:U]
\|x_1\|_{\Knorm}^{2} \|x_2\|_{\Knorm}^{2}$. Writing
$\prod_i \langle gx_i, x_i \rangle $ as $$\langle g (x_1 \otimes x_2)^U, (x_1 \otimes x_2)^U \rangle \langle g x_3, x_3 \rangle,$$ we deduce from bounds on matrix coefficients (see \S \ref{mc-local}, esp. footnote):

$$\prod_i \langle gx_i, x_i \rangle  \leq \begin{cases} \Xi(g)^2 [K:U]^2 \|x_1\|_{\Knorm}^2 \|x_2\|_{\Knorm}^2 \|x_3\|^2,  \\ \ \Xi(g)^3 \prod_{i} \Sob^{\pi_i}(x_i)^2,  \end{cases}$$
To analyze $|L(x_1, x_2, x_3)|$, we split the $g$-integral into $B_{\eta}$ and its complement,
applying the former bound on $B_{\eta}$ and the latter on its complement. This leads to the stated bound;
the bounds on $\Xi(g)$ are obtained using only the fact that 
$$\int_{g \in G} |\Xi(g)|^{2+\eps} \ll_{\eps} 1.$$

\begin{Rem}
When $\pi_3$ is not necessarily tempered, an analysis of the prior argument shows that, if $x_i$ are invariant by subgroups $U_i$, then $|L(x_1, x_2, x_3)|$ is bounded above by
\begin{equation} \label{non-tempered-bound} 
[K:U_3]  q^{m \theta/2}   \|x_1\|_{\Knorm} \|x_2\|_{\Knorm} \|x_3\|  +  
\prod_{i} [K:U_i]^{1/2} q^{(\theta- 1/2) (m+1)} \|x_1\| \|x_2\| \|x_3\|\end{equation}
\end{Rem}

\subsection{Choice of test vectors in the trilinear form.}
In this section, we estimate $L$ and $I_{\pi_{3}}$ for specifically chosen test vectors.

\begin{Prop} \label{lowerbounds}
Let $\pi_1, \pi_2$ be the two generic unitary representations of $\GL_2(k)$. 
Suppose that $\pi_3$ is the principal series representation $1 \boxplus \chi_{3}$,
where $\chi_{3}^{-1}$ is the product of central characters of $\pi_1, \pi_2$. 

Write the analytic conductors $C_i = C(\pi_i)$. 
Then, for any $\epsilon > 0$, 
there exists vectors $x_i \in \pi_i$ with the following properties:
\begin{itemize}
\item[-] $\|x_i\| = 1$ ;
\item[-] For $d\geq 0$, $\Sob^{\pi_2}_d(x_2) \ll_{d,\pi_2} 1$;
\item[-] For $d \geq 0$, $\Sob^{\pi_3}_d(x_3) \ll_{d, \pi_1, \pi_2} 1$; 
\item[-] $\|x_3\|_{\Knorm} \ll_{\pi_2,\epsilon}  C_1^{-1/4+\epsilon}$;
\item[-] $|L(x_1, x_2, x_3)| \gg_{\pi_2, \epsilon}  C_1^{-1/2-\epsilon}$;
\end{itemize}
\begin{itemize}
\item[-] If $k$ is non-archimedean and all the $\pi_i$ are unramified, we may take $x_1, x_2, x_3$
to be the new vectors. In this case,  $$\frac{|L(x_1, x_2, x_3)|^2}{\prod \|x_i\|^2}  = \zeta_k(2)^2 \frac{ L(\pi_1 \otimes \pi_2 \otimes \pi_3,\frac{1}{2})}{\prod_{i} L(1, \pi_i, \Ad)}.$$
\end{itemize}

\end{Prop}
The rest of this section indicates the choice of test vector and the proof of the assertions.
The unramified case has already been noted.

\subsubsection{Case of $k$ nonarchimedean}
Write $f$ to be the larger of the valuations of the conductors of $\pi_1, \pi_2$. We set:\begin{enumerate}
\item $x_1 \in \pi_1$ to be the new vector in the Whittaker model $\Whit(\pi_{1},\psi)$, normalized
to have norm $1$; 
\item $x_2 \in \pi_2$ to be the new vector in the Whittaker model $\Whit(\pi_{2},\ov\psi)$, normalized
to have norm $1$; 
\item $x_3 \in \pi_3$ to be the vector in the principal series representation corresponding to the function $f_{3}\in\mcI_{3}$ of norm $1$ which upon restriction to $K$, is supported on the compact subgroup $K_0[f]$, and
is a $K_0[f]$-eigenvector: explicitly
\begin{equation} \label{f3nonarchdef} f_{3|K}:\left( \begin{array}{cc} a & b \\ c & d \end{array} \right) \mapsto \vol(K_0[f])^{-1/2}\chi_{3}(d)1_{K_{0}[f]}.\end{equation}
\end{enumerate}
The ratio $\frac{|W(1)|}{\|W\|}$, when $W$ is the new vector, is well-known (e.g. \cite[\S 7]{VenAdv}), 
and in particular is absolutely bounded above and below at nonarchimedean places;
the assertion on $\Sob^{\pi_i}_{d_i}(x_i) \ (i=2,3)$ follows from the fact that $x_i$
is invariant under $K[f_i]$, with $f_i$ the conductor of $\pi_i$. 

Any $K$-translate of $x_3$ is either orthogonal or proportional to $x_3$, and therefore
$$\|x_3\|_{\Knorm} = \vol(K_0[f])^{1/4}\ll_{\pi_2}C_1^{-1/4}.$$ The lower bound
for $L$ follows from \eqref{RSlocal}
taking into account that the function $k\mapsto W_1(a(y)k)W_2(a(y)k)f_3(k)$ is supported on $K_0[f]$ and is $K_0[f]$-invariant:
\begin{eqnarray*}L(x_{1},x_{2},x_{3})&&=\int_{K}\int_{k^\times}W_{1}(a(y)k)W_{2}(a(y)k)f_{3}(k)|y|^{-1/2}\dti ydk \\
&&=\vol(K_{0}[f])^{1/2}\int_{k^\times}W_{1}(a(y))W_{2}(a(y))|y|^{-1/2}\dti y \gg C_1^{-1/2}
\end{eqnarray*}
Thus in the non-archimedean case, the proposition is valid with $\epsilon=0$.

\subsubsection{Case of $k$ archimedean: the tempered case} We assume that $\pi_{2}$ is {\em tempered}.
 Let $\varphi$ be a smooth function on $k^{\times}$, with $\varphi(1)=1$, which is
supported in a fixed compact neighbourhood of the identity and so that $\int |\varphi|^2 d^{\times}y = 1$
and $\int |\varphi|^2 |y|^{-1/2} d^{\times}y =1$. 
Fix $\delta, M > 0$; later, we will choose $\delta$ to be ``small'', to an extent depending on the parameter
$\epsilon$ chosen in Proposition \ref{lowerbounds}, whereas $M$ will be a ``large'' but absolute constant. 

We set \begin{equation} \label{Cdef} C := C_2^{M} \max(C_1,C_{2})^{1+\delta}. \end{equation}
Intuitively, the reader should think of $C$ as being ``a tiny bit larger than $C_1$.'' 
It will be convenient to define
\begin{equation} \label{KCdef}
\mcK_C :=  \{ g\in K : \ \ |g-k|\leq C^{-{1}} \mbox{ for some } k\in K\cap B\} \end{equation}
where $|g-k|$ simply denotes the largest absolute value of any matrix entry of $g-k$. 
Thus $\mcK_C$ is a small neighbourhood of $K \cap B$ whose volume, with respect
to the Haar measure on $K$,  is $\asymp C^{-1}$.  

If $C$ is large enough, a matrix $g= \left(\begin{array}{cc} a & b \\ c & d \end{array} \right) \in K$ belongs to $\mcK_C$ iff $|c| \leq C^{-1}$; indeed, if this is so, $|b| = |c|$ and $|a| = |d| = \sqrt{1-|c|^2}$, from where we see
that $|g -k| \leq C^{-1}$, where $k$ is the diagonal matrix with entries $a/|a|, d/|d|$. %

The vectors $x_{1},x_{2},x_{3}$ are as follows:
\begin{enumerate}

\item $x_1 \in \pi_1$ is $y \mapsto \varphi(y)$ in the Kirillov model $\kirill(\pi_{1},\psi)$:  $W_{1}(a(y))=\varphi(y)$ in the Whittaker model.
\item $x_2 \in \pi_2 $ is $\overline{\varphi(y)}$ in the Kirillov model $\kirill(\pi_{2},\ov\psi)$: $W_{2}(a(y))=\ov\varphi(y)$ in the Whittaker model.
\item We give an explicit construction of $x_{3} \in\pi_{3}$ below: for the moment it is sufficient to say that  correspond to  a function $f_{3}\in \mcI_{3}$  which upon restriction to $K$
is a smooth bump function around $K \cap B$, supported in $\mcK_C$. 

We advise the reader
to skip the details of the explicit construction below; the properties just mentioned
are all that is needed to follow the proof. \end{enumerate}

\begin{Defi} \label{X3def} (Construction of $x_3$). 
Let $\phi$ be a (non-negative) bump function on $\Rr$, taking value $1$ at zero,
and supported in $|x|\leq \delta_0$. 
Set $\Psi_{3}(x,y)=\phi_{1}(x)\phi_{2}(y)$,  where
\begin{equation} \label{phiiarch} \phi_{1}(x)=\phi(C|x|),\ \phi_{2}(y)=\chi_{3}(y).\phi(|y|-1).\end{equation}
Then, for sufficiently small (absolute) $\delta_0>0$, the restriction to $K$ of
 \begin{equation} \label{f3intdef}
f(g)=|\det g|^{1/2}\int_{k^\times}\Psi_{3}((0,t)g)\chi_{3}^{-1}(t)|t|\dti t.\end{equation}
is supported in $\mcK_C$,  and has $\|f\|_{L^2(K)} \sim C^{-1/2}$. We take $f_3 = \frac{f}{\|f\|_{L^2(K)}}$,
and $x_3 \in \pi_3$ the vector corresponding to $f_3$.  
\end{Defi}

\proof 
For $g=\left(\begin{array}{cc} a&b \\c&d\end{array}\right)\in K$, we note
 that $\Psi_{3}(ct,dt)\not=0$ implies (for sufficiently small $\delta_0$) that $|dt|\geq (9/10)^2$ so $|t|\geq (9/10)^2$ (since $|d|\leq 1$)
  and then $|c|\ll C^{-1}$. From this we deduce -- for sufficiently values of $\delta_0$ --  that $f_{3|K}$ is supported in $\mcK_C$.

 Moreover, 
 $f(g) =  \zeta_k(1) \chi_{3}(d) \int_{t} \phi(C|c t|) \phi(|dt|- 1) dt$. From 
 this expression, it is easy to verify $\|f\|_{L^2(K)} \sim C^{-1/2}$;
 later on, we shall also use the fact that, at least if $M$ is sufficiently large, 
\begin{equation} \label{usefulsoon} |\int_{\mcK_C} f(k) | \geq \frac{1}{2} \int_{\mcK_C} |f(k)|dk \asymp C^{-1}.\end{equation}
 The point in this inequality is that $\omega(d)$ remains close to $1$ so long as
 $g \in \mcK_C$.  It follows that corresponding inequalities hold for $f_3$, 
but with the last expression replaced by $C^{-1/2}$.  
     \qed

We remark for later usage that, if $v$ is nonarchimedean and we take
\begin{equation} \label{phiinonarch} \phi_{1}=\vol(K_{0}[f])^{-1/2}1_{\varpi^f\mfo},\ \phi_{2}=\omega^{-1}.1_{\mfo^\times},\end{equation}
then the function $f_3$ defined by \eqref{f3intdef} coincides with that
defined in \eqref{f3nonarchdef}. 
\footnote{By construction, the integral belongs to $\mcI_{3}$ so it is sufficient to check the equality for $g= \left(\begin{array}{cc} a & b \\ c & d \end{array}\right)\in K$. In that case
the right hand side equals
$$\int_{k^\times}\Psi_{3}(ct,dt)\chi_{3}^{-1}(t)|t|\dti t.$$
Note that  $\max(|ct|, |dt|) = |t|$ since $g \in K$; so the integrand is nonvanishing only when
$|t| = 1$. In that case it is nonvanishing only if $d \in \order^{\times}$ and
 $c \in \varpi^{f} \order$, in which case the integral equals $\vol(K_{0}[f])^{-1/2}\chi_{3}(d) $ proving our claim. }

The assertions concerning $\|x_1\|, \|x_2\|$ are immediate; the
assertions concerning the Sobolev norms of $x_1, x_2$ follow from 
discussion of \S \ref{seckirillov}.  

By \refs{RSlocal}, we have
\be\label{Larchirep}L(x_{1},x_{2},x_{3})=\int_{K}\ell(k.W_{1},k.W_{2})f_{3}(k)dk
\ee
where $\ell(W, W')$ denote  the (Borel-equivariant) functional on 
$\Whit(\pi_{1},\psi)\otimes\Whit(\pi_{2},\ov\psi)$ given by
$$\ell(W, W')  := \int_{y \in k^{\times}} W(a(y)) {W'(a(y))}  |y|^{-1/2} d^{\times}y.$$

Write
\begin{multline*}
\ell(k.W_{1},k.W_{2})=\ell(W_1, W_2) + \ell(W_1, k.W_2 - W_2)\\ + \ell(k.W_1-W_1, W_2)+\ell(k.W_1-W_1, k.W_2-W_{2})
\end{multline*}

We need to bound the terms only for $k \in \mcK_C$; they are individually bounded thus:
\begin{enumerate}
\item We have set things up so that $\ell(W_1, W_2) \geq 0.9$;
\item By Cauchy-Schwarz and \eqref{meanvalue}\footnote{using the fact that $y\mapsto W_{i}(a(y))$ is supported on a fixed compact set of $k^\times$},
 \begin{eqnarray*}|\ell(W_1, k.W_2 - W_2)| &&\ll \|W_{1}\|\|k.W_2 - W_2\|\\
 &&\ll \|W_{1}\|C^{-1/\deg(k)}C_{2}^{1/\deg(k)}\\
 &&\ll_{\pi_{2}} C^{-1/\deg(k)}.
 \end{eqnarray*}
\item 
Similarly we have $|\ell(k.W_1-W_1, W_2)| \ll \|W_{2}\|C^{-1/\deg(k)}C_{1}^{1/\deg(k)}$. 
\item As for the last term: set $$W(y)=(k.W_1-W_1)(a(y)),\ W'(y)=(k.W_2-W_2)(a(y)).$$ Given $\eps>0$
small enough, we have, since {\em  $\pi_{2}$ is tempered},
by Proposition \ref{zeroboundprop}, 
\begin{align*}\bigl|\int W(y)W'(y)|y|^{-1/2}\dti y\bigr|&\leq \big(\int|W(y)|^2|y|^{-\eps}\dti y\big)^{1/2} \big(\int|W'(y)|^2|y|^{-1+\eps}\dti y\big)^{1/2}\\
&\ll\Sob^{\pi_{2}}(W')\left(\|W\|+\|W\|^{1-O(\eps)}\Sob^{\pi_{1}}(W)^{O(\eps)}\right)\\
&\ll  (C_{1}/C)^{(1-d\eps)/\deg(k)}
 C_{1}^{d'\eps}\end{align*}
 for some $d' \geq 0$. 
 Here, to bound $\int|W(y)|^2|y|^{-\eps}\dti y$, we split into $|y| \geq 1$,
 where the bound $\|W\|$ suffices, and $|y| \leq 1$, which gives --
 by Proposition \ref{zeroboundprop} and interpolation -- a term 
 $\|W\|^{1-O(\eps)}\Sob^{\pi_{1}}(W)^{O(\eps)}$; 
we also used $\Sob^{\pi_2}(W') \ll_{\pi_2} 1$, as well as 
\eqref{meanvalue}.
\end{enumerate}

It follows by taking $M$ large enough, we will have $|\ell(k.W_{1},k.W_{2})-1|\leq \frac{1}{100}$
for all $k \in \mcK_C$. 
On the other hand, $|\int_{\mcK_C} f_3(k)| \geq  \frac{1}{2} \int_{\mcK_C} |f_3(k)|$.
Therefore, 
$$L(x_{1},x_{2},x_{3})\gg \left|\int_{\mcK_C} f_3(k)\right| \gg_{\pi_{2}} C_{1}^{-1/2-\delta/2}.$$

We choose $\delta \leq \epsilon$ to conclude. 

\subsubsection{The non-tempered case}

If $\pi_{2}$ is {\em not tempered}, then $\pi_{2}$ is an unramified principal series representation
and may be written into the form $\pi^0_{2}.|.|^{it}$ with $t\in\Rr$ and $C(\pi^0_{2})\leq 9$. 
 Let $W^0_{2}$ be the $K$-invariant vector in the Whittaker model of $\pi^0_{2}$,
 normalized to have $L^2$-norm $1$. Since the 
 set of possibilities
 for the isomorphism class of $\pi^0_{2}$ is precompact,
 there exists a finite set $\mathcal{F}$ of smooth non-negative function $\varphi$, each compactly supported in $[1/2,2]$ and with $\varphi_i(1)=1$, with the following property:
 
For any such $\pi_2$, there exists $\varphi \in \mathcal{F}$ such that:
 $$\int_{k^\times}\varphi_i(y)W^0_{2}(a(y))|y|^{-1/2}\dti y\gg 1$$
 where the constant implied is absolute.

 Let $\pi^0_2$ be as above, and let $\varphi \in \mathcal{F}$
 be the corresponding function.  We take $x_{1}$ to be the vector
 corresponding to $y \mapsto \varphi(y) y^{-it}$ in the Kirillov model of $\pi_{1}$;
 take $x_2$
 to be the vector corresponding to $W_2(y) |y|^{it}$ in the Kirillov model of $\pi_2 = \pi_2^{0} |\cdot|^{it}$.   In the previous notations, we have by definition
 $$\ell(W_{1},W_{2})\gg{1}$$
 Now for $k\in\supp(f_{3|K})$ we have by $K$-invariance of $W_{2}$, 
 $$\ell(k.W_{1},k.W_{2})=\ell(W_{1},W_{2})+\ell(k.W_{1}-W_{1},W_{2})$$
Since the term $\ell(k W_1-W_1,  W_2)$ is bounded by
$O((C_{1}/C)^{1/\deg(k)} C_2^{m})$, for some $m \geq 1$, we deduce, just as in the prior argument, that
 $$L(x_{1},x_{2},x_{3})\gg C_{1}^{-1/2-\delta/2}.$$
  \qed 

\subsection{Bounds for trilinear functionals, II: computations} \label{endpartIII}

In this section we shall establish the following result, 
{\em which can be replaced by the results of \S \ref{softsection}  under the assumption of the Ramanujan-Petersson conjecture}. 
Indeed central bound \eqref{equpperboundexplicit} 
follows from Lemma \ref{tempbound} if $\pi$ is tempered.   
The reader may, indeed, wish to omit the present section at a first reading, because it contains ``ugly'' computations.

\begin{Prop}\label{upperboundexplicit} Let $\pi$ be a generic unitary representation with trivial central character. 

Let $\pi_1, \pi_2$ be two generic unitary representations and $\epsilon > 0$. 

We denote the product of the central characters of $\pi_1, \pi_2$ by $\chi_{3}^{-1}$. 
Let $\pi_{3}=1\boxplus\chi_{3}$ be the corresponding unitary principal series representation.

Let $x_3 \in \pi_3$ be the vector constructed in Proposition \ref{lowerbounds} -- applied to $\pi_1, \pi_2, \pi_3, \epsilon$; let $\tilde x_{3}$ be the conjugate vector in 
in the contragredient representation $\widetilde\pi_{3}=1\boxplus\omega$ 
(i.e., the vector representing the functional $v \mapsto \langle v, x_3 \rangle$).
\begin{itemize} \item[-] For $\pi,\pi_j$ all unramified, and $x \in \pi$ the new vector of norm $1$, one has
$$\frac{ | L(x,x_{3},\tilde x_{3})|^2}{\|x\|^2 \|x_3\|^2 \|\tilde x_3\|^2}  ={\zeta_{k}(2)^2}\frac{L(\pi \otimes \pi_3 \otimes \tilde\pi_3,\frac{1}{2})}{L(\pi,\Ad,1)L(\pi_{3},\Ad,1)L(\tilde\pi_{3},\Ad,1)} $$
\item[-] In general, there is an absolute constant $d_{0}\geq 0$ (independent of $k$) such that: 
\be\label{equpperboundexplicit}| L(x,x_{3},\tilde x_{3})|^2 \ll_{\epsilon, \pi_{2}}C_{1}^{d_0 \epsilon} C(\chi_{3})^{-1/2}\left( \frac{C_{1}}{C(\chi_{3})}\right)^{\theta-1/2} \ \ \Sob_{d_{0}}^\pi(x). \ \ (x \in \pi). 
\ee
\end{itemize}
\end{Prop}
Averaging \refs{equpperboundexplicit} over an orthogonal basis $\mcB(\pi)$ of $\pi$ we obtain
\begin{Corollary*} With the previous notations, if $\pi,\pi_j$ are all unramified, one has
$$\frac{ \|I_{\pi}(x_{3}\otimes\tilde x_{3})\|^2}{\|x_3\|^2 \|\tilde x_3\|^2} ={\zeta_{k}(2)^2}\frac{L(\pi \otimes \pi_3 \otimes \tilde\pi_3,\frac{1}{2})}{L(\pi,\Ad,1)L(\pi_{3},\Ad,1)L(\tilde\pi_{3},\Ad,1)}.$$
In general, there is  an absolute constant $d_{0}\geq 0$ (independent of $k$) so that: \be\label{equpperboundexplicitII}\Sob^\pi_{-d_{0}}(I_{\pi}(x_{3}\otimes\tilde x_{3}))\ll_{\epsilon, \pi_{2}}C_{1}^{d_0 \epsilon} C(\chi_{3})^{-1/2}\left( \frac{C_{1}}{C(\chi_{3})}\right)^{\theta-1/2}.
\ee
\end{Corollary*}
We observe that the proof is parallel in nonarchimedean and archimedean cases,
and we could indeed have given an integrated treatment. The unramified assertion has
already been discussed.

\proof
Let $C = C_{2}^{M} \max(C_{1}, C_{2})^{1+\delta}$, where $M, \delta$
are as in the prior section (cf. \eqref{Cdef};
recall that $\delta$ is chosen depending on $\epsilon$, whereas $M$ is an absolute constant) when $k$ is archimedean,  and $C = q^f$ when $k$ nonarchimedean; here $f$ is the larger
of the conductors of $\pi_1$ and $\pi_2$. 
Let $\mcK_{C}$ be as in \eqref{KCdef} when $v$ is archimedean, and $K_0[f]$ when $v$ is nonarchimedean. 

We realize $\pi$ and $\pi_{3}$ in their  respective Whittaker models, i.e., fix
isometric intertwiners between $\pi, \pi_3$ and these models. Similarly, we realize $\tilde\pi_{3}$ in the induced model $\tilde\mcI_{3}$. Denote by $W,\ W_{3}$ and $\tilde f_{3}$ the corresponding vectors: $W_{3}$ is obtained from $f_{3}\in\mcI_{3}$ by means of the interwiner \eqref{principaltokirillov}. 
By \refs{RSlocal} and \eqref{usefulsoon}, we have
\begin{align*}L(x,x_{3},\tilde  x_{3})
&=\int_{K}\int_{k^\times}W(a(y)\kappa) W_{3}(a(y)\kappa)\tilde f_{3}(\kappa)|y|^{-1/2}\dti yd\kappa\\
&\ll C^{-1/2} \max_{\kappa\in\mcK_{C}}\int_{k^\times}|W(a(y)\kappa)W_{3}(a(y)\kappa)||y|^{-1/2}\dti y\\
&\ll C^{-1/2}\Sob^{\pi}(W)\max_{\kappa\in\mcK_{C}}\int_{k^\times}|W_{3}(a(y)\kappa)||y|^{-\theta}(1+|y|)^{-N}\dti y,
\end{align*}
where we applied Proposition \ref{zeroboundprop} for the last inequality. 
In the nonarchimedean case, it suffices to specialize to $\kappa =1$, 
since $W_3$ is in any case $\mcK_C$-invariant. 

 We now apply the bounds of the next Lemma to conclude.  Loosely speaking,
 it says that the Whittaker function $W_3$ is peaked near $y \sim C(\chi_{3})/C$
 and it takes a value of size $\sim 1$ there.  The reason for this
 can already be discerned from \eqref{principaltokirillov}: $W_3(y)$ is, approximately,
 the Fourier transform of a oscillating function restricted to an interval
 of the shape $|y| > C$.   The ``frequency'' of the function near the endpoints of this interval
 are approximately $C(\chi_{3})/C$. 
 \qed

 \begin{Lemmat} \label{uglyupperWhittakerbound}
Notations as above, for any $M,\eps> 0$, one has
 $$|W_3(a(y) \kappa)| \ll_{M,\eps,\psi,\pi_{2}}(C/|y|)^\eps\left( \frac{C}{C(\chi_{3})}\right)^{1/2 }|y|^{1/2} (1+\frac{|y|}{C(\chi_{3})/C})^{-M}.$$

 \end{Lemmat}
 
The proof of Lemma \ref{uglyupperWhittakerbound} is, 
 regrettably, an explicit calculation. 
 Let us note that 
$$\whitconst^{1/2} W_3(g) = |\det g|^{1/2}\int_{k^\times}\int_{k}\Psi_{3}((t,tx)g)\psi(x)|t|\chi_{3}^{-1}(t) dx\dti t,$$

where $\Psi$ is as in Definition \ref{X3def} in the archimedean case;
in particular $\Psi_1(x,y) = \phi_1(x) \phi_2(y)$, where $\phi_i$ 
are as in \eqref{phiiarch}; in the nonarchimedean case the same
is true with $\phi_i$ as in \eqref{phiinonarch}; 
 finally
$\whitconst$ is as in \eqref{principaltokirillov}, and is a harmless constant. 
In particular $\whitconst^{1/2} W_3(a(y)) $ may be expressed as
$$y^{1/2}  \int_{k^{\times}} \int_{k} \phi_1(tx) \phi_2(ty) \psi(x)|t|\chi_{3}^{-1}(t) dx\dti t.$$
\proof (of Lemma \ref{uglyupperWhittakerbound} -- nonarchimedean case.)  
It suffices to check $\kappa = 1$. 
\begin{align*}\whitconst^{1/2} W_{3}(a(y))&=\vol(K_{0}[f])^{-1/2}|y|^{1/2}\int_{k^\times}1_{\varpi^f\order}(yt)|t|\dti t\int_{k}\chi_{3}(x)1_{\order^\times}(tx)\psi(x) dx\\
&\asymp C^{1/2}|y|^{1/2}\int_{|t| |y| C \leq 1}|t| \dti t \int_{|x| |t| = 1}\chi_{3}(x)\psi(x) dx.
\end{align*}
If $\chi_3$ is ramified, this equals, in absolute value, $\asymp_{\psi} (C/C_{\chi})^{1/2} |y|^{1/2}$
when $|y|  \leq \frac{D C(\chi_3)}{C}$, where $D = q^{d_{\psi}}$, and $0$ otherwise (cf. \S \ref{gausssumcomputation}).
 On the other hand, 
if $\chi_3$ is unramified, it is zero for $|y| > q D/C$, and otherwise is bounded by
$\ll_{\psi} C^{1/2} |y|^{1/2} \log(1+ \frac{qD}{|y||C|})$.  \qed 
\proof 
(of Lemma \ref{uglyupperWhittakerbound} -- archimedean case.) 
We prove the bound first when $\kappa=Id$ as it is notationally more pleasant; we will comment on the necessary modifications afterwards.

We have \begin{equation} \label{question} W_{3}(a(y)) \asymp C^{1/2} |y|^{1/2}\int_{k^\times} \dti t \phi(C |t||y|) |t|\int_{k}\chi_{3}(x)\phi(|tx|-1)\psi(x) dx.\end{equation} 
Set $G(t)= |t| \int_{k}\chi_{3}(x)\phi(|tx|-1)\psi(x) dx=\chi_{3}^{-1}(t) \int_{k}\phi(|x|-1)\psi(x/t)\chi_{3}(x) dx $;  then
by Lemma \ref{almostgauss}, one has for any $N\geq 0$ and $\eps>0$
$$ |G(t)| \ll_{N,\eps,\phi} C(\chi_{3})^{-1/2+\eps}\min({C(\chi_{3})|t|},\frac{1+|t|^{-1}}{C(\chi_{3})})^{N}.$$

We deduce that $|W_3(a(y)|$ is bounded, up to an implicit constant, by
\begin{equation} \label{Completionofbound} C^{1/2} |y|^{1/2} \int_{k^{\times}} \phi(C |t| |y|) G(t)\dti t   \ll_{\eps,\pi_{2},N}   C^{\eps}\left( \frac{C}{C(\chi_{3}} \right)^{1/2} \frac{ |y|^{1/2-\eps}}{
 (1+\frac{|y|}{C(\chi_{3})/C})^{N}}.\end{equation} 
For general $\kappa=\left(\begin{array}{cc}a & b \\c & d\end{array}\right)\in\mcK_{C}$,
replace $\Psi_{3}$ by the function $\kappa.\Psi_{3}$ with
$$\kappa.\Psi_{3}:\ (x,y)\mapsto \Psi_{3}((x,y)\kappa)=\Psi_{3}(ax+cy,bx+dy)=\phi(C|ax+cy|)\phi(|bx+dy|-1).$$
In \refs{question} we then need to replace $\phi(C |t||y|) |t|\int_{k}\chi_{3}(x)\phi(|tx|-1)\psi(x) dx$ by 
\begin{multline*}|t|\int_{k}\phi(C|aty+ctx|) \phi(|bty+dtx|-1)\psi(x)\chi_{3}(x) dx\\=
\chi_{3}^{-1}(t)\int_{k}\phi(C |aty+cx|) \phi(|bty+dx|-1)\psi(x/t)\chi_{3}(x) dx\end{multline*}
Observe that since $|b|,|c|\ll C^{-1},\ |a|-1,|d|-1\ll C^{-1}$, the integral is zero unless $|ty|\ll C^{-1}$, while the $x$ variable satisfies
$\big||x|-1\big|\leq \delta_{0}+O(C^{-1}).$ 

The above proof carries on {\em mutatis mutandis}. At the first stage,
we compute the $x$-integral, noting that $x \mapsto \phi(C |aty+cx|) \phi(|bty+dx|-1)$
satisfies similar smoothness bounds to $\phi(|x|-1)$, since $c |C| \ll 1$.  Then \eqref{Completionofbound} uses in addition
only the fact that we may restrict to $|ty| \ll C^{-1}$, and so also goes through. 
\qed

\section{Integral representations of $L$-functions: global computations.}

We now turn to global aspects of the analysis of standard $L$-functions on $\GL(2)$ 
(\S \ref{PGL2}) and of the triple product $L$-function (\S \ref{triple}). 

\subsection{Notation} \label{notation-global}

\subsubsection{Subgroups}Henceforth $F$ is a number field, $\Aa:=\Aa_{F}$ and we often denote $|.|$ for $|.|_{\Aa}$. For $\G = \GL_{2}$ or $\PGL_{2}$,
we set $\bfX_{\bfG} = \G(F) \backslash \G(\adele)$.  We will often write simply $\mathbf{X}$
for $\mathbf{X}_{\PGL_2}$. 

For the various usual subgroups of $\G$ we adopt the notations and parametrizations of \S \ref{notationsubgroups}. We denote by $K$ the usual maximal open compact subgroup of $\bfG(\Aa)$.
  
We put $H:=\rmA(\Aa)=a(\Aa^\times)$ and $$H^{(1)} =\rmA(\Aa^{(1)})= \{ a(y): y \in \adele^{\times}: |y|_{\Aa} = 1\} \subset
H.$$ 
$H$ and $H^{(1)}$ are {\em closed, non-compact subgroups}
of $\GL_2(\adele)$. By an abuse of notation, we regard  them as closed, non-compact subgroups of $\PGL_2(\adele)$.  
 We let $Y$ be the $H$-orbit of the identity coset in $\bfX_{\PGL_{2}}$.  Then $Y$ carries a $H$-invariant measure of infinite volume. 

\subsubsection{Measures} \label{measures-global}

We adopt on each of the groups the product of those measures specified in 
\S \ref{measures-local}.

  Although not strictly necessary, it is a useful
check to be aware of the relation with a Tamagawa measure.
A routine computation shows that, for the measures 
on $\PGL_2$, \begin{equation} \label{tamag} dg =  (\disc F)^{-1} \xi_F(2) d^{\tau}g,\end{equation}
where $d^{\tau}g$ is the Tamagawa measure. 
In particular, the volume of the quotient $\bfX$ is $2 \xi_F(2) (\disc F)^{1/2}$.

\subsubsection{Additive characters and Fourier expansion} Let $\psi=e_{F}=\otimes_{v}\psi_{v}$ be the non-trivial additive character of $\Aa/F$ given as $e_{F}(.)=e_{\Qq}(\tr_{F/\Qq}(.))$ with $e_{\Qq}$
the unique additive character of $\Aa$ whose restriction to $\Rr$ coincide with $e^{2 \pi i x}$. 
 
 For $f(g)$ a function on $\bfX_{\GL_{2}}$, we denote its constant term $f_N$
 and Whittaker function $W_f$ by 
\begin{equation} \label{AW} f_{N}(g):=\int_{\Aa/F}f(n(x)g)dx, \ W_{f}(g)=\int_{\Aa/F}\psi(x)f(n(x)g)dx\end{equation} 
and one has the Fourier expansion
\begin{equation} \label{fourierexpansion} f(g)=f_{N}(g)+\sum_{y\in F^\times}W_{f}(a(y)g).\end{equation}

\subsubsection{Representations} 
In this section, $\pi$ will be a standard automorphic representation on $\GL(2)$, 
and $C(\pi)$ will be its analytic conductor, i.e. if $\pi = \otimes \pi_v$,
then $C(\pi) = \prod_{v} C(\pi_v)$, the right-hand quantities as in \S \ref{sec-localanalyticconductor}. 
If $\pi$ is such, then $\pi$ carries a canonical inner product, and we shall always regard
it as a unitary representation with respect to this inner product; we fix moreover
an inner product on each $\pi_v$ compatibly with this choice. 

Therefore, if we say, ``let $f =\otimes f_v \in \pi = \otimes \pi_v$,'' 
we always have $$\langle f, f \rangle = \prod_{v} \langle f_v, f_v \rangle.$$

\begin{Rem} 
Here is a note of warning in connection with this terminology. 
Suppose $\pi = \otimes \pi_v$ is a standard automorphic representation. 
The map $f \mapsto W_f$ intertwines $\pi$ with $\otimes_{v} \mathcal{W}(\pi_v)$. 
This map is {\em not} an isometry when $\pi$ is endowed with the canonical inner product,
and each Whittaker model is endowed with the norm \eqref{whitinner}.  See \S \ref{canGL2}. 
\end{Rem}

\subsubsection{Infinite product notation} \label{ipnot}
We denote by $\Lambda(\pi, s)$ resp. $L(\pi, s)$ the completed $L$-function of $\pi$, resp. 
the $L$-function omitting archimedean factors. 

Let us adopt the following notational conventions, to be held in force through the rest of the paper:
\begin{itemize}
\item [-]
First of all, if $\pi$ is an automorphic representation, we shall understand the notation $A \ll_{\pi} B$
to mean ``there exists absolute constants $a,b$ so that $|A| \leq a C(\pi)^b |B|$. 
\item [-] We continue with the already established convention: if we write $|\ell(f)| \ll \Sob(f)$, we mean
that there exists a constant $d$, depending only on $[F:\Q]$,  so that $|\ell(f) | \ll \Sob_d(f)$. 
\item [-]  (Regularizing products over places.)  Suppose that $\Lambda_1, \Lambda_2(s)$ are (completed) global $L$-function with local factors $L_{1v}(s), L_{2v}(s)$;
let $s_0$ be so that $L_{1v}(s_0) \neq 0$ and $L_{2v}(s_{0}) \neq 0$ for all $v$,
and so that $\Lambda_1$ is holomorphic at $s_0$; and suppose that $E_v$ is a function on places of $F$
with the property that $E_v =L_{1v}(s_0)/L_{2v}(s_0)$ for almost all $v$.  

Put $\Lambda^*(s_0) $ to be $\Lambda_1(s_0)/\Lambda_2^*(s_0)$, where
$\Lambda_2^*(s_0)$ is the first nonvanishing Laurent coefficient of $\Lambda_2$ at $s=s_0$. 
In particular, $\Lambda^*(s_0) =\Lambda_1(s_0)/\Lambda_2(s_0)$
if $\Lambda_2$ is holomorphic and nonvanishing at $s=s_0$; we define $L^{S,*}(s_{0})$ in the same way,
omitting the factors at $v \in S$. 
We shall then write:

$$\prod_v^* E_v \stackrel{\mathrm{def}}{=}  \Lambda^*(s_0) \prod_v \frac{E_v}{L_{1v}(s_0)/L_{2v}(s_0)}
= L^{S,*}(s)\prod_{v \in S }E_{v}\prod_{v \notin S}\frac{E_v}{L_{1v}(s_0)/L_{2v}(s_0)}.
$$
where the second equality holds for any finite set of places $S$; we shall often use
it with $S$ the set of all archimedean places. \end{itemize}

\subsubsection{Eisenstein series} \label{notation-global-Eis} Given two characters $\chi^+,\chi^-$ of $F^\times\bash \Aa^\times$ whose product is unitary, we denote by $\mcI(\chi^+,\chi^-)$ or $\chi^+\boxplus\chi^-$ the representation of $\GL_{2}(\Aa)$ unitarily induced from the corresponding representation on $\bfB(\Aa)$:
the $L^2$-space of functions $f$ on $\GL_{2}(\Aa)$ such that
$$f(\left(\begin{array}{cc}a & b \\0 & d\end{array}\right)g)=|a/d|_{\Aa}^{1/2}\chi^+(a)\chi^-(d)f(g), \ \ \ \   \peter{f,f}=\int_{K} |f(k)|^2 dk.$$

Given $f$ in such a space, we denote by $E(f)=\Eis(f)$ the corresponding Eisenstein series (defined by analytic continuation, in general), i.e. $$\Eis(f) := \sum_{\rmB(F) \backslash \GL_2(F)} f(\gamma g).$$

For $s$ a complex parameter and $\pi=\mcI(\chi^+,\chi^-)$ we set 
\begin{equation} \label{pidefdef} \pi_s = \pi(s)=\mcI(\chi^+|.|^s_{\Aa},\chi^-|.|^{-s}_{\Aa})\end{equation} and for $f\in\pi$ we define $f_s = f(s)$ to be the unique function in $\pi(s)$ whose restriction to $K$ coincide with $f$.

\subsubsection{Fourier coefficients of Eisenstein series} \label{Eis-fourier} 
Given an Eisenstein series, $\Eis(f)$, its constant term is given by:
\begin{equation} \label{eisenconstant} \Eis(f)_N(g) = f(g) + \int_{ \adele} f(w n(x) g)dx \end{equation} 

The Eisenstein series has a pole at the point $\chi^{+}/\chi^{-} = |\cdot|_{\Aa}$ coming from its constant term; the residue, with respect to the coordinate $s$ as in \eqref{pidefdef}, is given by
\begin{equation} \label{rescomp} \prod_{v}^* \int f_v(w n(x)) dx  =  \frac{1}{2} (\mathrm{disc} F)^{-1/2} \frac{ \xi^*_F(1)}{\xi_{F}(2)} \int_{k \in K} f(k);\end{equation}
with $\xi_F^*(1):=\mathrm{res}_{s=1} \xi_F(s)$;
see \eqref{ntok} and discussion thereafter.

The other Fourier coefficients of Eisenstein series are determined by the Whittaker function:
\begin{align} \label{whit-eis} W_{\Eis(f)}(g)&=\int_{F\bash \Aa}f(wn(x)g)\psi(x)dx =  \prod_{v} W_{f,v},\\
W_{f,v}(g)&= \int_{x \in F_v} f_v(w n(x)g) \psi_v(x) dx.\nonumber\end{align}
$$W_{f,v}(a(y))= |y|^{1/2} \chi^{-}(y) \int_{x \in F_v} f_v(w n(x)) \psi_v(xy) dx.$$
 
 Let us recall that in \S \ref{Eisenstein}, we have discussed two possible norms on Eisenstein series: the Eisenstein norm
 $\peter{ \Eis(f), \Eis(f)}_{Eis}=\peter{f,f}$ and the canonical norm
$\peter{ \Eis(f), \Eis(f)}_{can}$ formed out of the inner products of the local Whittaker functions $W_{f,v}$ (cf. \refs{canreg} and \refs{globalinner}). In view of \S \ref{ps-local} and the remark in \S \ref{meascomp1}, we have the following relation -- if $\Lambda(\pi, \Ad, 1)$ 
has a simple pole -- (recall that ${\whitconst}_{v}=q_{v}^{-d_{\psi_{v}}/2}\zeta_{v}(1)^2/\zeta_{v}(2)$)
\begin{align} \nonumber \langle \Eis(f), \Eis(f) \rangle_{can} &=   \Lambda^*(\pi, \Ad, 1) \frac{ 2 \xi_F(2) (\disc F)^{1/2} }{\xi_F^*(1)} \prod_{v}
\frac{ \langle f_v, f_v \rangle {\whitconst}_{v}  }{\zeta_v(1) L_v(\pi, \Ad, 1) / \zeta_v(2)}  \\ & = 
2 \xi_F(2) \langle f, f \rangle = 2 \xi_F(2) \langle \mathrm{Eis}(f), \mathrm{Eis}(f) \rangle_{Eis}.\label{eis-norm-eqn}
\end{align}

\subsubsection{The intertwiner operator} \label{intertwining} 

The map $M: f \mapsto \int_{x \in \adele} f(w n(x) g)dx$ that occurs in \eqref{eisenconstant}
is the {\em standard intertwining operator}.  It is an isometry on the ``unitary axis,'' i.e.
when the characters $\chi^{\pm}$ are unitary. We shall need at several points to study its analytic behavior off the analytic axis. 

$$ M = \frac{\Lambda(\chi^{+}/\chi^{-},0)}{\epsilon(\chi^+/\chi^-, 0) \Lambda(\chi^{+}/\chi^{-},1)} \prod \bar{M}_v, 
\ \ \bar{M}_v := \frac{\epsilon(\chi_v^{+}/\chi_v^{-}, \psi_v, 0) L(\chi_v^{+}/\chi_v^{-},1)}{ L(\chi_v^{+}/\chi_v^{-},0)} M_v.$$
Here $M_v: \mathcal{I}(\chi_v^+, \chi_v^{-}) \mapsto  \mathcal{I}(\chi_v^-, \chi_v^{+}), f \mapsto \int_{x \in F_v} f(w n(x) g)dx$; the operator $\bar{M}_v$ has the advantage of being holomorphic
in $\chi_v^+, \chi_v^{-}$.    In view of the functional equation, the global correction factor
$ \frac{\Lambda(\chi^{+}/\chi^{-},0)}{\epsilon(\chi^+/\chi^-, 1) \Lambda(\chi^{+}/\chi^{-},1)} $
has absolute value $1$ (it is the ``scattering matrix'') although it may not be equal to $1$. 

 By \S \ref{meascomp1},  $\bar{M}_v$
takes the spherical vector with value $1$ on $K$, to the spherical vector with value $1$ on $K$,
for almost all $v$. More generally, it preserves norms
up to a scalar depending only on $\psi$: \cite[Section 4]{GJPSP}. 

\begin{Lemma*}
Let $f_v \in \mathcal{I}(\chi_v^+, \chi_v^{-}) $.
Suppose that $|s| \leq \delta_0$,  and the deformation $f_{v,s}   \  (s \in \C)$ is as in \S \ref{sss:deform}. 
Then there exists $d$ so that:
\begin{equation} \label{intertwining-bounds-local} \sup_{k \in K} |\bar{M}_v f_{v,s}(k)| \ll \Sob^{\mcI}_d(f), \int_{k \in K} |\bar{M}_v f_{v,s}(k)|^2 \ll \Sob_{d \delta_0}(f_v)^2.\end{equation}
\end{Lemma*}

\proof 
Because of  \eqref{wfmellin}, 
and with notation as contained there, 
$\bar{M_v}f(1)$ is proportional to 
$\bar{\ell}^{\chi'} (W_f)$, where 
$\chi' = \alpha^{-1/2}/\chi^{-}$, and 
where $\bar{\ell}^{\chi}$ is the ``normalized'' functional $\frac{ \int W(y) \chi(y) d^{\times}y}{L(1/2, \pi \otimes \chi)}$, for $W \in \Whit(\pi, \psi)$. More precisely,
$\bar{M_v}f(1) = \horrid \bar{\ell}^{\chi'}(f)$, where
$$\horrid = \zeta_k(1)^{-1}  \whitconst^{1/2}  
\frac{L_v(1, \chi_v^{+}/\chi_v^{-}) \epsilon(0, \psi_v, \chi_v^{+}/\chi_v^{-})\zeta_v(1)}{ \epsilon(\alpha, \psi, 0)}.$$
Applying the local functional equation shows 
that $\bar{M_v}f_v(1)$ is given by:
$$\bar{\ell}^{\chi''}(W_f) \frac{\horrid}{\epsilon(\pi, \chi', 1/2)}, \chi'' = \alpha^{1/2}/\chi^{+}.$$

The local bound $\sup_{k \in K_v} |\bar{M_v} f_s(k)| \ll \Sob^{\mcI}_{d_0}(f)$,
now follows from  Proposition \ref{zeroboundprop} together with the following observation:
the smallest eigenvalue of the local Laplacian $\Delta_v$ on $ \mathcal{I}(\chi_v^+, \chi_v^{-}) $ is bounded
below by a positive power of $ (\Cond_v(\chi^{+}) + \Cond_v(\chi^{-}))$. This (a simpler form of
Lemma \ref{gan}) allows us to absorb dependencies on $\chi^{\pm}$ into the Sobolev norm. 

To obtain
the stated bound for $\int_{k \in K} |\bar{M}_v f_{v,s}(k)|^2 $, we
interpolate, taking into account the equality:
$\| \bar{M}_v f_{v,s}\|_{L^2(K_v)} \asymp_{\psi} \|f_{v}\|_{L^2(K_v)}$ whenever $\Re(s) = 0$. 
This argument is analogous to the interpolation in the proof of Proposition \ref{zeroboundprop-def}. 
\qed

We now give a global analogue of this statement:

\begin{Lemmat}
Let $f \in \mcI(\chi^{+}, \chi^{-})$. 
Suppose that $|s| \leq \delta_0$,  and the deformation $f_s   \  (s \in \C)$ is as in \eqref{pidefdef}. 
Then there exists $d$ so that \begin{equation} \label{intertwining-bounds} \sup_{k \in K} |M f_s(k)| \ll \Sob^{\mcI}_d(f), \int_{k \in K} |Mf_s(k)|^2 \ll \Sob^{\mcI}_{d \delta_0}(f)^2.\end{equation}
\end{Lemmat}
\proof 
An application of (S1d) to the previous Lemma gives the first statement.   From this, it follows that
$\int_{k \in K} |M f_s(k)|^2 dk$ is bounded by $\Sob^{\mcI}_{d_0}(f)$. 
In order to deduce the second statement, we observe that
$\int_{k \in K} |M f_s(k)|^2 dk  = \|f\|_{\mcI}^2$ when $\Re(s) = 0$, 
and then interpolate, as in the proof of Proposition \ref{zeroboundprop-def}. 
\qed 
\subsubsection{The Eisenstein series at a singular parameter}\label{singulareisenstein}

It will also be of interest to consider the Eisenstein series associated to $\chi^{+} = \chi^{-} = 1$. 
At this point, the map $f \mapsto (g \mapsto \int_{x \in \adele} f(w n(x) g) dx)$
is the negative of the identity map, and accordingly the Eisenstein intertwiner vanishes to order one.  
Write $\Eisreg$ for the (normalized) intertwiner
from $\mcI(1, 1)$ to $C^{\infty}(\bfX)$, given by:
$$f\mapsto \xi_{F}(1+2s)\mathrm{Eis}(f_s)|_{s=0}.$$
For $f=\prod_v{f_{v}}$ factorizable, let us define $\peter{\Eisreg(f),\Eisreg(f)}_{can}$ via \refs{canreg} and \refs{globalinner}.

Then, if $f \in \mathcal{I}(1,1)$, then
the constant term
\begin{multline} \label{eis-singular} \Eisreg(f) _N (a(y) g)   \\ =
- \frac{1}2\xi_{F}^*(1)
\frac{d}{ds}|_{s=0}  \left( \int_{x \in \adele} |y|^{1/2 +s} f_s(w n(x) g )dx  +  |y|^{1/2-s} f_s(g) \right) \\ = 
\frac{1}2\xi_{F}^*(1)(2 f(g)   |y|^{1/2} \log y - M^* f (g) |y|^{1/2})
\end{multline}
where $M^* f$ is the derivative, at $s=0$, of $f_s + \int_{x \in \adele} f_s (w n(x) g)dx$;
as usual, this needs to be interpreted by a process of analytic continuation.

Trivial bounds on $M^*$ follow from bounds on intertwining operators in \S \ref{intertwining} (in particular,
one may bound the derivative of an analytic function by bounding its value on a small circle around
the point of interest). This shows that, 
 for every $\varphi $ belonging to the space of $\Eisreg(\mathcal{I}(1,1))$, we have
\begin{equation} \label{eis-singular-2} \varphi_N(a(y) g) \ll |y|^{1/2} \log|y| \Sob^{\pi}(\varphi),\end{equation} 
Now take $f =f^0$ to be the spherical vector in $\mcI(1,1)$, normalized so that $f^0 | K=1$.
The function $M^* f^0$ is continuous, and thus bounded on compact sets;  it follows from \eqref{eis-singular} that
there exists $X_0$ so that
\begin{equation} \label{eis-sing-lowerbound}
\Eisreg(f^0) (x)  \gg \height(x) \log \height(x), \ \ \height(x) \geq X_0. \end{equation}

\subsection{Hecke integrals on $\PGL_2$}\label{PGL2}
In this section we shall study the Hecke integral for the standard $L$-function on $\GL_2$,
and we shall give some bounds on the Hecke integral of a translate of a given vector,
using the method described in \S \ref{MC0}, Remark.

\subsubsection{The Hecke-Jacquet-Langlands integral}  \label{sec:hjl} 
We recall in this section the integral representation
for the standard $L$-function on $\GL_2$, following Jacquet and Langlands. 
 Let $\chi$ be a character of $\Aa^\times/F^\times$;
the integral
\be\label{JLint}\ell^\chi(\varphi):=\int_{\adele^{\times}/F^\times} (\varphi - \varphi_N)(a(y))\chi(y)  d^{\times} y
\ee
defines a functional on the space of any standard generic automorphic representation $\pi$. Indeed,
it is absolutely convergent if $\pi$ is cuspidal, and, in general, can be interpreted by analytic continuation 
in the $\chi$ variable.

It was observed by Hecke, and generalized by Jacquet and Langlands, that
the period $\ell^{\chi}$ is very closely related to the standard $L$-function:
if $\varphi = \otimes \varphi_v$ is a pure tensor in $\pi$,
and we factorize the associated Whittaker function (see \eqref{AW}) as $W_\varphi = \prod_{v} W_{\varphi,v}$
, one has:
\be\label{JLperiod}\ell^\chi(\varphi)=\prod_v^* \ell^{\chi_v}(W_{\varphi,v}) \left( =\Lambda(\pi\otimes\chi,1/2)\prod_v \frac{\ell^{\chi_v}(W_{\varphi,v})}{L(\pi_{v}\otimes\chi_{v},1/2)}\right), 
\ee
for $\ell^{\chi_v}(W_{\varphi,v})$ the local Hecke functionals defined in \S \ref{ugly}.
The verification of \eqref{JLperiod} is ``unfolding.'' It is valid for $\pi$ Eisenstein, even in the singular case.

We will use also use  the following ``canonically normalized'' expression, given in terms of matrix coefficients and which follows from \S \ref{matrixcoeff} and \S \ref{canGL2}) (cf. the Remark in \S \ref{MC0}):
 for $\chi$ a unitary character, and $\pi$ standard, 
\be\label{periodsquare}\frac{ |\ell^\chi(\varphi)|^2}{\langle \varphi, \varphi \rangle_{can}} =  
\left( \whitintertwinerconstant \right)^{-1} \prod_{v}^* \frac{h^{\chi}(\varphi_v)}{\langle \varphi_v, \varphi_v \rangle}
\ee
where almost all local factors in the last product are equal to 
$\frac{|L_{v}(\pi\otimes\chi,1/2)|^2}{\zeta_{v}(1)L_{v}(\pi, \Ad, 1)/\zeta_{v}(2)}$. 
Again, this formula is valid for $\pi$ Eisenstein, even at singular points;
of course, at singular points, the regularization implicit in $\prod^*$ involves
taking a higher residue of $L(s, \pi, \Ad)$ at $s=1$. 

It should be noted that the right-hand sides of \eqref{periodsquare} and \eqref{JLperiod}
make sense for all $\chi$, even when the left-hand sides must be interpreted by analytic continuation. 
\subsubsection{Bounds for the Hecke integral of the translate of a function} \label{sec:udtt-cor}
\begin{Lemma*} \label{udtt-cor} 
Let $\pi$ be a generic automorphic representation of $\PGL_{2}$. For $\varphi \in \pi$, 
any $g \in \PGL_2(\adele)$, and unitary $\chi$: there is an absolute constant $d$ so that
$$ |  \ell^{\chi}(g. \varphi) |
\ll_{\pi,\chi} \disc(Y g)^{-\frac{1-2\theta}{32}} \Sob_{d}^{\pi}(\varphi),$$
  Here $\disc(Yg)$ is the discriminant of the adelic torus orbit $Yg$ 
  as defined in \cite{ELMV3} -- see below --  and the unitary structure on $\pi$ is the canonical norm. 
\end{Lemma*}

The paper \cite[4.1--4.2]{ELMV3} attaches to the adelic torus orbit $Yg$
local and global discriminants, denoted, respectively, as $\disc_v(Yg)$ and $\disc(Yg)$.  Although most of that paper deals with the case of adelic points of anisotropic tori, the definition is perfectly applicable to the split adelic torus $H$. 
 For our purposes, it is enough to know
that, for the special case $g_{0} = n(T)$, $|T|_{\Aa}\geq 1$ we have:
\begin{equation} \label{wtf} \disc_v (Yg)  \asymp \max(1,|T|^2_v), \disc(Yg) \asymp \prod_{v}\max(1,|T|^2_v)\geq |T|^2_{\Aa}.\end{equation}

\begin{Rem} The exponent $-\frac{1-2\theta}{32}$ is certainly not best possible: it is taken from the general computations of {\em loc. cit.}, which are in no way optimal for $\PGL_{2}$ (see for instance \cite{CU} for better bounds of similar integrals.) What is important to us, in this paper, is that this exponent is negative. 
\end{Rem}

\proof 
By  \refs{periodsquare} and \eqref{localglobalbounds}, we have:
\begin{gather*} |\ell^{\chi_{v}}(g_v \cdot \varphi_{v})|^2=\int_{F_v^{\times}} \langle a(y)g. \varphi_v , g. \varphi_v  \rangle  \chi_v(y) d^\times y=  \int_{F_v^{\times}} \langle g^{-1} a(y)g. \varphi_v , \varphi_v  \rangle \chi_v(y)   d^\times y\\
\ll \Sob^{\pi_v}(\varphi_v)^2\int_{F_{v}^\times}\Xi_{v}( g^{-1}a(y)g)^{1-2\theta}d^\times y\\
\ll \disc_v(Y g)^{-\frac{1-2\theta}{16}}  \Sob^{\pi_v}(\varphi_v)^2.
\end{gather*}

The bound in the last step is carried out, in a more general setting in  \cite[Lemma 9.14]{ELMV3}, and
we do not reproduce it here.  

Write $L_{v}  = L(\pi_{v}\otimes\chi_{v},1/2)/L(\pi_v, \Ad, 1)^{1/2}$. 
We are going to apply Property (1d) of Sobolev norms to $\prod \ell'_v$, where $\ell'_v$ is
the ``normalized'' functional on $\pi_v$ so that $|\ell'_v|^2 = \frac{|\ell^{\chi_v}|^2}{|L_v|^2}$. 
It enjoys the following properties:

\begin{enumerate}
\item 
For any place for which $\varphi_{v}\otimes\chi_{v}$ 
is spherical, $|\ell'(\varphi_v)|^2 = \langle \varphi_v, \varphi_v \rangle = \mathcal{S}^{\pi_v}_{0}(\varphi_v)^2$. 
\item  There exists $A, d_0$ so that, for any $v$, the operator norm
of $\ell'$ w.r.t. $\Sob^{\pi_v}_{d_0}$ is $\leq A |L_v|^{-1} \disc_v(Y g)^{-\frac{1-2\theta}{32}}.$
\item 
There exists an absolute constant $C$ so that, for nonarchimedean $v$, 
$|L_v|^{-1} \leq C$. 
\end{enumerate}

Thus property (S1d) applied to $\ell' := \prod_{v} \ell'_v$
shows that, for $\varphi \in \pi$, 
$$\frac{ |\ell'(\varphi)|^2}{\Sob^{\pi}_{d'}(\varphi)^2}  \ll_{\varepsilon} \disc(Yg)^{-\frac{1-2\theta}{16}+\eps} \prod_{v | \infty} |L_v|^{-2}
,$$ for some $d > 0$. 
Applying trivial or convexity bounds
for all local factors or $L$-functions,  we arrive at:
$$ \frac{ |\ell^{\chi}(\varphi)|^2 }{\Sob^{\pi}(\varphi)^2} \ll_{\pi, \chi, \eps} \disc(Yg)^{-\frac{1-2\theta}{16}},$$
as required.  \qed

This Lemma admits the following mild refinement where we improve
the $\chi$-dependence at the implicit cost of weakening the $\pi$-dependence: 

\begin{Lemmat} \label{udtt-cor-refined} 
Notation as in the prior Lemma, for any $N\geq 0$, there is $d=d(N)$ so that
$$ |  \ell^{\chi}(g. \varphi) |
\ll_{\pi} (\Cond(\chi))^{-N} \disc(Y g)^{-\frac{1-2\theta}{32}} \Sob_{d}^{\pi}(\varphi),$$
\end{Lemmat}

\proof
By \eqref{rt1}, it suffices to prove such an assertion for $\varphi \in \pi[\underline{m}]$.
(One also verifies that the implicit dependence of $A'$
on $A$ in \eqref{rt1} is independent of $\pi$).  
The assertion for $\varphi \in \pi[\underline{m}]$ follows by integration by parts.
\qed

\subsection{Regularization}\label{regularsection}

In this section, we define a regularization process to define integrals of non-necessarily decaying functions on $\bfX_{\PGL_{2}}=\PGL_{2}(F)\bash\PGL_{2}(\Aa)$.  Such a regularization was 
given by Zagier \cite{Zagier}; for our purposes, an alternate way of defining it using convolution will be
of use. 

\subsubsectionind{} \label{Rreg}
For motivational purposes, we consider first a toy example, namely, the corresponding situation on $\mathbb{R}_{>0}$ 
(one could even consider the case of the integers).  We shall regard $\mathbb{R}_{>0}$
as a multiplicative group in what follows. 

A {\em finite} function on $\mathbb{R}_{>0}$ is one whose
translates, under the action of the multiplicative translations:
$$\tau_y f (x) = f(xy)$$ span
a finite dimensional representation; equivalently, it is a linear
combination of functions $x^\alpha (\log x)^{b}$, for $\alpha \in \C$ and $b \in \mathbb{N}$. 

We call a finite function {\em admissible} if the exponent $\alpha = 0$ never occurs. A more
intrinsic and useful way of formulating this is: $f$ is admissible if the span $\langle \tau_{y} f \rangle$
of all multiplicative translates
does not contain the trivial representation of $\mathbb{R}_{>0}$. 

Let $V_+$ be the space of continuous functions on $\mathbb{R}_{>0}$ of rapid decay as $x \rightarrow \infty$ (i.e. $|f(x)| \ll_N |x|^{-N}$ for all $N$); let $V_{-}$ be the space of rapid decay as $x \rightarrow 0$;
and let $V_0 = V_{+} \cap V_{-}$. 
Let $V$ be the space of all continuous functions $f$ on $\mathbb{R}_{>0}$
so that there exists {\em admissible} finite functions $f_1, f_2$ so that $f-f_1 \in V_{+}, f-f_2 \in V_{-}$. 
Clearly $V_0 \subset V$.  The following Lemma is well-known:
 
\begin{Lemmat}
There exists a unique functional (the ``regularized integral'') on $V$ which extends
integration $f \mapsto \int f(x) \frac{dx}{x}$ on $V_0$, and is invariant under  multiplicative translation. 
\end{Lemmat}

This functional may be defined in multiple ways:
\begin{enumerate}
\item 
Given $f \in V$, we may find finite collections $y_i \in \mathbb{R}_{>0}$, $c_i \in \mathbb{C}$
so that $f' := \sum c_i \tau_{y_i} f \in V_0$, and so that $\sum c_i \neq 0$.  We then define the regularized integral of $f$
to be $\frac{ \int f'(x) \frac{dx}{x}}{\sum c_i}$. 
\item Consider $\int_{1/T}^{T} f(x) \frac{dx}{x}$; this has the form $g(T) + h(T)$, 
where $g(T)$ is an admissible finite function, and $h(T)$ has a limit as $T \rightarrow \infty$. 
We define the integral to be $\lim_{T \rightarrow \infty} h(T)$. 
\item Consider $F_{+}(w) = \int_{1}^{\infty} f(x) x^w \frac{dx}{x}$; 
define similarly $F_{-}(w)$. The functions $F_{\pm}(w)$ are convergent
for $\mp \Re(w) \gg 1$; they extend to meromorphic functions on the plane, and
are holomorphic at $w=0$. We define the integral to be $F_{+}(0) + F_{-}(0)$. 
 \end{enumerate}

The regularized integral of any admissible finite function is zero. This follows, without computation,
because of invariance under multiplicative translation.
Thus, if there exists finite function $f_0$ so that $f - f_0  \in V_0$, then
  the regularized integral of $f$ equals $\int (f-f_0)(x) \frac{dx}{x}$.

  The remarks of this section adapt without change to define a regularized 
  integral on $\adele^{\times}/F^{\times}$.  

\subsubsectionind{}
Given a function $\varPhi$ on $\GL_{2}(F)\bash\GL_{2}(\Aa)$, with unitary central character $\chi$ (i.e. which transform by $\chi$ under translation by $\rmZ(\Aa)$), we say that $\varPhi$ is of controlled increase if there exists a function $$f: \rmN(\adele) \rmA(F) \backslash \GL_2(\adele) \rightarrow \C,$$ spanning
a finite-dimensional space under the translation action of $\GL_2(\adele)$, and with central character $\chi$, so that, 
for every $N \geq 0$

\be\label{decay-cond} \varPhi\left( \left( \begin{array}{cc} 1 & x \\ 0 & 1 \end{array}\right)\left(\begin{array} {cc}y_1 & 0 \\0 & 1 \end{array}\right) k  \right) = f \left( \left(\begin{array} {cc}y_1 & 0 \\0 & 1 \end{array}\right) k  \right) + O(|{y_1}|^{-N}) \mbox{ as $|{y_1}| \rightarrow \infty$.}\ee
In other terms, the difference is of rapid decay.  In more explicit terms,
there must exist a finite collection of functions $\chi_{i} : \adele^{\times}/F^{\times} \rightarrow \mathbb{C} \ (i \in I)$, each finite under the left translation action of $\adele^{\times}/F^{\times}$, 
as well as a corresponding collection of $K$-finite functions $\mathscr{K}_i: K \rightarrow \mathbb{C}$,
so that the left-hand side is well-approximated by $\sum \chi_i(y_1) \mathscr{K}_i(k)$. 

{\em A basic example to bear in mind is any sum or product of Eisenstein series.}

The expression $f$ is uniquely determined. We denote it by $\varPhi_N^*$; it need not coincide with the true
constant term of $\varPhi$. 
The set of {\em exponents} of $\varPhi$ (or of $\varPhi_N^*$) -- denoted $S_{\varPhi}$  -- 
is the set of characters of $\adele^{\times}/F^{\times}$  which
are (generalized) eigenvalues for the translation action of $\adele^{\times}/F^{\times}$
on the space spanned by $\varPhi_N^*$ and its translates. Of course if $\varPhi_N^*=0$ ($\varPhi$ is of rapid decay)
we set $S_{\varPhi}=\emptyset$.

\begin{Exam} If $\chi^+\not=\chi^-$, the exponents of a (non-zero) Eisenstein series $E\in\Eis(\mcI(\chi^+,\chi^-))$ are $\{\chi^+|.|^{1/2},\chi^-|.|^{1/2}\}$. If $\chi^+=\chi^-=\chi$ the same is true for $E\in\Eisreg(\mcI(\chi,\chi))$; the exponent $\chi|.|^{1/2}$ is now a generalized eigenvalue, i.e., has multiplicity. 
 If $E_{s}$ denote the Eisenstein series on the modular group with eigenvalue $1/4-s^2$, the exponents of $E_{s_{1}} E_{s_{2}} E_{s_{3}}$
is the set of characters $|\cdot|^{3/2+s}$, where $s = \pm s_{1} \pm s_{2} \pm s_{3}$. 
\end{Exam}

In the sequel, we will identify the complex numbers $\Cc$ with a subset of 
$\widehat{\adele^{\times}/F^{\times}}$ via $z\mapsto |.|^z$. For instance, given
 $S\subset \widehat{\adele^{\times}/F^{\times}}$, 
we say that ``$1 \in S$'' when $S$ contains the character $x \mapsto |x|$. 
We will also use an additive notation: given two such subsets $S_{1},\ S_{2}$ we denote by $S_{1}+S_{2}$ the set of pairwise products of the characters of $S_{1}$ and $S_{2}$. Of course $S_{1}+\emptyset=\emptyset$.

The operation $\varPhi \mapsto \varPhi_N^*$ is multiplicative: given two functions $\varPhi_1, \varPhi_2$ with exponents $S_1, S_2$, $(\varPhi_1 \varPhi_2)^*=\varPhi_{1,N}^*\varPhi_{2,N}^*$ and 
 $\varPhi_1 \varPhi_2$ has exponents in $S_1 + S_2$. 
 Finally, the complex conjugate $\overline{\varPhi}$ has exponents in $\bar{S}$ (i.e.
 the set of conjugates of those characters in $S$). 

The set of characters whose real part is $1/2$ (the real part being defined by $|\chi(.)|=|.|^{\Re(\chi)}_{\Aa}$) will be called the {\em unitary axis}: this is the set of exponents of the automorphic Eisenstein series.

\subsubsection{Regularized integral and regularized innerproduct}
Let $\Space_S$ be the vector space of smooth functions that are of controlled increase with trivial central character 
and whose exponents belong to $S$; and $\Space$ the union of $\Space_S$, where $S$ is
taken through all finite subsets of characters that do not contain any character of square $|\cdot|^2$
(i.e., any quadratic twist of $|\cdot|$). 

\begin{Lemma*}
There's a unique $\PGL_{2}(\Aa)$-invariant functional on $\Space$
extending integration on $L^1(\bfX_{\PGL_{2}})$. 
\end{Lemma*}

\proof One definition, due to Zagier, is given as follows: let $\mathbf{E}=\Eis(\varPhi^*_{N})$ be the Eisenstein
series induced from all exponents of $\varPhi$ that are of real part $> 1/2$ (or suitable derivatives thereof) and define
$$\int^{reg}_{\bfX_{\PGL_{2}}} \varPhi := \int_{\bfX_{\PGL_{2}}} (\varPhi - \mathbf{E}).$$ 
The right hand side makes sense, for $\varPhi - \mathbf{E}$ lies in $L^1$. Since it is clear on representation-theoretic grounds that the regularized integral of $\mathbf{E}$ must be zero (because the exponents of $\varPhi$ do not contain $1$), the uniqueness follows. \qed

As a corollary of the previous Lemma, we can define 
the {\em regularized inner product} for  $\varPhi_{1},\varPhi_{2}$ of controlled increase with the same central character and such that
$1 \notin S_1 + \overline{S_2}$:
$$\langle \varPhi_1, \varPhi_2 \rangle_{reg} = \int^{reg}_{\bfX_{\PGL_{2}}} \varPhi_1 \overline{\varPhi_2}.$$

\subsubsection{Regularization via convolution with measures} Here is an alternate definition that will
be, in fact, more suited for our purposes. (It also works better in higher rank.)

For every place $v$, we may choose a compactly supported measure $\mu_v$ on $\PGL_{2}(F_v)$ with the property that, for $\varPhi \in \Space_S$,
$\varPhi \star \mu_v \in L^1$. If $\int \mu_v  \neq 0$, then $\varPhi \mapsto\frac{ \int (\varPhi \star \mu_v)}{\int \mu_v}$ defines
a functional such as in the Lemma:   this functional is independent of $v$ or $\mu_v$, as we see by choosing a different place $w$ and
measure $\mu'_w$, and noting that $\star\mu_v$ and $\star \mu_w'$ commute. This functional is
 $\prod_{v' \neq v} \PGL_{2}(F_v)$-invariant; since $v$ is arbitrary, it is $\PGL_2(\adele)$-invariant. 
 Other definitions utilize truncation or related ideas; the disadvantage
of these is that the $\PGL_{2}(\adele)$-invariance is less clear.

\subsubsection{Explicit choice of a regularizing measure} \label{tripleexample}   The following
special case will occur:

Let $v$ be a finite place with a residue field of size $q_v$. 
Suppose that $\varPhi$ is of controlled increase and  unramified at $v$, and all $\chi \in S_{\varPhi}$ have real part $1$. 
Suppose that $\varphi$ is a cusp form unramified at $v$. 

Then one may choose a $K_v$-bi-invariant (signed) measure $\mu_v$ so that $\varPhi \star \mu_v$ has rapid decay, so that the total mass of $|\mu_v|$ is at most $4^{|S_{\varPhi}|}$, and 
so that $\varphi \star \mu_v = \lambda \varphi$, $|\lambda| \geq (1 - \frac{2 q_v^{\theta}}{q_v+1})^{|S_{\varPhi}|}$.
In other words, we may ``kill the growth of $\varPhi$ whilst only wounding $\varphi$;''
it should be noted that here the set $|S_{\varPhi}|$ is counted ``with multiplicity.''

 Indeed, there exists a finite function $f = \sum_{\chi \in S_{\varPhi}} f_{\chi}$
 on $\rmN(\adele) \rmA(F) \backslash \GL_2(\adele)$ 
 so that $\varPhi - f$ is rapidly decaying, as in \eqref{decay-cond}. 
The standard Hecke operator
$T_v := 1_{K_{v}a(\varpi_{v})K_{v}}$ acts on each $f_{\chi}$
by a (generalized) eigenvalue $\lambda_v(\chi)$, satisfying $q_v+1 \geq |\lambda_v| \geq q_v - 1$ in absolute value;
on the other hand, $T_v$ acts on $\varphi$ by an eigenvalue that is at most $2 q_v^{\theta}$. 

The measure  $(1 - \frac{T_v}{\lambda_v(\chi)})$ therefore 
annihilates $f_{\chi}$, and has total mass  $(1 + \frac{q_v+1}{|\lambda_v(\chi)|} ) \in [2,4]$. 
It acts on $\varphi$ by an eigenvalue that is $\geq 1 - \frac{2 q_v^{\theta}}{q_v+1}$. 
Take $\mu_v$ to be the convolution of these measures, for all $\chi \in S_{\varPhi}$.

\subsubsection{Regularized IP formula}  \label{regularIP}
A simple form of the Plancherel formula is the following: if $\varPhi_{1},\varPhi_{2}$ are functions on $\bfX_{\PGL_{2}}$
with rapid decay then
$$\peter{\varPhi_{1},\varPhi_{2}}=\int_{\pi\in\widehat{\PGL_{2}}^{Aut}}\peter{\Proj_{\pi}\varPhi_{1},\Proj_{\pi}\varPhi_{2}}d\mu_{P}(\pi)$$
where $\Pi_{\pi}$ denote the orthogonal projection on the space of $\pi$.

We describe now a version of that formula for functions $\varPhi_{1},\varPhi_{2}$ on $\bfX_{\PGL_{2}}$ of controlled increase.
Let $\pi$ be a standard automorphic representation of $\PGL_{2}$, and $\varPhi \in \mathcal{V}_S$.  If $S$ does not intersect the unitary axis and $\pi$ is {\em generic} (i.e. not equal to a quadratic character) then, for $\mcB(\pi)$ an orthonormal basis of $\pi$, we set
\be\label{projreg}\Pi_{\pi} \varPhi  = \sum_{\varphi \in \mcB(\pi)} \langle \varPhi, \varphi \rangle_{reg} \varphi \in \pi. \ee
Likewise, we define similarly $\Pi_{\pi}\varPhi$  for any $\pi$ which is {\em nongeneric} (i.e. one-dimensional), so long as $S$ does not contain any exponent whose square is $|.|^2$.

\begin{Proposition*}
Given $\varPhi_1$ and $\varPhi_2$  of controlled increase with exponents of real part $>1/2$. Let $S_1, S_2$ denote the respective sets of exponents . Suppose that $S_{1},\ S_{2}$  are disjoint and that $S_{1}\cup S_{2}\cup S_1 +\overline{S_2}$ does not contain any character whose square is $|.|^2$, then:
\begin{equation}\label{RIP}   \langle \varPhi_1, \varPhi_2 \rangle_{reg}  = \int_{\pi} \langle \Proj_{\pi} \varPhi_1, \Proj_{\pi} \varPhi_{2}\rangle d\mu_{P}(\pi) + 
\langle \varPhi_1, \mathscr{E}_2 \rangle_{reg} + \langle \mathscr{E}_1, \varPhi_2 \rangle_{reg}
 \end{equation}
 where $\Proj_{\pi}$ is the {\em regularized} projection onto the space of  automorphic $\pi$, 
and $$\mathscr{E}_i = \Eis( (\varPhi_i)^*_N)$$
\end{Proposition*}
We will call the additional contribution, $ 
\langle \varPhi_1, \mathscr{E}_2 \rangle_{reg} + \langle \mathscr{E}_1, \varPhi_2 \rangle_{reg}$
, in the regularized Plancherel formula  the {\em degenerate contribution}.
\begin{Rem} If $\varPhi_{1}$ is of rapid decay, the formula continues to hold with $\langle \mathscr{E}_1, \varPhi_2 \rangle_{reg}=0$
\end{Rem}
\proof Firstly our assumptions insure that all the terms of \refs{RIP} are well defined.
Moreover since the exponents of $\varPhi_{1}$ and $\varPhi_{2}$ are $>1/2$ and not of the form $\chi.|.|$ with $\chi$
quadratic, the representations underlying $\mathscr{E}_i$ have no subquotient isomorphic
to a standard automorphic representation, and thus
$\Proj_{\pi} \mathscr{E}_i = 0 $ for $i=1,2$ and all $\pi$. Similarly $\peter{\mathscr{E}_1,\mathscr{E}_2}_{reg}=0$ by our assumption that $S_{1},\ S_{2}$ are disjoint and have real parts
larger than $1/2$.  So it is enough to check,
with $\bar{\varPhi}_i = \varPhi_i - \mathscr{E}_i$, that
$$\langle \bar{\varPhi}_1, \bar{\varPhi}_2 \rangle_{reg} = \int_{\pi} \langle \Proj_{\pi}\bar{ \varPhi}_1, \Proj_{\pi} \bar{\varPhi}_2 \rangle d\mu_{P}(\pi)$$
but then $\bar{\varPhi}_i$ belongs to $L^2(\bfX_{\PGL_{2}})$.  \qed

\subsection{(Regularized) triple products.} \label{triple}

In this section we establish the following:  Suppose that $\pi_i$ are generic standard automorphic
representations, at least one of which is Eisenstein. 
Then, for each factorizable vector $\varphi_i = \otimes_{v} \varphi_{i,v}$, we have:
\be\label{ichino-symmetric} \frac{ \left| \int_{\bfX} \varphi_1 \varphi_2 \varphi_3(g)dg  \right|^2}{\|\varphi_{1}\|^2_{can}
\|\varphi_{2}\|^2_{can} \|\varphi_3\|_{can}^2}  = \frac{1}{8 (\disc F)} \prod_{v}^*\frac{ \zeta_v(1)}{\zeta_v(2)^3}  \frac{ |L_{W}(\varphi_{1,v},
\varphi_{2,v}, \varphi_{3,v})|^2}{\prod_{i=1}^3 \langle \varphi_{i,v}, \varphi_{i,v} \rangle } 
\ee
where the local factors are equal at almost all places, 
to $$L_v := \frac{\zeta_v(1)}{\zeta_v(2)} \frac{L(\frac{1}{2}, \pi_1 \otimes \pi_2 \otimes \pi_3)}{\prod_{i=1}^{3} L(1, \Ad, \pi_i)}. $$
If all $\pi_i$ are Eisenstein, the integration on the left-hand side
is to be understood by regularization.

\begin{Rem}
Let us compare this result with that in \cite{Ich}. 
Recall first the relation $dg =  (\disc F)^{-1} \xi_F(2) d^{\tau}g$ between
our measure and Tamagawa measure. It follows that,
if we replace $dg$ by $d^{\tau}g$, then \eqref{ichino-symmetric}
holds with a factor $\frac{1}{8}\prod_{v} \frac{\zeta_v(1)}{\zeta_v(2)^4}$. 
In the work of \cite{Ich}, the $\zeta_v(1)/\zeta_v(2)^4$ does not occur. 
\end{Rem}

From this, we deduce -- by summation through an orthogonal basis of $\pi_1$,
with respect to the canonical norm -- the following:
\begin{equation} \label{ichino-alternate} \frac{ \|\Pi_{\pi_{1}} (\ov{\varphi_{2}\varphi_{3}})\|^2}{ \|\varphi_2\|_{can}^2 \|\varphi_3\|_{can}^2} =  \frac{1}{8} (\disc F)^{-1} \prod_{v}^* \frac{ \zeta_v(1)}{\zeta_v(2)^3} \frac{ \| I_{\pi_{1,v}} (\varphi_{2,v}\otimes\varphi_{3,v} )\|^2}{ \langle \varphi_{2,v}, \varphi_{2,v} \rangle \langle \varphi_{3,v} , \varphi_{3,v} \rangle},\end{equation}
with a.e. local factor equal to $L_{v}$ above.

Before we embark on the proof, we note that it is sufficient -- by a continuity argument --
to treat the case where no $\pi_i$ is a singular Eisenstein series (i.e. of the form $\chi \boxplus \chi$.) 
For instance, let us consider the case when $\pi_1, \pi_2$ are cuspidal, $\pi_{3}=1\boxplus 1$
and let us take a family $\varphi_3(t) \in \pi_3(t)$, where $\pi_3(t) = |\cdot|^t \boxplus |\cdot|^{-t}$
deforming $\varphi_3$ (i.e., $\varphi_3(t) \rightarrow \varphi_3(0)$ pointwise as $t \rightarrow 0$). 
 Then the left-hand side and right-hand side
of \eqref{ichino-symmetric}, denoted $L(t)$ and $R(t)$ respectively, do not necessarily 
depend continuously on $t$ when $t=0$.  
However, both $L^*(1, \Ad, \pi_3(t)) L(t)$  and $L^*(1, \Ad, \pi_3(t)) R(t) $
extend to continuous functions around $t=0$, and are equal for $t \neq 0$. If we denote by $LR$ their common limit, then $L(0) $ and $R(0)$ are both given by $L^*(1, \Ad, \pi_3) LR$,
and are therefore equal. 

\subsubsection{Upper bounds} \label{ss:ub}
Unfortunately, we shall use the beautiful formula \refs{ichino-symmetric} only for upper bounds.
We now explicate the bounds that are derived from it.

The local factors of \eqref{ichino-alternate} are equal, 
to $L_v := \frac{\zeta_v(1)}{\zeta_v(2)^3} \frac{L(\frac{1}{2}, \pi_1 \otimes \pi_2 \otimes \pi_3)}{\prod_{i=1}^{3} L(1, \Ad, \pi_i)} $
, at almost all places. 
Observe that -- taking into account bounds towards the Ramanujan conjecture -- 
$L_v$ is absolutely bounded above and below at nonarchimedean $v$. 
Let $S$ be the set of places where the local factor is not equal to $L_v$, together with all archimedean places, 
and suppose we are supplied with the estimate $$\frac{\Sob_{-d}^{\pi_{3,v}}(I_{\pi_{3,v}}(x_1 \otimes x_2))^2}{\langle x_1, x_1 \rangle \langle x_2, x_2 \rangle}
\leq B_v$$ for $v \in S$. 
We conclude --
inserting $\Delta_{\Aa}^{-d}$ and applying the uniqueness of trilinear functionals -- 
\begin{equation} \label{ichino-upper-bound} \frac{ \|\Delta_{\Aa}^{-d}\Pi_{\pi_{1}} (\ov{\varphi_{2}\varphi_{3}})\|^2}{ \|\varphi_2\|_{can}^2 \|\varphi_3\|_{can}^2}  \ll  A^{-|S|}  \frac{L(\pi_1 \otimes \pi_2 \otimes \pi_3,\frac{1}{2})}{\prod_{i=1}^{3} {L}^{*}(\Ad, \pi_i,1)}  
\prod_{v \in S}  B_v,\end{equation}
where $A$ is an absolute constant.

Let us note that:   given $d$, there exists $d'$ so that:
\be\label{projectioncoarse}
\int_{\pi \, \mathrm{generic}} \Sob_d(\Pi_{\pi}(\varphi_2 \varphi_3)) d\mu_P(\pi)  \ll \Sob_{d'}(\varphi_2) \Sob_{d'}(\varphi_3)
\ee  
(the Sobolev norms are relative to the canonical inner product.)
Note that this is easy if $\pi_2,\pi_{3}$ are cuspidal;
in that case it can be deduced directly from (S3b). 
In the remaining cases, by virtue of the bounds of \S \ref{S3cproof},
it suffices to check that for for $\pi$ generic and any $d\geq 0$ 
\begin{equation} \label{frida} \Sob_d(\Pi_{\pi}(\varphi_2 \varphi_3))  \ll \Sob_{d'}(\varphi_2) \Sob_{d'}(\varphi_3)\end{equation} 
for some $d'$ depending on $d$ only. 

For this
we appeal to the prior formula. The bound on $B_v$ is supplied
by \eqref{ipibound}; we use also
 the fact that $L(\pi\otimes\pi_{2}\otimes\pi_{3},1/2)$ is bounded polynomially in the $C(\pi_{i})$ and that $L(\pi_{i},\Ad,1)=C(\pi_{i})^{o(1)},\ i=1,2$. This yields a bound as in \eqref{frida}, 
 but only for factorizable $\varphi_2 = \otimes_{v} \varphi_{2,v} , \varphi_{3} = \otimes_v \varphi_{3,v}$,
 and where the bound on the right-hand side is instead
 $$C(\pi_1)^A C(\pi_2)^A \prod_{v} \Sob_{d}(\varphi_{2,v}) \Sob_{d}(\varphi_{3,v}).$$
 We now apply (S1d) (see also Remark \ref{s1dremark}), together with \eqref{gan}
 to absorb $C(\pi_1)^A, C(\pi_2)^A$ into the Sobolev norms.  Since we did not prove \eqref{gan},
 we draw attention to the fact that we do not need this last step in any application;
 it would be fine to retain the factor $C(\pi_1)^A C(\pi_2)^A$. 

\subsubsection{The Rankin-Selberg integral} \label{Subsec:rs}
In this section, we shall prove the main results under the assumption
that either $\pi_1$ or $\pi_2$ are cuspidal. In this case,
we may proceed by the usual Rankin-Selberg method. 

We recall here the definition and basic properties of the Rankin/Selberg integral: let $\pi_i,\ i=1,2,3$ be generic automorphic representations of $\GL_2(\Aa)$  such that the product of their central characters 
$\chi_{1}\chi_{2}\chi_{3}$ is trivial.
We assume that $\pi_3$ is Eisenstein, (say $\pi_{3}=\chi^+_{3}\boxplus\chi^-_{3}$  for a pair of characters satisfying
$\chi^+_{3}.\chi^-_{3}=\chi_{3}$) and $\pi_1$ is cuspidal. 

The (absolutely convergent) integral
$$L(\varphi_1, \varphi_2, \varphi_3):=\int_{\bfX} \varphi_1 \varphi_2 \varphi_3(g)dg,\ \varphi_{i}\in V_{\pi_{i}},\ i=1,\dots,3,\ $$
defines a linear functional on the space of the representation $\pi_{1}\otimes\pi_{2}\otimes\pi_{3}$. If the $\varphi_{i}$
 are factorisable vectors, so that $W_{\varphi_{i}}=\prod_{v}W_{i,v}$ and $\varphi_{3}=\Eis(f_{3})$, with $f_{3}=\otimes_{v}f_{3,v}$
 the Rankin-Selberg unfolding method yield the following factorization (if $\chi_{3}^+\not=\chi_{3}^-$)
 \be\label{RSwhit}\int_{\bfX} \varphi_1 \varphi_2 \varphi_3(g)dg=\prod_{v}^*\frac{ L_{RS,v}(W_{1,v}, W_{2,v}, f_{3,v}) }{\zeta_v(1)^{1/2}}
 \ee
 where (see \eqref{RSlocal}) $$\frac{L_{RS,v}(W_{1,v}, W_{2,v}, f_{3,v})}{ \zeta_{v}(1)^{1/2}}= \int_{N(F_{v})\bash \PGL_{2}(F_{v})}W_{1,v}W_{2,v}f_{3,v}(g)dg;$$ and, for almost every $v$,  
 \begin{equation} \label{rs-ae} \frac{L_{RS,v}(W_{1,v}, W_{2,v}, f_{3,v})}{W_{1,v}(1) W_{2,v}(1) f_{3,v}(1)} =\zeta_{v}(1)^{1/2} 
\frac{L_{v}(\pi_1 \otimes \pi_2\otimes\chi_3^{+},{1/2})}{L_{v}(\chi_3^{+}/\chi_3^{-},1)}.\end{equation}

 Taking residue at the pole point $\chi_{3}^+=|.|^{1/2}$, $\varphi_{2}=\ov\varphi_{1}$, and using \eqref{rescomp}, yields: 
 $$  \peter{\varphi_{1},\varphi_{1}}_{\bfX} = 2 (\disc F)^{1/2}\frac{\xi_{F}(2)}{\xi_F^*(1)}
 \prod_{v}^*  \langle W_{1,v} , {W_{1,v}}\rangle
 $$ with $ \langle W_{1,v} , {W_{1,v}}\rangle=\zeta_{k}(1)L(\pi,\Ad,1)/\zeta_{k}(2)$ for a.e. $v$. This proves \refs{canreg} in the cuspidal case.
 
 From these remarks,  \refs{RSwhit}, together with the definition
 of the canonical norm in \S \ref{canGL2}, we deduce that:
\be\label{ichino} \frac{ \left| \int_{\bfX} \varphi_1 \varphi_2 \varphi_3(g)dg  \right|^2}{\|\varphi_{1}\|^2_{can}
\|\varphi_{2}\|^2_{can}}  = \frac{1}{4} (\disc F)^{-1} \prod_{v}^*\frac{ \zeta_v(1)}{\zeta_v(2)^2}  \frac{ |L_{RS,v}(W_{1,v}, W_{2,v}, f_{3,v})|^2}{\prod_{i=1}^2 \langle W_{i,v} , W_{i,v} \rangle } 
\ee
with almost all factors equal to $\frac{L(\pi_{1,v} \otimes \pi_{2,v} \otimes \pi_{3,v},\frac{1}{2})}{\zeta^{-1}_v(1) \prod_{i=1}^3 L(\pi_i, \Ad, 1)}$. Note that, since $\pi_3$ is Eisenstein, the product
would not converge without the inclusion of the $\zeta_v(1)$ factor. 
Now \eqref{ichino-symmetric} -- in the case where one of $\pi_1, \pi_2$ is cuspidal, 
and $\pi_3$ is not isomorphic to $1 \boxplus 1$ --
follows, taking into account the equality between $L_{RS}$ and $L_W$ that was already established 
in \S \ref{subsec:wii}, together with the relation \eqref{eis-norm-eqn}
between $\prod_{v} \langle f_v, f_v \rangle$ and the canonical norm on $\pi_3$.

\subsubsection{A regularized triple product}\label{regulartriple}
We discuss now the situation when {\em all} $\pi_i$ are Eisenstein, i.e. $$\pi_{i}=\Eis(\mcI(\chi_{i}^+,\chi_{i}^-)),\ \varphi_{i}=\Eis(f_{i}),$$
we the $\chi_{i}^\pm$ are unitary.
There are two equivalent definitions of the regularized integral $\pi_1 \otimes \pi_2 \otimes \pi_3 \rightarrow \mathbb{C}$; we define them
and prove their equivalence: 

Set
\begin{equation} \label{defn-1} L_a: \varphi_1 \otimes \varphi_2 \otimes \varphi_3 \mapsto  \int^{reg}_{\bfX_{\PGL_{2}}} \varphi_1 \varphi_2 \varphi_3,\end{equation} 
and set $L_c : \pi_1 \times \pi_2 \times \pi_3 \rightarrow \mathbb{C}$ to be the value at $s=0$ of the meromorphic continuation (from $\Re(s) \gg 1$)
of the following expression:
\begin{equation} \label{alternate} L_s:  \varphi_1 \otimes \varphi_2 \otimes \mathrm{Eis}(f_3) \mapsto \int_{\mathbf{N}(\Aa)\rmA(F) \backslash \PGL_2(\adele_F)} ((\varphi_1 \varphi_2)_{N} - (\varphi_1)_{N} (\varphi_2)_{N}) {f_{3}(s)}dg. \end{equation}
Note that the later expression is convergent for $\Re s\gg 1$ due to the rapid decay of 
$(\varphi_1 \varphi_2)_{N} - (\varphi_1)_{N} (\varphi_2)_{N}$ and unfolds to
\be\label{alternateunfold}s\mapsto \int_{\mathbf{N}(\adele) \rmA(\Aa) \backslash \PGL_2(\adele)} W_{\varphi_1} W_{\varphi_2}
f_{3}(s)dg\ee
which extends to an holomorphic function  in a neighborhood of $s=0$.

\begin{Lemma*}
 $L_a = L_c$. Moreover, if we write $L$ for the common value of these expressions,
\begin{equation}
\frac{\left| L(\varphi_1 ,\varphi_2 ,\varphi_3) \right|^2}{\|\varphi_{1}\|^2_{can}\|\varphi_{2}\|^2_{can} \|\varphi_3\|^2_{can}}  = 
 \frac{1}{8} (\disc F)^{-1}  \prod_{v}^*\frac{ \zeta_v(1)}{\zeta_v(2)^3}  \frac{ |L_v(\varphi_{1,v}, \varphi_{2,v}, \varphi_{3,v})|^2}{\prod_{i=1}^3 \langle \varphi_{i,v} , \varphi_{i,v} \rangle } 
\end{equation}
\end{Lemma*}
In this way, we have established \eqref{ichino-symmetric} in the remaining case also. 

\proof Considering the central characters $\chi_{i},\ i=1,2,3$ fixed, the pairs of characters $(\chi_{i}^+,\chi_{i}^-)$ $i=1,2,3$ such that $\chi_{i}^+\chi_{i}^-=\chi_{i}$ will be referred as the {\em parameters}. The set of parameters has the structure of a $3$-dimensional complex manifold
with infinitely many connected components.

The reasoning by which we derived \eqref{ichino-symmetric} in the ``at least one $\pi_i$ cuspidal'' case
may be applied to $L_c$, at least on an open set of parameters which intersects every connected component.  By analytic continuation this shows that $|L_c|^2$ is indeed given by \eqref{ichino-symmetric} everywhere. 

To establish $L_a =L_c$,  we consider the parameters in a given connected component: in other words, assuming the characters $\chi_{i}^\pm,\ i=1,2,3$ unitary, we consider any analytic deformation $$\varphi_{1}(s_1)\otimes\varphi_{2}(s_2)\otimes\varphi_{3}(s_{3})\\ \in\pi_{1}(s_1)\otimes\pi_{2}(s_2)\otimes\pi_{3}(s_3)$$indexed by the complex parameters
$(s_{1},s_{2},s_{3})\in\Cc^3$. 
Clearly, the integral \eqref{alternate} extends (via \eqref{alternateunfold}) to an holomorphic function on an open subset of $\Cc^3$ (containing $(0,0,0)$)
\be\label{alternate3}
(s_{1},s_{2},s_{3})\mapsto\int_{\mathbf{N}(\Aa)\rmA(F) \backslash \PGL_2(\adele_F)} ((\varphi_{1}(s_1) \varphi_{2}(s_2))_{N} - (\varphi_{1}(s_1))_{N} (\varphi_{2}(s_2))_{N}) {f_{3}(s_{3})}dg.\ee
 Let $\bfE({s_{1},s_{2}})$ be the Eisenstein series formed out of the exponents of the product $\Phi := \varphi_{1}(s_1).\varphi_{2}(s_2)$ which are of real part $>1/2$. Explicitly: let $S_{\Phi}$ be
 the set of exponents of $\Phi$; let $\Phi_{N, > 1/2} ^*$ be the part of $\Phi_N^*$
 that corresponds to exponents $\chi \in S_{\Phi}$ with $\Re(\chi) > 1/2$,
 and let $\bfE({s_{1},s_{2}}) := \mathrm{Eis}(\Phi_{N, >1/2}^*)$, interpreted by analytic continuation
 if there exists $\chi \in S_{\Phi}$ with $\Re(\chi) \in (1/2, 1)$. 
    
The map $(s_1, s_2) \mapsto \bfE(s_1, s_2)$ defines a meromorphic function on an open subset of $\Cc^2$; in an open subset of $\Cc^3$, one has
$$\int^{reg}_{\bfX_{\PGL_2}}
\varphi_{1}(s_1) \varphi_{2}(s_2) \varphi_{3}(s_{3}) = \int^{reg} (\varphi_{1}(s_1)\varphi_{2}(s_2) - \mathbf{E}(s_{1},s_{2})) \varphi_{3}(s_{3}). $$
Considering Fourier expansions, one see that, given any $N>1$, one has
\begin{equation} \label{fitano} \varphi_{1}(s_1)\varphi_{2}(s_2)(x) - \mathbf{E}({{s_{1},s_{2}}})(x)\ll_{N} \Ht(x)^{-N}, x \in \mathbf{X}_{\PGL_2},\end{equation} 
as long as $\Re s_{1} \gg_N \Re s_{2}\gg_{N} 1$.

 Therefore, there is an open subset of $\Cc^3$ (in which $s_{3}>1/2$) so that
the previous integral is absolutely convergent. It unfolds to
\begin{equation} \label{alternate2} \int_{\mathbf{N}(\adele)\rmA(F) \backslash \PGL_2(\adele)} ((\varphi_{1}(s_1) \varphi_{2}(s_2))_N - \mathbf{E}({s_{1},s_{2}})_{N}) f_{3}(s_{3}) dg. \end{equation}
Indeed, \eqref{fitano} implies that the constant term $((\varphi_{1}(s_1) \varphi_{2}(s_2))_N - \mathbf{E}({s_{1},s_{2}})_{N}$ is {\em bounded}, and from this one justifies the unfolding process.

The definitions \eqref{alternate3} and \eqref{alternate2} are {\em a priori} convergent
in different regions and they cannot be compared directly. Nonetheless, they coincide on the intersection of their domains of holomorphic continuation. Indeed,  there exists a nonempty open set of parameters, intersecting
every connected component of either domain of holomorphic continuation, 
so that both \eqref{alternate3} and \eqref{alternate2} can 
be defined by regularizing the integral over $\mathbf{N}(\adele)\rmA(F) \backslash \PGL_2(\adele)$;
this being done, their difference vanishes by invariance of the regularized integral. 

To be more specific, we have, whenever $|s_1| +|s_2| \leq A$, the bound
$$(\varphi(s_1) \varphi(s_2))_N(a(y) k) \ll \max(|y|,|y|^{-1})^{1+A},$$
with bounds of similar nature for $\varphi(s_1)_N \varphi(s_2)_N$ and also 
$\mathbf{E}({s_{1},s_{2}})_{N})$. Therefore, if $\Re(s_3) \gg_{A} 1$,
and we write out
the integrals \eqref{alternate2} and \eqref{alternate3} in the Iwasawa decomposition, 
we obtain functions of $a(y) k$ which are integrable in the region $|y| \leq 1$.
By contrast, in the region $|y| \geq 1$ they are asymptotic to sums of 
finite functions of $y$, i.e., functions whose translates span a finite-dimensional vector space;
these may be regularized 
as in \S \ref{Rreg},   and our conclusion follows. 
 \footnote{The reader
may wish to consider the following simpler example of this reasoning: the characteristic function of $[0,1]$ and
the characteristic function of $[1, \infty]$ have Mellin transforms, respectively, $\frac{1}{s} \ (s > 0)$
and $-\frac{1}{s} \ (s < 0)$. However, the fact that the meromorphic extensions are negative to each other can be deduced without computation: the sum of $f_1 + f_2$ is the constant function,
and its Mellin transform in any regularized sense must be zero wherever defined, by invariance. }

\qed

\subsection{An example: $F = \mathbb{Q}$ at full level.} \label{fulllevel}
For $F=\Q$, we have in particular $\xi_F^*(1) = \disc F = 1$. Let us describe
some explicit consequences of the foregoing remarks. In particular, we elaborate
on the remark, contained in the introduction, that the ``evident'' identity \eqref{commute} 
gives rise to an identity between families of $L$-functions.

\subsubsectionind{}

 Let $\mathcal{M}$ be the set of even Maass cusp forms on $\SL_2(\Z) \backslash \H$.  Let 
 $\xi(s) =\pi^{-s/2} \Gamma(s/2) \zeta(s)$. We shall use $\Lambda$ to denote the completed $L$-function.  For $\varphi \in \mathcal{M}$ define
  $$\eta_{\varphi}(x,y) = \Lambda( \varphi,1/2+ x+y) 
\Lambda( \varphi,1/2 + x - y),$$
the corresponding definition for an Eisenstein series $E_{s}\in\Eis(|.|^s,|.|^{-s})$ is: 
$$\eta(s,x,y) = \xi(1/2 + s +x+ y ) \xi(1/2 - s +x+ y) \xi(1/2 +s+x-y) \xi(1/2-s+x-y);$$

The quotient of $\mathbf{X}$ by the maximal compact $K$ may be identified
with the quotient of $\SL_2(\Z) \backslash \mathbb{H}$,
by $z \mapsto -\overline{z}$.  The measure of the quotient is $\frac{\pi}{3}$, i.e.
it arises from $2 \frac{dx \ dy}{y^2}$ (but if we work on $\SL_2(\Z) \backslash \mathbb{H}$,
it is simply $\frac{dx \ dy}{y^2}$). 

For $s \in \C^\times$ let $E_s$ resp. $E^*_s$ be the result of applying the Eisenstein intertwiner
to the vector $f \in |\cdot|^{s} \boxplus |\cdot|^{-s}$ whose restriction to $K$ is $1$
resp. $\xi(1+2s)$. Then $E_s$ and $E_s^*$ descend to functions on $\SL_2(\Z) \backslash \mathbb{H}$;
 their constant terms are $$y^{1/2+s} + \frac{ \xi(2s)}{\xi(1+2s)} y^{1/2-s}\ \hbox{
resp. }\xi(1+2s) y^{1+2s} + \xi(1-2s) y^{1-2s}.$$  For $s\in i\Rr^\times$, we have by \eqref{eis-norm-eqn}, 
$\langle E_s, E_s \rangle_{can}  = 2 \xi(2) =  \pi/3$ while $$\langle E_s^*, E_s^* \rangle_{can} =2  \xi(2)  \xi(1+2s) \xi(1-2s)\ra \infty,\ s\ra 0$$

By \eqref{ichino-symmetric},  if $\psi \in \mathcal{M}$ or is an Eisenstein series, 
$$\langle E_{t_1}^* E_{t_2}^* , \psi \rangle^2 = \frac{\|\psi\|_{can}^2}{2 \xi(2)}  \frac{\eta_{\psi}(t_1, t_2)^2}{\Lambda^*(1, \Ad, \psi)}$$

For $t_1 \neq \pm t_2 \in i\mathbb{R}$:
\be\label{structure-constants} E_{t_1}^* E_{t_2}^* = \sum_{t_1' = \pm t_1, t_2' = \pm t_2} \frac{ \xi(1+  2 t_1') \xi(1 + 2 t_2') }{\xi(2 + 2t_1' + 2t_2')} E_{1/2 +t_1' +t_2'}  \ee
$$+ \sum_{\varphi \in \mathcal{M}} \frac{ \eta_{\varphi}(t_1, t_2)}{\sqrt{2 \xi(2) L( \Ad, \psi,1)}\|\varphi\|_{can}} \varphi  + \frac{1}{2\pi}  \int_{s=0}^{i \infty} \frac{  \eta(s, t_1, t_2) }{\xi(1+2s) \xi(1-2s) } E_s^* $$

The first and third term on the right-hand side are, of course, intimately related,
through a process of contour-shifting; 
indeed, if we compute the constant term of the right-hand side,
one finds that (after shifting contours) the residues of the third term
match exactly, and indeed cancel, part of the constant term of the first term. 

It is interesting to substitute the point $i$ (``period over a non-split torus with class number one'').  Noting that,
$$E_s^*(i) = 2^{s+3/2}\xi_{\Q(i)} (1/2+s),$$  the resulting formula
relates $\xi_{\Q(i)} (1/2+t_1) \xi_{\Q(i)} (1/2+t_2)$ with
$\xi_{\Q(i)}(1 \pm t_1 \pm t_2)$ and integrals over the critical line.  The contribution of cusp
forms may be expressed in terms of Fourier coefficients of half-integral weight. 
In this case, {\tt PARI/GP} evaluates (with
$t_1 = 0.9 i, t_2 = -1.3 i$) the noncuspidal part of the right hand side to
$0.10553970$, whereas the left-hand side comes to $0.10554092$.

\subsubsection{The associativity formula}
\eqref{structure-constants} expresses a ``multiplication table'' for forms; 
this is of course constrained by the associativity of multiplication.
These constraints lead (among other things) to identities generalizing that of Kuznetsov:

Set  $$f_1(t_1, t_2, t_3, t_4) = \frac{ \xi(1 + 2t_3 ) \xi(1+2t_4) }{ \xi(2 + 2 t_3 + 2 t_4)} \eta(t_1, t_2, 1/2  + t_3 + t_4)$$
and take $\Delta (t_1, t_2, t_3, t_4) = \sum_{\pm, \pm} ( f(t_1, t_2, \pm t_3, \pm t_4) + \sum_{\pm, \pm} 
f(t_3, t_4, \pm t_1, \pm t_2)$.  
Then the function $\Sigma$ defined by:
$$\Sigma = \int_{s=0}^{i \infty} \frac{ds}{2 \pi} \frac{\eta(s, t_1, t_2) \eta(s, t_3, t_4)}{\xi(1+2s) \xi(1-2s)} 
+ \frac{3}{\pi} \sum_{\varphi \in \mathcal{M}} \frac{\eta_{\varphi}(t_1, t_2) \eta_{\varphi}(t_3, t_4)}{\Lambda(\varphi, \Ad,1)} 
+\Delta(t_1, t_2, t_3, t_4)$$
is {\em invariant under all permutations of coordinates.} This is a spectral identity
between families of $L$-functions; a version was first discovered by Kuznetsov. 

We tried to test this numerically.  For $(t_1, t_2, t_3, t_4) = (1.2 i, 1.5i, 3i, 4i)$,  the difference
between the {\em noncuspidal} parts of $\Sigma(t_1, t_2, t_3, t_4) - \Sigma(t_1, t_3, t_2, t_4)$. 
This is estimated by  {\tt PARI/GP} to be  $8.29 \times 10^{-11}$.   (To get a sense 
of size, each of the degenerate terms have size $\sim -0.005$, and the difference
between the two degenerate terms has size $\sim 2.6 \times 10^{-9}$.) Nonetheless, this quantity  -- while very small -- is in fact substantially larger
than the contribution of the first Maass form. 
It is quite possible that error in the numerical integration is responsible for most of the difference.  
 
\subsubsection{The Motohashi formula}\label{motid}
One may also take the period of \eqref{structure-constants} over a split torus, i.e.,
integrate $\int_{y \in \mathbb{R}_+} E_{t_1}(a(y)) E_{t_2}(a(y)) d^{\times}y$;
the integral does not converge, but can be regularized as in \S \ref{Rreg}. 

The resulting formula relates, on the left-hand side, $$\int_{\nu \in i \mathbb{R}} \xi(1/2 + t_1 + \nu) \xi(1/2 - t_1 +  \nu) \xi(1/2 + t_2 + \lambda - \nu) \xi(1/2 - t_2 + \lambda - \nu)$$
to, on the right hand side, 
\begin{gather*}\sum_{\mathcal{M}} \Lambda(1/2+ t_1 +t_2 , \varphi) 
\Lambda(1/2 +t_1 - t_2, \varphi )  \Lambda(1/2+\lambda, \varphi)\\\ \  + \hbox{\ Eisenstein and degenerate terms}
\end{gather*}  This is an example of a formula of Motohashi \cite{Mot};
it is perhaps most interesting to let $t_1, t_2, \lambda =0$. 

The possibility of thus deriving the Motohashi formula was remarked in \cite{MV-ICM}
in \S 4.3.3 together with the remark that the divergence of the resulting integrals would
``cause considerable technical difficulty.'' The method of using regularized integrals
thus settles that issue in a satisfactory way.  We would like to remark
that it is necessary to be wary of the following point:
the regularized integral
{\em cannot} be interchanged naively
with the continuous integral over the Eisenstein series that occurs on the right
of \eqref{structure-constants}. In ``commuting'' the two,
we obtain an extra factor of $\frac{ \eta(-1/2-\lambda, t_1, t_2)}{\xi(2+2\lambda)}+
\frac{ \eta(-1/2+\lambda, t_1, t_2)}{\xi(2-2\lambda)}$.

\subsubsection{Commentary}
For the purpose of analytic number theory the above formulas, although beautiful, are insufficient;
one needs a family with more flexibility of ``test functions.'' This extra flexibility is provided by the work of
Motohashi and Kuznetsov. From the point of view of the method above, this can 
be obtained by varying the choice of test vector in the representations underlying $E_{t_1}^*$ and $E_{t_2}^*$. 

The associativity formula, and that of Motohashi, have played an important -- though not always explicit -- role
in analytic number theory on $\GL_2$. We invite the reader to see ``shadows'' 
of these formulas hiding in various other papers on the subject. 

The primary advantage of the formalism above seems to be that it generalizes
immediately to ramified settings. One may, for instance, replace $E_{t_1}^*$ by an Eisenstein series
in the representation $\chi \boxplus 1$ to get a new formula; or one may replace $\Q$ by a number field. 
It is a very interesting question to investigate more deeply the ``test function'' version
of such formulas in a general setting, and to study the integral transforms -- both archimedean and $p$-adic -- that intervene.

Such general formulas would probably lead
to, among other things, a good exponent in the subconvexity theorem. 
In the present paper we have taken a ``short cut'' to subconvexity.

\section{Proof of the theorems.}\label{part:thmproof}

\subsection{Subconvexity of character twists.}\label{proof}

In the present section we shall prove the following theorem which is a special case
of our main result, Theorem \ref{GL12}, but which is also necessary in its proof. 

  We continue with notation as in \S \ref{notation-global}.

\begin{Theorem}\label{thmcharacter} There are absolute constants $\delta, A\geq 0$ such that for $F$ a number field, $\pi$ an automorphic representation of $\GL_2(\Aa_{F})$, and $\chi$ a unitary (Hecke) character of $F^\times\bash\Aa^\times_F$, 
$$L(\pi \times\chi,\frac{1}{2} )\ll_{F, \pi} C(\chi)^{1/2-\delta}, $$
where $C(\chi)$ denotes the analytic conductor, and our conventions about $\ll$ are as in
\S \ref{ipnot}. 
In particular,
$L(\chi,\frac{1}{2})\ll_F  C(\chi)^{1/4-\delta}.$
\end{Theorem}

The proof follows the description of \S \ref{classical}; we try below to give a ``translation''
of the adelic steps to the steps in \S \ref{classical}. 

 Let us recall that over $\Qq$, the first result in the direction of Theorem \ref{thmcharacter} are due to Good \cite{Good} (in the $t$-aspect) and to Duke-Friedlander-Iwaniec in the conductor aspect \cite{DFI1}. Over general number fields
 the conductor aspect was obtained by the second author \cite{Ve} and Diaconu-Garrett for the $t$-aspect \cite{DG}. The various methods used in these papers, although superficially rather different
 and having different strengths and weaknesses, nonetheless are closely linked;
 they are all, in various ways, related to versions of the identity described in \S \ref{motid}. 
In particular, it is possible to redo the proof below
in a way that is much closer to the proof of Theorem \ref{GL12}, i.e. removing the ergodic theory
and substituting explicit spectral expansions. 

\subsubsection{Notation} 
Let $H^{(1)} = \{a(y): |y| = 1\}$. 
Given $\chi$ a character as in the Theorem,  we define the following signed measure $\mu_{\chi}$ supported on a closed orbit of $H^{(1)}$ on $\bfX$:
\begin{equation} \label{muchidef} \mu_{\chi}(\varphi) = \int_{F^\times\bash \Aa^{(1)}} \varphi(a(y))
\chi(y) d^\times y,  \ \ \ \mu = \mu_{1}.\end{equation}

Let us note that \eqref{muchidef} makes sense for any function on 
 $\rmB(F) \backslash \PGL_2(\Aa)$, i.e. $\mu_{\chi}$ has a canonical lifting $\tilde{\mu}_{\chi}$ to that space.
Given  $g \in \GL_2(\Aa)$, we denote by $\mu_{\chi}^g$ the translate of $\mu_{\chi}$ by $g$
 (i.e. $\mu_{\chi}^g(\varphi)=\mu_{\chi}(g.\varphi)$), and similarly define $\tilde{\mu}_{\chi}^g$. 
For $t\in\Rr_{>0}$, we pick $y_{t}\in\Aa^\times$ such that $|y_{t}|_{\Aa}=t$. 

We shall be primarily interested in the translates $\mu_{\chi}^g$
when $g$ is of the form $a(y) n(T)$, for $T \in \adele$. 

In trying to get some geometric intuition for these measures, we suggest
that the reader bear in mind the following simple example: 
 $F = \mathbb{Q}$, $T  = ``p \in \Q_p \hookrightarrow \adele$'', the projection of $\supp(\mu_{1}^{a(q^{-1}) n(T)})$ to the space of lattices is the set $\{\Lambda_x: x \in (\Z/q\Z)^{\times}\}$.
Here $\Lambda_x$ are as in \S \ref{classical}; recall in particular that $\Lambda_x = n(x) \Lambda_1$. Therefore, in classical terms,
the measures $\mu_{\chi}^{a(y) n(T)}$ will correspond to  {\em certain orbits
of the discrete horocycle flow}; however, the {\em signs} of the measure $\mu_{\chi}$
encode arithmetic information (e.g., whether or not $x$ is a quadratic residue).  

\subsubsection{A sketch of the proof in the simplest case: $F \neq \Q$ and $\pi$ cuspidal} \label{outline}
Unfortunately, there are two cases which introduce extra notation and small complexities: $F = \Q$, 
and $\pi$ Eisenstein.  Thus we advise the reader to read
carefully the present subsection, where we sketch the proof in the case when
neither of these complications exist.

The reasoning has the following three stages. Fix $\varepsilon > 0$ and $\kappa \in (0,1)$
and write, for typographical simplicity, $Q = C(\chi)$. 
\begin{enumerate}
\item[Step 1.]  Theorem \ref{thmcharacter} $\impliedby$  equidistribution result for $\mu^{g}_{\chi}$, for suitable $g$:

\underline{Claim 1:} There's $\varphi \in \pi$, $T=(T_{v})_{v} \in \adele$
with  $|T|  \gg Q^{1-\eps}$ and $t \in [Q^{-1-\kappa}, Q^{-1+\kappa}]$: 
$$ Q^{-1/2-\eps} \left| L(\pi \otimes \chi,1/2) \right|  \ll_{\eps,\kappa,\pi} 
  | \mu_{\chi}^{a(y_t) n(T)}(\varphi)|
 + Q^{-\kappa/2}.$$

\item[Step 2.]   \S \ref{udtt-cor} proves
 equidistribution of translates of $\mu_{\chi}$ in the case $\chi=1$. 

\underline{Claim 2:}  For some absolute $\delta > 0$, and $t$ as above
 $$\mu^{a(y_t) n(T)} (f)  \ll \Sob^\bfX(f)(|t|_\Aa^{1/2}+|\shift|_{\Aa}^{-\delta}),
 \mbox{$f \in L^2_{0}(\bfX_{\PGL_2})$ and smooth} $$

\item[Step 3.]  An application of \S \ref{Prinzip} allows us to bound
$\mu_{\chi}^{a(y_t) n(T)}$ starting from \underline{Claim 2}. 

\underline{Claim 3:}  $|\mu_{\chi}^{a(y_t) n(T)}(\varphi)|   \ll_{\pi} Q^{-\delta'} \Sob^{\bfX}(\varphi)$. 

The utility of $F \neq \mathbb{Q}$ comes in here:
we take advantage of the fact that we can find many pairs of distinct prime
 ideals with the same norm. We use this to construct a suitable measure $\sigma$ with which to apply \S \ref{Prinzip-quant}, with $\sigma$ supported entirely on the group $H^{(1)}$
and commuting with $n(T)$\footnote{Let $K$ be a parameter to be chosen as a {\em fixed} positive power of $Q$. Consider all finite primes of $\Q$ that are contained in $[K,2K]$, are split in $F$ and  above which $\chi$ is unramified. 
  The number of such primes $p$ is then $\gg_{F} K/\log K$.  Above each such prime $p$, 
 let $v_1(p), v_2(p)$ be two {\em distinct places}.  Let $\sigma$ be the probability measure on $H^{(1)}$ which is the average of the Dirac measures at the $a(\varpi_{v_1(p)}^{-1} \varpi_{v_2(p)})$ for all such primes
 $p \in [K,2K]$.}. When $F = \mathbb{Q}$ one can only make such a measure supported ``near''
 $H^{(1)}$.  
\end{enumerate}

In terms of the discussion in \S \ref{classical}, and in particular \S \ref{classicaled}, 
Step 1 amounts to the remark
that (C) or (C2) implies subconvexity and Step 2 amounts to the implication $(B) \implies (C)$.

\subsubsection{The general proof}

For $T\in\Aa$ and $h$ a smooth compactly supported function on $\Rr_{>0}$ we define the measure
$$\mu^{n(T)}_{\chi,h}(\varphi):=\int_{\Rr_{>0}}h(t)\chi(y_{t})\mu_{\chi}^{a(y_{t})n(T)}(\varphi)\dti t.$$

\begin{Lemmat}  \label{lrep} (variant of \underline{\rm{Claim 1}}) Notation as in the theorem, set $Q:=C(\chi)$. For any $\eps>0, \kappa \in ]0,1[$, there exists a smooth vector $\varphi \in \pi$, $T=(T_{v})_{v} \in \adele$ and a smooth, non-negative, bounded by $1$, function on $\Rr_{>0}$,
$h=h_{Q,\kappa}$ say, supported in the interval $[Q^{-1-\kappa}, Q^{-1+\kappa}]$ such that
\label{last}  
\be\label{truncated}Q^{-1/2-\eps} \left| L(\pi \otimes \chi,1/2) \right|  \ll_{\eps,\kappa,\pi} 
  | \mu_{\chi,h}^{n(T)}(\varphi)|
 + Q^{-\kappa/2}\ee
 (if $\pi$ is not cuspidal,  replace $\mu_{\chi}(\varphi)$ by $\tilde{\mu}_{\chi}(\varphi - \varphi_N)$) 
 and moreover:
\begin{enumerate}
\item $T_{v}=0$ unless $v$ is archimedean or $\chi$ is ramified at $v$;
\item $\frac{\log |T|_{\Aa} }{\log Q} \in [1-\eps, 1+\eps]$ if $Q \gg_{\eps} 1$; 
\item  $\Sob^{\pi}_d(\varphi) \ll_{d,\pi} 1$ for any $d$. In particular, if $\pi$ is cuspidal, then,  for any $d$, $\Sob^\bfX_d(\varphi)\ll_{d,\pi} 1$. 
\end{enumerate} 
\end{Lemmat}
 \proof
For finite $v$, we take $W_{\varphi,v} \in \Whit_{\pi_v}$ to be the new vector and $T_v = \varpi_v^{-r}$,
where $r$ is the conductor of $\pi_v$, just as in Lemma \ref{lem:hjllocalnonarch}.  For archimedean $v$, we choose $W_{\varphi,v}$ and $T_v \in F_v$ according to Lemma \ref{lem:hjllocal}.
Put $T = (t_v)_{v} \in \adele$ and let $\varphi \in \pi$ be the preimage of $\otimes W_{\varphi,v}$
under the canonical intertwiner from $\pi$ to its Whittaker model.  

The third assertion 
$\Sob_d^{\pi}(\varphi) \ll_{d,\pi} 1,$  follows since $\Sob_d(W_{\varphi,v}) \ll \Cond(\pi_v)^{N(d)}$, for some $N(d)$ depending on $d$ -- noted in Lemma 12.3 in the archimedean case, and
immediate in the nonarchimedean case -- 
 and from  \S \ref{canGL2}. 

It follows then from Lemma \ref{lem:hjllocalnonarch},  Lemma \ref{lem:hjllocal}, and the results
of \S \ref{PGL2} that
\begin{equation} \label{startbound} L(\pi \otimes \chi,1/2) \ll_{\eps} (C(\pi)Q)^\eps Q^{1/2} \int_{F^\times\bash\Aa^\times} \varphi(a(y) n(T)) \chi(y) d^{\times}y\end{equation}
We need to pass from this statement to the desired property (3) by truncating the range of the $y$-integral. 
That is carried out using similar reasoning to \cite[Lemma 11.9]{Ve}
; it can be considered the geometric form of the approximate functional equation. 
Let $h$ be a fixed smooth function on $\Rr_{>0}$ with values in $[0,1]$, which is $1$ on $(0,1]$, $0$ on $[2,\infty)$. Take $A = Q^{-\kappa-1}$ and $h_{A}:t\ra h(t/A)$.

Write  as a shorthand,
$$f(t) = \chi(y_{t})\mu^{a(y_t) n(T)}_{\chi}(\varphi)$$
so that the integral on the right-hand side of \eqref{startbound} is given by $\int_{0}^{\infty} f(t) d^{\times}t$, and, more generally, the Mellin transform $F(s) = \int f(t) t^{s} d^{\times}t $
is given by $F(s)=\ell^{\chi |\cdot |^s}(n(T) \varphi) $.

Given $\eps>0$ small, we will need to bound $F(s)$ for $\Re s=-1/2-\eps$.  
We claim that
\begin{equation} \label{Fbound}
|F(s)| \ll_{\pi, \eps} (1+|s|)^2 Q^{1/2+3 \eps} \mbox{ when }\Re(s) = -1/2 - \eps.  \end{equation} Crudely, 
the $L$-function in front contributes $Q$, whereas the ramified factors
contribute $Q^{-1/2}$, and the rest is of size $1$. 

To be more precise, 
let $S$ denote the (finite) set of places where either $v$ is infinite
 or $\varphi_{v}$ is not spherical, or $T_{v}\not=0$, or $\psi$ is ramified. 
 Then $F(s)$ may be expressed as:
 $$ L^{(S)} (\pi \otimes \chi,1/2+s)  \prod_{v \in S} 
\ell^{\chi_v|.|^s_{v}}(n(T_{v}) W_{v,\varphi}) \prod_{v\notin S} \frac{\ell^{\chi_v|.|^s_{v}}(n(T_{v}) W_{v,\varphi})}{L(\pi_v \otimes \chi_v,1/2+s)}.$$
The product over $v \notin S$ is equal to $1$.   For $v \in S$, each factor is bounded by 
$\ll_{\pi, \eps} \Cond(\chi_v)^{-1/2+\eps}$, as follows from 
 Lemma \ref{lem:hjllocalnonarch} and Lemma \ref{lem:hjllocal}.
 Finally, by the convexity bound together with bounds towards Ramanujan\footnote{More precisely, 
 it is necessary to bound each factor $L(\pi_v \otimes \chi_v, 1/2+s)^{-1}$.  Suppose $v$ is finite
 and that $\psi$ is unramified at $v$.  If $\pi_v$
 is tempered, this is bounded by $(1+q_v^{\eps})^2$.  If $\pi_v$ fails to be tempered, 
 then $\pi_v$ is a twist of a spherical representation, as is $\pi_v \otimes \chi_v$. If the local $L$-factor
 is not identically $1$, then necessarily $\pi_v \otimes \chi_v$ is spherical. Since $v \in S$, 
 it must be that $\pi_v$ was ramified at $v$; the claimed bound follows.}  
  $$L^{(S)} (\pi \otimes \chi,s) \ll_{\eps} (1+|s|)^2C(\pi\otimes\chi)^{1/2+\eps} \ \ (\Re s = - \eps).$$
Noting that $|S|=o(\log(C(\pi\otimes\chi)+\disc(F)))$ as $C(\pi\otimes\chi)\ra\infty$,
and since $C(\pi\otimes\chi)\ll C(\pi)C(\chi)^2$, we get \eqref{Fbound}. 

By Mellin inversion ($H$ being the Mellin transform of $h$),  
$$  \int_{\Rr_{>0}}h_{A}(t) f(t) d^{\times}t  = \left| A^{1/2+\varepsilon}  \int_{\Re(s) = -1/2- \varepsilon} H(-s) 
F(s) \frac{ds}{2 \pi i}  \right| \ll_{\pi, \eps, h }
 Q^{\frac{-\kappa  }{2} + \eps} 
$$
The effective content of this statement is that the range $t \lesssim Q^{-\kappa-1}$
contributes very little to the integral $\int_0^{\infty} f(t) d^{\times} t$. 
A corresponding analysis yields the same statement for the range $t \geq Q^{\kappa-1}$, namely: 
$$\int_{\Rr_{>0}}(1-h_{Q^{\kappa-1}}({t}))\chi(y_{t}) \mu^{a(y_t) n(T)}_{\chi}(\varphi) d^{\times}t\ll_{\pi,\eps, h}Q^{-\kappa/2+\eps}.$$
We take $h = h_{Q^{\kappa-1}} - h_{Q^{-\kappa-1}}$ to conclude. 

The following will be useful later: If we replace $h$ by any translate
$y \mapsto h(y \omega)$, where $\omega$ remains within (say) the set $[1/4, 4]$, 
then \eqref{truncated} remains valid.  This is evident from the above proof, replacing
e.g. $Q^{\kappa -1}$ by $\omega Q^{\kappa-1}$. 
  \label{end of proof} \qed

\begin{Prop}\label{udtt}  (variant of \underline{\rm{Claim 2}}). Take $g= a(t) n(T)$ for $t \in \mathbb{R}_{>0} \hookrightarrow \Aa^{\times}, \shift=(T_{v})_{v}\in\Aa$. For $f \in L^2_{0}(\bfX_{\PGL_2})$ and smooth, one has
 $$\mu^g (f)  \ll \Sob^{\bfX_{\PGL_2}}(f)(|t|_\Aa^{1/2}+|\shift|_{\Aa}^{-\delta})$$
 for some absolute $\delta>0$.
\end{Prop}

\proof 
By property (3c) of Sobolev norms, it is sufficient to prove that for any automorphic representation $\pi$ equipped with its canonical norm, and $f\in\pi$ a smooth function,
 we have the inequality 
\begin{equation} \label{repbyrep} |\mu^{g}(f)| \ll \Sob^\pi(f)(|t|_\Aa^{1/2}+|\shift|_{\Aa}^{-\delta}).\end{equation}
We may decompose $\mu^{g}(f) $ as $\tilde{\mu}^{g}(f_{N})+\tilde{\mu}^{g}(f-f_{N})$. 

By \eqref{horcyclebound}
the first term is bounded by $\Sob^\pi(f)|t|_\Aa^{1/2}$  (and is even zero if $\pi$ is cuspidal). The second term equals $\intc_{\Re s=0}\ell^{|\cdot|^s} (g. f) ds,$ by inverse Mellin transform;
 here $\ell^{|.|^s}(f)$  is the linear functional  given in \refs{JLperiod} and associated with the character $|.|_{\Aa}^s$.  Applying Lemma \ref{udtt-cor-refined},  together with \eqref{wtf},  we see that
$(1+ |s|)^4 \ell^{|\cdot|^s} (f) $ is bounded above by $\Sob^{\pi}(f) |T|_{\Aa}^{-\delta} $, with $\delta>0$ absolute, 
and our result follows.  \qed

\subsubsection{The cuspidal case}  
We now prove \underline{Claim 3} and conclude the proof of Theorem \ref{thmcharacter} for $\pi$ cuspidal.  We advise
the reader to first consider \S \ref{outline} which gives a somewhat cleaner version
of the proof in the case $F \neq \Q$. 
	
Let $t, T, \varphi, h$ be as in Lemma \ref{lrep}. 
We need to show that $|\mu_{\chi,h}^{ n(T)}(\varphi)|$ decays with $Q$. 
The basic idea is this: since
 $\mu_{\chi,h}^{n(T)}$ is $\chi$-equivariant under the subgroup of elements of $H^{(1)}$
 which commute with $n(T)$, we can reduce this to the corresponding fact for $\chi=1$ , already known by Proposition \ref{udtt}, using the ergodic principle (\S \ref{Prinzip}).

Noting that $\mu^{a(y_{t}) n(T)}$ is orthogonal to all one-dimensional
 automorphic representations on $\PGL_2$  except the constants, Proposition \ref{udtt} implies that for $f$ non-negative on $\mathbf{X}_{\PGL_2}$, 
\begin{equation} \label{above} |\mu^{n(T)}_{1,h}(f)| \leq {\int h}  \int_{\bfX_{\PGL_2}}  f + \epsilon \Sob^{\bfX_{\PGL_2}}(f), \mathrm{with}\ \epsilon \ll_{\eps} Q^{\eps} (Q^{\frac{\kappa-1}{2}}+Q^{-\delta}),\end{equation}
  for any $\eps>0$. Here
 $\kappa$ is as in the statement of Lemma \ref{lrep}.

Let $\sigma$ be the averaged sum of the Dirac measures which are supported at $$a(\varpi_{v}^{-1})a(\varpi_{v'})\in a(F_{v})a(F_{v'})\subset H$$ for $(v,v')$ ranging over pairs of distinct non-archimedean places at which $\chi$ is unramified and for which $q_{v},q_{v'}$ are
 contained in $[K,2K]$  (to be choose later; in any case $K \leq Q$.) The number of such pairs
is $\gg K^2/\log^2 K$; $\sigma$ is not supported on $H^{(1)}$ but rather on $$H^{([1/4,4])}=\{a(y),\ y\in \Aa^\times,\ |y|\in [1/4,4]\}.$$ 
Since $\| \Ad(a(\varpi_{v}^{-1})a(\varpi_{v'})) \|  \asymp K^2$, we see (notation of \S \ref{Prinzip-quant}):
 $$\|\sigma\|_d \ll K^{2 d}, \|\sigma \star \check{\sigma}\|_{\theta} \ll (\log K)^2\bigl( K^{-2}+K^{-1-2\theta}+K^{-4\theta}\bigr). $$

 The measure $\mu_{\chi, h} (\varphi)$ is not exactly invariant under
the substitution $\varphi\mapsto \varphi\star_{\chi}\sigma$ but almost: indeed since the support of $\sigma$ commute with  $n(T)$ one has
$$\mu_{\chi, h} (\varphi\star_{\chi}\sigma)=\mu_{\chi,h\star\eta}(\varphi),$$
where $\eta$ denote the average of the Dirac measures at $q_{v}/q_{v'}$ on $\Rr_{>0}$ for $(v,v')$ as above.
Therefore, the reasoning of \S \ref{Prinzip-quant} shows that:
\begin{equation} \label{eqA} | \mu_{\chi,  h \star \eta}(\varphi) |^2 \ll_{\varepsilon, \pi}  Q^{\varepsilon} (K^{4d} (Q^{\frac{\kappa-1}{2}}+Q^{-\delta}) + 
 K^{-2}+K^{-1-2\theta}+K^{-4\theta} )\end{equation} 
 Note that $|\varphi \star_{\chi} \sigma|^2$ descends to $\mathbf{X}_{\PGL_2}$,
 and so it was admissible to apply \eqref{above} to it. 

 It follows from the (the last line of the) proof of Lemma \ref{lrep} that \refs{truncated} holds with $h$ replaced by $h\star\eta$ (which is supported on $[\frac{1}{4}Q^{-1-\kappa}, 4Q^{-1+\kappa}]$) :
\begin{equation} \label{eqB} Q^{-1/2-\eps}L(1/2, \pi \otimes \chi) \ll_{\eps,\pi, \kappa}  \left| \mu_{\chi, h\star\eta} (\varphi) \right| 
 + Q^{-\kappa/2}\end{equation} 
 Taking $K$ to be a suitable small power of $Q$ and combining \eqref{eqA} and \eqref{eqB}, we conclude.  \qed

\subsubsection{The Eisenstein case}  \label{gl2eis}

In this section, we prove Theorem \ref{thmcharacter} for $\pi$ noncuspidal by utilizing a simple regularization.  By factorization, 
it suffices to consider the case where $\pi= 1 \boxplus 1$ is induced
from two trivial characters. 

 Let $\varphi$ be as furnished by Lemma \ref{lrep}. 
 Let $k$ be a fixed smooth compactly supported function on $\G(\adele)$ of integral $1$. 
 Fix a sufficiently large parameter $X \geq 1$ (to be a fixed power of $Q$),  and split
 $\varphi = \varphi_1 +\varphi_2$, where 
 $$\varphi_1(g) := \wedge^X \varphi \star k,  \  \varphi_2 = (\varphi - \wedge^X \varphi) \star k, $$
say.   Thus
 $\varphi_1$ is of rapid decay high in the cusp, whereas $\varphi_2$
 is supported high in the cusp.  Let us observe,

\begin{enumerate}
\item $|\varphi_2(x)| \ll \Sob^\pi(\varphi) \height(x)^{1/2} \log \height(x)$. \item $ \left| \int_{\bfX} \varphi_i \right| \ll \Sob^{\pi}(\varphi)  X^{-1/4} \ i=1,2 .$
\item For every $d$, there exists $N(d)$ so that  $\Sob^{\bfX}_d(\varphi_1) \ll X^{N(d)} \Sob^{\pi}(\varphi)$.
\end{enumerate}
The first property follows from \eqref{eis-singular-2}. 
The second property follows from the first, since $\varphi_2$ is supported
in $\height(x) \gg X$ and $\int\varphi_1+\int\varphi_2=0$.    To verify the third property, it suffices to check that 
(for any $\phi \in \pi$ and any $m \geq 1$) we have the bound $$|\wedge^X \phi(x) | \ll_m X^{A(m)}
\height(x)^{-m} \Sob^{\pi}(\phi).$$
By \eqref{fourierexpansion} together with \eqref{eis-singular-2}, 
it suffices to verify that for $y \in \adele^{\times}$, 
$$ \sum_{\alpha \in F^{\times}} W_{\phi}(a(\alpha y)) \ll_m |y|^{-m} \Sob^{\pi}(\phi).$$
This follows from Proposition \ref{zeroboundprop}, together
with the definition of the norm on $\pi$ given by \eqref{canreg}. It is necessary only to observe that
$$\sum_{\alpha \in F^{\times}} \prod_{v} \max(|\alpha y|_v, 1)^{-n} \ll_n |y|^{-n'}.$$ 

Now,  let $g = a(y_t) n(T)$, where $t,T$ are as in Lemma \ref{lrep}. 
 We need to bound 
 $$\tilde{\mu}^g_{\chi}(\varphi - \varphi_N) =  A + B +C +D, \mbox{ where}$$ 
 \begin{itemize}
 \item[-] $A = \mu^g_{\chi}\left( \varphi_1 -  \frac{1}{\vol(\bfX) }\int_{\bfX} \varphi_1 \right)$
 is bounded as in the prior argument by
 $$Q^{\varepsilon} (K^{2d} (Q^{\frac{\kappa-1}{4}}+Q^{-\delta/2}) + 
 K^{-1}+K^{-1/2-\theta}+K^{-2\theta} ).$$ 
 \item[-]  $B =  \mu^g_{\chi}\left(  \frac{1}{\vol(\bfX)} \int \varphi_1  \right)$ is bounded by $\Sob^{\pi}(\varphi) X^{-1/4}$;  
\item[-]  $C =  \tilde{\mu}^g_{\chi} (\varphi_N) $ is bounded, in view of 
\eqref{eis-singular-2}, by $|t|^{1/2-\varepsilon} \Sob^{\pi}(\varphi) =  Q^{\frac{\kappa-1}{2}+\varepsilon} \Sob^{\pi}(\varphi)$. 
 \item[-]  $D = \mu^g_{\chi}(\varphi_2) $.

\end{itemize}
 
Now let us bound $D$.  Let $\varphi^{\circ}$ be the spherical vector of norm $1$
in the representation $\pi$. 
Decompose $\varphi^{\circ} = \varphi^{\circ}_1 + \varphi^{\circ}_2$ just as above.
By \eqref{eis-sing-lowerbound},  $\varphi_2^{\circ} \geq 0$.
By {\em loc. cit.}, it is also true that
-- for $X$ sufficiently large -- $|\varphi_2| \ll \Sob^{\pi}(\varphi) \varphi_2^{\circ}$. 
Therefore, 
\begin{equation} |D| \leq \Sob^{\pi}(\varphi) (|\mu^g(\varphi_1^{\circ}) | + |\mu^g(\varphi^{\circ})|) 
   \ll_{\eps}  \Sob^{\pi}(\varphi) \left( X^{-1/4} + X^{A} (|t|_\Aa^{1/2}+|\shift|_\Aa^{-\delta}  )  \right) 
   \end{equation}
   for some absolute constant $A$; we applied \eqref{repbyrep} 
   to the function $\varphi^{\circ}$, and Proposition \ref{udtt} to the function $\varphi_1^{\circ} - \int_{\bfX_{\PGL_2}} \varphi_1^{\circ}$, observing that $\left| \int_{\bfX_{\PGL_2}} \varphi_1^{\circ}  \right|
   \ll_{\eps} X^{-1/4}$.      We conclude that:

\begin{multline*}Q^{-1/2}L(\chi,1/2)^2=Q^{-1/2}L(\pi \otimes \chi,1/2)\\ \ll_{\eps}  Q^\eps \bigl ( Q^{-\frac{\kappa}{2}}+  X^{-1/4}  + 
K^{-\theta}+X^{A}K^{2d}(Q^{\frac{\kappa-1}4}+ Q^{-\delta/2})
\bigr)
\end{multline*}
and we conclude by an appropriate choice of the parameters.
\qed

\subsection{Subconvex bounds for Rankin/Selberg $L$-functions}
In this section we prove Theorem \ref{RSthm}. Let $\pi_{1},\pi_{2}$ be the two automorphic representations of $\GL_{2}(\Aa)$ considered there and let $\chi_{1},\chi_{2}$ denote their respective central characters. 

\subsubsectionind{} 

Fix a parameter $\defparam\in i\Rr$  whose modulus will be a non-negative power of $C(\pi_{1})$ to be determined later. We shall prove our result on certain slight constraints
on the parameters of $\pi_i$ so as to stay away from singular Eisenstein series. The general case will reduce to this constrained case, as we explain now.%

If $\pi_2$ is not cuspidal, it will be sufficient to take $\pi_2 = 1\boxplus|.|^{\defparam}$.
 Let $\pi_3$ be the Eisenstein series $1\boxplus\chi_{3}$ where $\chi_{3}=(\chi_{1}\chi_{2})^{-1}$;
 we shall suppose that  $\chi_3^2$ is not of the form $|.|^{2it}$ for $|t|<  |\defparam|/3$. 
   Under such constraints we will establish the bound\footnote{ of course the exponent $8$ is not optimal for specific
  configuration (e. g. if $\pi_{2}$ is cuspidal) but this will be sufficient for our present needs}
 $$L(\pi_1 \otimes \pi_2 \otimes \pi_3,\frac{1}{2})  \ll_{\pi_{2}} C(\pi_1 \otimes \pi_2)^{1/2-\delta} |\defparam|^{-8},$$
 for some absolute $\delta>0$. 
 
 This implies our main result: from Theorem \ref{thmcharacter}, we may assume that $\pi_{1}$ is cuspidal. 
  In that case, we wish a subconvex bound for
 $L(\pi_{1}\otimes\pi_{2},1/2)$ with either $\pi_{2}$ cuspidal or $\pi_{2}=1\boxplus 1$ (which yield a subconvex bound for $L(\pi_{1},1/2)$).
We take $\defparam$ such that $|\defparam|=C(\pi_{1})^{-\delta/9}$, and apply the previous bound to a triple
 $(\pi_{1},\pi_{2}',\pi_{3})$ where  $\pi'_{2}$ equals $\pi_{2}$ or $\pi_{2}\otimes|.|^{z}$ if $\pi_{2}$ cuspidal, and equals $1\boxplus|.|^{z}$ if not, while $\pi_{3}=1\boxplus (\chi_{1}\chi'_{2})^{-1}$, choices being made so that $(\pi_{1},\pi_{2}',\pi_{3})$ fulfills the above constraints.
A subconvex bound follows for $L(\pi_{1}\otimes\pi_{2},1/2)$ since, by convexity, we have, for any $t \in [-1,1]$, 
\begin{eqnarray*} \left| L(\pi_{1},1/2)-L(\pi_{1},1/2+it) \right| \ll C(\pi_{1})^{1/4+o(1)}|t|, \\
|L(\pi_{1}\otimes\pi_{2},1/2)-L(\pi_{1}\otimes\pi_{2},1/2+it)|\ll_{\pi_{2}} C(\pi_{1})^{1/2+o(1)}|t|. \end{eqnarray*}

\subsubsection{Choice of the test vectors}\label{choicetest} 
Factorize $\pi_i = \otimes_{v} \pi_{i,v}$, and choose
unitary structures on each $\pi_{i,v}$ so that the product coincides
with the canonical norm on $\pi_i$. Given $\epsilon>0$ small and
 $i=1,2,3$, let $\varphi_{1}\in\pi_{1}$, $\varphi_{2}\in\pi_{2}$, $E \in\pi_{3}$ 
be the tensor product of the test vectors $(\varphi_{1,v})_{v}, (\varphi_{2,v})_{v},\ (E_{3,v})_{v} $ constructed in  Proposition \ref{lowerbounds} for each place $v$ (applied with the parameter $\epsilon$).

The canonical norm of $\varphi_1, \varphi_2,\ E$ are easy to estimate:
they are, by definition, the product of the norms of these local vectors,
and are therefore equal to $1$; moreover we have for any $d$, $$\Sob_d(\varphi_{2}) \ll_{\pi_{2}}1.$$
Applying \eqref{ichino-symmetric}, we have, from the hypotheses made in the above section,
 \begin{eqnarray} \nonumber \frac{L(\pi_{1}\otimes\pi_{2},1/2)}{C(\pi_{1})^{1/2+\epsilon}}&&\ll_{\epsilon, \pi_{2}} |L^*(\pi_2, \Ad,1) L^*(\pi_3, \Ad,1)|\left| \int_{\bfX}\varphi_{1}\varphi_{2}E(g)dg.\right|\\ 
 \label{zump}
 &&\ll_{\epsilon, \pi_{2}}|\defparam|^{-4}\peter{\varphi_{1},\ov{\varphi_{2}E}}. \end{eqnarray}

Set $Q:=C(\pi_{1}) C(\pi_2)$. We want to check that $\peter{\varphi_{1},\ov{\varphi_{2}E}}  \ll_{\pi_2}Q^{-\delta}$. 
Roughly, we use Cauchy-Schwarz to bound the square of this quantity by
$\peter{\varphi_2 E,\varphi_2 E}$. The later equals $\peter{\varphi_2 \ov{\varphi_2}, E\ov E}$ which may be decomposed
 along the automorphic representations of $\PGL_{2}(\Aa)$. It turns out that the contribution of each of these is small {\em except} for one-dimensional representations.
Naively speaking, their contribution would be $\|\varphi_{2}\|^2\|E\|_{can}^2$; the truth
is a little more complicated because we need to use regularization. In any case,
to reduce this ``large'' term, we use Friedlander-Iwaniec's amplification method.

\subsubsection{The amplification method}\label{amplificationsection}We choose a real signed measure $\sigma$ compactly supported on $\GL_{2}(\Aa_{f})$ which 
 satisfies
\begin{enumerate}
\item  \label{ampli0} $u\in\supp(\sigma)\ \Longrightarrow\ \|u\|\leq K$ where $K>0$ is some parameter which will be chosen to be a small fixed positive power of $Q$.   
\item  $\supp(\sigma)$ commutes with $\GL_{2}(F_{v})$ at all archimedean places $v$ and at all places $v$
for which either $\pi_1$ or $\pi_2$ is not spherical or the chosen additive character $\psi$ is ramified.\footnote{i.e. the places dividing the discriminant of $F$.}
\label{ampli1}
\item $(\varphi_2 E ) \star \check\sigma$ is of rapid decay on $\bfX$. 
\item Let $|\sigma|$ denote the  total variation measure. Then
the total mass of $|\sigma|$ is bounded above by $K^{B}$, for some absolute constant $B$. 
Moreover, with $|\sigma|^{(2)}=|\sigma|\star|\check\sigma|$, one has for any 
$\gamma>1/2$
$$\int\|\Ad(u)\|^{-\gamma}d|\sigma|^{(2)}(u)\leq K^{-\eta},\ \hbox{for some }\eta=\eta(\gamma)>0.\label{ampli2}$$\label{greaterthan1/4}
\item $\varphi_{1}\star\sigma=\lambda_1.\varphi_{1}$ with $\lambda_{1} \gg_{F,\eps} Q^{-\eps} $ for any $\eps>0$.\label{ampli3} 

\end{enumerate}

The construction of such a measure (which is inspired by \cite{DFI2,DFI8}) is given in \cite[\S 4.1]{Ve},
except for the remark that $(\varphi_2 E) \star\check\sigma$ is of rapid decay.  However, this follows
by convolving with a measure as in the remark of \S \ref{tripleexample}; 
one may choose the place $v$ so that $q_v \ll (\log Q)^2$. 

By property \refs{ampli3} of $\sigma$, stated above and Cauchy-Schwarz, we have
\begin{align*} |\lambda_1|^{2} \left|\peter{\varphi_{1},\ov{\varphi_{2}E}}\right|^2&=  \left|\peter{\varphi_{1}\star\sigma,\ov{\varphi_{2}E}}\right|^2=
 \left|\peter{\varphi_1, \ov{(\varphi_2.E)\star{\check\sigma}}}\right|^2\\
 & \leq \peter{(\varphi_2.E)\star\check\sigma,(\varphi_2.E)\star\check\sigma};
 \end{align*}
 Thus far, the integrals considered are convergent. However, we shall now expand
 the integral implicit in convolving with $\check\sigma$; at this point, we need
 to make use of regularized integrals. This is possible since, for $u\in\supp(|\sigma|^{(2)})$, the set of exponents of 
 $\varphi_{2}^{u}\ov{\varphi_{2}}\ov{E^{u}}E$ is either $\emptyset$ (if $\pi_{2}$ is cuspidal) or is comprised of characters with real part equal to $2$.
Noting that 
\be\label{RScommute}\peter{\varphi_{2}^{u}.E^{u},\varphi_{2}.E}_{reg}=\peter{\varphi_{2}^{u}\ov{\varphi_{2}},\ov{E^{u}}E}_{reg}\ee
 we obtain that
 \be\label{lastbound} |\lambda_1|^{2} \left|\peter{\varphi_{1},\ov{\varphi_{2}E}}\right|^2\leq \int \left|\peter{\varphi_{2}^{u}\ov{\varphi_{2}},\ov{E^{u}}E}_{reg}\right|d|\sigma|^{(2)}(u)\ee
Now, applying properties  \refs{ampli0} and \refs{ampli2} of $\sigma$ to \refs{lastbound}, 
we conclude that it suffices to prove:

 \begin{Prop}\label{prop:boundfourthintegral} 
For $u \in \GL_2(\adele)$, put
 $\varPhi_{2,u}:=\varphi_{2}^{u} \ov{\varphi_{2}}$ and 
$\varPhi_{3,u}:=\ov{E^{u}} E$.  For $u\in\supp(|\sigma|^{(2)})$ 
$$| \peter{\varPhi_{2,u},\varPhi_{3,u}}_{reg}| \ll_{\pi_2,\epsilon} |z|^{-4}  \|\Ad(u)\|^{-\gamma}+\|u\|^AQ^{-\delta}$$
for some absolute positive constants $\gamma>1/2$ and $\delta,A>0$.
\end{Prop}

The proof will follow. Roughly, we evaluate the right-hand side
via the regularized Plancherel formula in  \S \ref{sec:proofpropfourthintegral};
we handle the spectral sum (cuspidal and Eisenstein) in \S \ref{genericbound},
and we bound the ``degenerate'' terms that arise in \S \ref{sec:degeneratebound}.

\subsubsection{Good and bad places}
Fixing now $u\in\supp(|\sigma|^{(2)})$, we make the following definitions of ``good'' and ``bad'' places:

Let $R$ be the set of finite places where $\varphi_{2,v}, E_{3,v}$ are spherical
and where $u_v = 1$; 
let $S$ be the set of places where $u_v \neq 1$: by the choice made in \S \ref{choicetest} 
$S$ consists of finite places, and the data $\varphi_{2,v},\ E_{3,v},\psi_v$ are all unramified for $v \in S$; let $T$ denote the complement of $R \cup S$.  

Observe that $|S| + |T| = o(\log C(\pi_1) + \log C(\pi_2))$, the bound
on $|S|$ arising from property (1) of the signed measure $\sigma$. 
This property will be used to control a product, over $v \in S \cup T$,
of ``implicit constants.''

\subsubsection{Deformation}
In this section, we shall deform $\pi_2, \pi_3$ in such a way that we will be able to apply the regularized Plancherel formula. These deformations will be parameterized by $s\in \C$ (for $\pi_2$) and $t \in \C$ (for $\pi_3$). We shall fix rather specific deformations of the vectors $\varphi_2 \in \pi_2, E \in \pi_3$. 
 We shall also fix factorizations into local constituents of various
vectors that will arise. 

Write $\pi_3=1\boxplus\chi_{3}$ and 
let $\pi_3(t):=|.|^{t}\boxplus \chi_{3}|.|^{-t}$, for $t \in \C$, be the one-parameter deformation of $\pi_3$ as described in \S 
\ref{notation-global-Eis}.
We choose $f_3$ so that $E =\Eis(f_3)$. 
We denote  the corresponding deformation of $E$ by $E(t)$: $E(t)=\mathrm{Eis}(f_{3}(t))$.
  
If $\pi_2=1\boxplus |.|^\defparam$, we consider the deformation $(|.|^s \boxplus |.|^{\defparam-s},\varphi_{2}(s))$ of $(\pi_2, \varphi_2)$: after choosing $f_2$ so that $\varphi_2 = \Eis(f_2)$, we set $\varphi_{2}(s)= \Eis(f_{2}(s))$.   (By convention, if $\pi_2$ is cuspidal we regard $(\pi_2(s), \varphi(s))$ as being {\em constant}.) Note that $\|E(t)\|_{can}$ and $\|\varphi_2(s)\|_{can}$ are both constant
by Lemma \ref{lem:cannorm}.

Factorizing $f_2 = c_{2,f} \prod_{v}  f_{2,v}$
and $f_3 = c_{3,f} \prod_{v} f_{3,v}$,
where $\|f_{2,v}\| = \|f_{3,v}\| = 1$ for all $v$; then 
\begin{equation} \label{cfbound} |c_{2,f}|\ll 1,\  |c_{3,f}| \ll 1 .\end{equation}

Let $W_{\varphi_2(s)}$ be the image of $\varphi_2(s)$ under the Whittaker intertwiner. 
We may then factorize
$$W_{\varphi_2(s)} =  c_{2,W}(s) \prod_{v} W_{\varphi_2, v}(s),$$ 
where we take $W_{\varphi_2, v}(s)$ to be the image of $f_2(s)$
under the intertwiner \eqref{principaltokirillov} for $v \in S \cup T$; 
and, for $v \in R$, we normalize so that $W_{\varphi_2,v}(s)$
takes the value $1$ at $1$.  In a similar way, we factorize $W_{\varphi_3(t)}= c_{3,W}(t) \prod_{v} W_{\varphi_3, v}(t)$, the normalizations being identical.  
Notice that $\|W_{\varphi_2, v}(s)\|  = \|W_{\varphi_3, v}(t)\| = 1$ for $v \in S \cup T$,
the intertwiner \eqref{principaltokirillov} being isometric. 

The constant $c_{2,W}(s)$ may be evaluated by computing norms. By \S \ref{canGL2}
$$\langle \varphi_2,  \varphi_2\rangle_{can} \asymp |c_{2,W}(s)|^2  \Lambda^*(\pi_2(s), \Ad, 1) \prod_{v \in S \cup T} \frac{1}{\zeta_v(1) L(\pi_{2,v}(s), \Ad, 1)/\zeta_v(2) }$$
and therefore, for any $\eps > 0$ ,
\begin{equation} \label{cwbound}  |c_{2,W}(s)|^2\ll_{\eps}  \frac{ C(\pi_1)^{\eps} C(\pi_2)^{\eps} }{L^*(\pi_2(s), \Ad, 1)},\end{equation}
 and similarly for $c_{3,W}(t)$. 

 \subsubsection{Application of the regularized Plancherel formula} \label{sec:proofpropfourthintegral} We note that $\varPhi_{2,u}$ and $\varPhi_{3,u}$ descend to function on $\PGL_{2}(\Aa)$.  
 Throughout this section we regard $u$ as fixed; 
set (we have suppressed the dependence on $u$ for typographical simplicity)
$$\varPhi_{2}(s) = \varphi_2^u \overline{\varphi_2(s)}, \varPhi_3(t) = \overline{E^u} E(t), 
\scE_{2}(s) = \Eis(\varphi^u_{2N} \overline{\varphi_{2N}(s)}) , \scE_3(t) = \Eis(\overline{E_N^u}E_N(t)).$$

\begin{itemize}
\item[-]
The set of exponents of $\varPhi_{2}(s)$ is
$S_{2}(s)=\{|.|^{1+ \ov s}, |.|^{1- \ov s},|.|^{1+ \defparam + \ov s  }, |.|^{1+ \ov \defparam - \ov s}\}$
(or the empty set if $\pi_{2}$ is cuspidal).
\item[-]The set of exponents of $\varPhi_{3}(t)$ is
$S_{3}(t)=\{|.|^{1+t},|.|^{1-t},\chi_{3}|.|^{1-t},\ov{\chi_{3}}|.|^{1+t}\}.$
\end{itemize}

 In particular, $(s,t)\mapsto\peter{\varPhi_{2}(s),\varPhi_{3}(t)}_{reg}$ defines an $(s,t)$-anti-holomorphic function in a neighbourhood of $(0,0)$ in $\mathbb{C}^2$. 
We will bound  the function $\peter{\varPhi_{2}(s),\varPhi_{3}(t)}_{reg}$ at the point $(0,0)$ by analyzing 
it along a suitable non-empty subset $\paramset\subset\mathcal{D}$ containing $(0,0)$ in its closure. 
The key point is that we should choose $\paramset$ so that \eqref{RIP} is applicable, 
and also so that $\pi_2(t), \pi_3(s)$ do not equal $1 \boxplus 1$ for $(t,s) \in \paramset$. For explicitness, take
\begin{equation} \label{paramsetdef} \paramset = \{(s,t)=(t/2,t),\ t\in i\R: 0<|t| < |z|/6\}.\end{equation}

By choice of $\mathcal{N}$, we can apply the regularized Plancherel formula \eqref{RIP}. 
Moreover, for $\pi$ belonging to the finite spectrum ($\pi=\chi$, $\chi^2=1$) and $(s,t)\in\paramset$ one has, again by invariance, $\Pi_{\pi}(\varPhi_{3}(t))=0$. (In other words,
there is no $\G(\adele)$-equivariant functional from the tensor product $\overline{\pi_3} \otimes \pi_3(t)$
to a one-dimensional $\G(\adele)$-representation). 
Therefore, for $(s,t)\in\paramset$, we have
 \begin{eqnarray}\label{regplancherel2}\peter{\varPhi_2(s),\varPhi_3(t)}_{reg}=&&\peter{\varPhi_2(s),\scE_3(t)}_{reg}+
\peter{\scE_{2}(s), \varPhi_3(t)}_{reg}\\ 
&&\ \ + 
\int_{\pi\ \mathrm{generic}}\peter{\Pi_{\pi}(\varPhi_2(t)),\Pi_{\pi}(\varPhi_3(s))}d\muP(\pi).\nonumber
\end{eqnarray}
where $\Pi_{\pi}$ is defined in \refs{projreg}.

We have already observed that $(s,t)\mapsto \peter{\varPhi_2(s),\varPhi_3(t)}_{reg}$ defines an $(s,t)$-anti\-holo\-morphic function in a neighbourhood of $(0,0)$, and, in particular,
a continuous function on $\mathcal{N}$;  since the same is true of the function $$(s,t)\mapsto\int_{\pi\ \mathrm{generic}}\peter{\Pi_{\pi}(\varPhi_2(t)),\Pi_{\pi}(\varPhi_3(s))}d\muP(\pi),$$ it follows from \refs{regplancherel2} that
\begin{Lemma*} The degenerate term $$(s,t)\in\mcN\mapsto \peter{\varPhi_2(s),\scE_3(t)}_{reg}+
\peter{\scE_{2}(s), \varPhi_3(t)}_{reg}$$ extends to a continuous function on $\mathcal{N}$. 
\end{Lemma*}

\subsubsection{Bounding the generic term}\label{genericbound}
\begin{Lemma*}
For $(s,t)\in\paramset$, the last term of \eqref{regplancherel2} is  $\ll_{\pi_{2}}\|u\|^AQ^{-\delta}$,
for absolute $\delta>0$ and $A$.
 \end{Lemma*}

 In fact, we need this bound only for $(s,t) = (0,0)$, and the reader is welcome
 to make this substitution in what follows.
 \proof Let $d_{0}$ be the  larger than any of the (absolute) Sobolev degrees occuring in Lemmata \ref{tempbound} and \ref{upperboundexplicit} as well as the degree of \eqref{trivbond}. 
 For any generic standard automorphic representation of $\PGL_2$, 
 the inner product  
 $\big|\peter{\Pi_{\pi}(\varPhi_2(s)),\Pi_{\pi}(\varPhi_3(t))}\big|$
 is bounded above
 by $\Sob^{\pi}_{d_0}(\varPhi_2(s)) \Sob^{\pi}_{-d_0}(\varPhi_3(t))$ (as a short cut we write $\Sob^{\pi}_{d}(\varPhi)$ for 
 $\Sob^{\pi}_{d}(\Pi_{\pi}(\varPhi))$), 
 by the duality between the Sobolev norms in question . It will be sufficient to verify that 
  $$\Sob^{\pi}_{-d_0}(\varPhi_3(t)) \ll_{\pi_{2},\pi,\epsilon} \|u\|^B Q^{-\delta},$$ 
  for absolute constants $\delta>0$ and $B$; indeed, for any constant $B>0$, one has
   $$\int_{\pi} C(\pi)^{B}\Sob^{\pi}_{d_0}(\varPhi_2(s))\ll_{B} \int_{\pi} \Sob^{\pi}_{d_0+O(B)}(\varPhi_2(s)) \ll_{\pi_{2},B} 1,$$ 
   the latter inequality following from \refs{projectioncoarse}.

We shall apply the results of \S \ref{ss:ub} especially \refs{ichino-upper-bound}. 
In the notations of that section, we utilize 
the corollary  to Lemma \ref{upperboundexplicit}
 \footnote{if $\pi_{v}$ is tempered, e.g., if $\pi$ is Eisenstein or under the Ramanujan-Peterson conjecture
 we may, more simply, apply Lemma \ref{tempbound}.} to see -- since $C(\chi_{3})\ll_{\pi_{2}}C(\pi_{1})$ and $t\in i\Rr$ -- 
$$ \label{boundA} B_v \ll_{\epsilon,\pi_2}
C(\pi_{1,v})^{d' \eps}C(\chi_{3,v})^{-1}\bigl(\frac{C(\pi_{1,v})}{C(\chi_{3,v})}\bigr)^{2\theta-1} \ \ (v \in T). 
$$ for some absolute $\delta>0$, $d' > 0$.  Note that $u_v$ does not contribute, since $u_v=1$ for $v \in T$.  On the other hand, we have 
\begin{equation}  \label{boundB} B_v \ll \|u_v\|^A \hbox{ for $v \in S$} \end{equation}
and some unspecified constant $A$: this follows from the ``trivial'' bound \eqref{trivbond} (see the sentence following that equation).

We take into account also
 the global subconvex bound of Theorem \ref{thmcharacter};
it implies $L(\overline{\pi_3} \otimes \pi_3 \otimes \pi,1/2 ) \ll_{\pi} C(\chi_{3})^{1-\delta}$. Thus, by
the results of \S \ref{ss:ub}, 
$$\Sob^{\pi}_{-d_0}(\Pi_{\pi}(\varPhi_3(t)) \ll_{\pi_{2},\pi,\epsilon} C(\chi_3)^{-\delta} \left(\frac{ C(\pi_1)}{C(\chi_3)} \right)^{2 \theta -1} \|u\|^A
\|E\|_{can}^2 \|E(t)\|_{can}^2.$$
Now, $\|E(t)\|_{can} = \|E\|_{can} =1$. 
Using again that $C(\chi_{3}) \ll_{\pi_{2}} C(\pi_1)$, we conclude the proof of the Lemma.
 \qed

\subsubsection{Bounding the degenerate term} \label{sec:degeneratebound}
\begin{Lemma*} For $(s,t)\in\paramset$ we have, for any $\eps>0$
\be\label{degeneratebound}
\left| \peter{\varPhi_2(s),\scE_3(t)}_{reg}+
\peter{\scE_{2}(s), \varPhi_3(t)}_{reg} \right| 
\ll_{\pi_{2},\epsilon,\eps} Q^{\eps} |z|^{-4} \|\Ad(u)\|^{-\gamma}\ee
for some absolute $\gamma>1/2$.
\end{Lemma*}

\proof We start with the first portion of degenerate term: $\langle \varPhi_2(s), \scE_3(t) \rangle_{reg}$.

Let us consider the linear functional on $$ \Pi=\Pi(s,t) := \pi_2 \otimes \overline{\pi}_2(\ov s) \otimes  \pi_3 \otimes \overline{\pi}_3(\ov t) \rightarrow \C$$ defined by the rule $$H=H_{s,t}: (\varphi_2, \varphi_2', \varphi_3, \varphi_3') \mapsto \int^{reg}_{\mathbf{X}_{\PGL_2}} \varphi_2 \varphi_2' \Eis(\varphi_{3,N} \varphi'_{3,N}).$$ 
Thus, $\langle \varPhi_2(s), \scE_3(t) \rangle = H(\varphi_2^u, \ov{\varphi_2(s)}, E^u, \ov{E(t)})$;
also, expanding the constant terms $\varphi_{3,N},\ \varphi'_{3,N}$, we may express $H_{s,t}$ as a sum of four terms 
$$H_{s,t} = \sum_{\pm,\pm} H_{s,t}^{\pm\pm},$$ where,
if we realize $\pi_2$ in its Whittaker model and $\pi_3$ in the principal series model, 
\be \label{Hstdef} H_{s,t} ^{++ }(\varphi^{u}_2, \ov{\varphi_2(s)}, \Eis(f_3)^{u}, \ov{\Eis(f_3(t))})  =  c \prod_{v}^*H^{++}_{v},\ee
$$ c =  c_{2,W}\ov{c_{2,W}(s)}|c_{3,f}|^2 , \ \ H^{++}_v :=  \int_{N(F_v) \backslash \PGL_2(F_v)}
W^{u_{v}}_{2, v} \overline{W_{2, v}(s)}  f^{u_{v}}_{3,v} \ov{f_{3,v}(t)}.$$

The other terms (e.g. $H_{s,t}^{+-}$,  $H_{v}^{+-}$) are defined similarly, by introducing standard intertwining operators $M$ (see \S \ref{intertwining}) 
in front of the $f_3$'s. The second portion of the degenerate $\langle \scE_2(s), \varPhi_3(t) \rangle_{reg}$, admits analogous expansion, 
of the shape $\sum_{\pm,\pm}J_{s,t}^{\pm\pm}$ and the evaluation of each $J^{\pm\pm}_{s,t}$ is entirely similar to that of $H^{\pm\pm}_{s,t}$. 
 We will examine in detail two terms: $H^{++}$ and $J^{--}$. We then complete the proof of the Lemma in \S \ref{alltogether}.

\subsubsection{Local bounds for $H$}\label{sec:localboundH}
We shall need bounds for $H_v^{++}$ as well as for its partial derivatives
w.r.t. $s, t$ on the unitary axes, ie. when $s$, $t$ are purely imaginary. In fact, bounds for the partial derivatives will be derived from bounds for $H_v^{++}$ in a small neighborhood of the unitary axes; however, for clarity, we begin by explaining the bounds for $H_v^{++}$ when $s,t$ are on the unitary axes. 

The evaluation of $H_v^{++}$, for $v \in R$, is the theory of local Rankin-Selberg integrals for $\GL_2$ 
(cf. \eqref{rs-ae}):
\begin{equation} \label{++eval} H_v^{++} = %
 \frac{ L_v(\pi_2 \otimes \tilde\pi_2(\ov s)  \otimes |\cdot|^{\ov t},1)}{\zeta_{v}(2+2\ov t)}\end{equation}

For $v \in S$, the vectors are spherical and  the additive character $\psi_{v}$ is unramified;
we apply Lemma \ref{jisometry} to conclude:
\begin{equation}\label{++Sbound}
|H_v^{++}| \ll \|\Ad (u_{v})\|^{-\gamma},\ \gamma=1-\theta>1/2. 
\end{equation}
For $v$ in $T$ we have again by Lemma \ref{jisometry}
\begin{eqnarray}
|H^{++}_{v}|&&\leq \|f_{3,v}W_{2,v}\|_{N_v\bash {G}_{v}}
\|f_{3,v}(t)W_{2,v}(s))\|_{N_v\bash {G}_{v}}\label{++Tbound}\\
&&=\|W_{2,v}\|\|W_{2,v}(s)\|\|f_{3,v}\|\|f_{3,v}(t)\|=\|f_{3,v}\|^2\|W_{2,v}\|\|W_{2,v}(s)\|
= 1.
\nonumber
\end{eqnarray}

More generally, for any fixed $i,j\geq 0$, $v\in S\cup T$ and any $\eps>0$, we have also 
\be\label{boundderivarives} |\partial_s^i \partial_t^j H_v^{++}| \ll_{i,j,\eps} C(\pi_{2,v})^\eps\|\Ad(u_{v})\|^{\eps-\gamma};
\ee

It is to verify \eqref{boundderivarives} that we consider $H_v^{++}$ off the unitary axes. 
The function $H_v^{++}$ being antiholomorphic in $s,t$, it suffices
by Cauchy's formula to bound it in a small neighbourhood of $(s,t) = (0,0)$. 
The required bound follows, for $v \in S$, from Lemma \ref{degboundlem}; 
for $v \in T$, it follows by a reasoning similar to the previous one that
for $|\Re(s)| + |\Re(t)| \leq \eps/2$, 
\begin{align*} |H_v^{++}(s,t)|^2 \leq
 \int | W_{2,v}(s) (a(y)) |^2  \max(|y|, |y|^{-1})^{\eps} d^{\times}y & \ll_{\eps}   \Sob_{d \eps}^{\pi_{2,v}} (W_{2,v})^2 \\ 
  & \ll_{\eps}   C(\pi_{2,v})^{d'\eps}, \end{align*}
where we applied Proposition \ref{zeroboundprop-def}, and $d,d'$ are absolute constants;
we also used the bound for $\Sob(W_{2,v})$ given in Proposition \ref{lowerbounds}. 

The same bounds apply to all the terms $H_v^{\pm, \pm}$.  We give
an example of how to handle the intertwining operators that intervene, in relation to the $J$-term. 

\subsubsection{The $J$ terms} \label{Jsubsec}
For the corresponding terms
$J^{\pm\pm}_{v}$ a bound similar to \refs{boundderivarives} applies. We shall study the term $J^{--}_v$ to make this reasoning clear; it will also give us the opportunity to explain the reasoning
involving intertwining operators. 

The second degenerate term $\peter{\scE_{2}(s), \varPhi_3(t)}_{reg}$ is zero unless $\pi_2$ is Eisenstein, so we shall suppose
that $\pi_2$ is Eisenstein.  Corresponding to \eqref{Hstdef} we have the identity:

$$ \label{Jstdef} J_{s,t} ^{-- }(\varphi^{u}_2, \ov{\varphi_2(s)}, \Eis(f_3)^{u}, \ov{\Eis(f_3(t))})  =  c' \prod_{v}^*J^{--}_{v},$$
 $$ |c'| =  |c_{2}\ov{c_{2}}c_{3,W}\ov{c_{3,W}(t)}| , \ \ J^{--}_v :=  \int_{N(F_v) \backslash \PGL_2(F_v)}
W^{u_{v}}_{3, v} \overline{W_{3, v}(t)}  \bar{M}_v f^{u_{v}}_{2,v} \ov{\bar{M}_v f_{2,v}(s)};$$
let us remind the reader that $f_{2}=\otimes_{v} f_{2,v}$ is so that $\varphi_{2}=\Eis(f_{2}),\ \varphi_{2}(s)=\Eis(f_{2}(s))$ and that $\bar{M}_v$ is as defined in \S \ref{intertwining}.  Recall also that $\bar{M}_v$ is isometric for every $v$ 
(at least up to a factor depending only on $\psi_v$). 

The unramified evaluation is similar to the prior one.

For $v \in S$, the vectors are spherical and  the additive character $\psi_{v}$ is unramified;
we apply Lemma \ref{degboundlem} to conclude:
\begin{equation}\label{++SboundJ}
|J_v^{--}| \ll \|\Ad (u_{v})\|^{-\gamma} ,\ \gamma=1-\theta>1/2. 
\end{equation}

For $v$ in $T$ we have again by Lemma \ref{jisometry} 
$$|J^{--}_{v}|\leq 1, \label{++Jbound}
$$

We need, again, a bound on derivatives. For fixed $i,j\geq 0$, $v\in S\cup T$ and any $\eps>0$, we have also 
\be\label{boundderivarivesJ} |\partial_s^i \partial_t^j J_v^{--}| \ll_{i,j,\eps} \left(C(\pi_{1,v}) C(\pi_{2,v})\right)^\eps\|\Ad(u_{v})\|^{\eps-\gamma}.\ee  
For this, just as as before, we bound $J_v^{--}$ for $(s,t)$ in a small neighbourhood of $(0,0)$; assuming that
$|\Re (s)|+|\Re(t)|\leq \eps/2$, we find, for $v\in T$, using Cauchy-Schwarz and Lemma \ref{jisometry}  that
$$|J_{v}^{--}|^2\leq \int_{F_{v}^\times\times K}|W_{3,v}(a(y)k)|^2\max(|y|,|y|^{-1})^\eps|\bar M_{v}f_{2,v,s}(k)|^2\dti y dk.$$
The bound \refs{boundderivarivesJ} follows the bounds for $\bar{M}_v$ furnished in \S \ref{intertwining} and the same reasonning as at the end of \S \ref{sec:localboundH}. For $v\in S$ the bound follows from Lemma \ref{degboundlem}.

\subsubsection{Putting it all together}\label{alltogether}

The degenerate term $$\peter{\varPhi_2(s),\scE_3(t)}_{reg}+
\peter{\scE_{2}(s), \varPhi_3(t)}_{reg} \ \ (s,t)\in\mcN$$ is given by
$\sum_{\pm,\pm} H^{\pm, \pm}_{s,t} + J^{\pm\pm}_{s,t}$ and, although the individual terms $H^{\pm\pm}_{s,t}, J^{\pm\pm}_{s,t}$ may not be regular in a neighbourhood of $(s,t) = (0,0)$, we have seen in the lemma preceding \S \ref{genericbound} that their sum is. In particular,  $$t\mapsto \peter{\varPhi_2(t/2),\scE_3(t)}_{reg}+
\peter{\scE_{2}(t/2), \varPhi_3(t)}_{reg}$$ defines an antiholomorphic function in the complex disc $|t|<0.1$.

Each quantity $H^{\pm\pm}_{t/2,t}, J^{\pm\pm}_{t/2,t}$ is the product of: a constant
$c$ as in \eqref{Hstdef}; a partial $L$-function at places outside $S \cup T$,
given e.g. by \eqref{++eval};  and local factors at $S \cup T$.
\begin{itemize}
\item[-] The constant $c$ satisfies $c \ll_{\pi_2, \eps} C(\pi_1)^{\eps}$ by \refs{cfbound} and \refs{cwbound}. 
\item[-]  The local factors at $S \cup T$, for $t\in i\Rr$, $|t|<0.1$, are bounded, along with their derivatives by  \eqref{++Sbound}, \eqref{++Tbound} and \refs{boundderivarives} and the corresponding bounds for $J$.  In particular,
the product of such factors is bounded by $Q^{\eps} \|\Ad(u)\|^{\eps-\gamma}$. 
\item[-]The partial $L$-functions in question extend to meromorphic functions of the $t$-variable with a pole at $t=0$ of order at most $4$. For $t\in i\Rr$ bounded away from $0$ these $L$-functions are bounded by $$\ll_{\eps} |z|^{-4} C(\pi_2)^{\eps} C(\pi_3)^{\eps}$$ for any $\eps>0$, and  at $t = 0$ the terms of bounded order in their Laurent expansion satisfies the same bound. 
\end{itemize}

It follows that for $(s,t)\in \mcN$, the degenerate contribution is bounded by
$$\ll_{\eps} Q^{\eps} |z|^{-4} \|\Ad(u)\|^{-\gamma},\ \gamma>1/2,$$ as required. 
\qed

This lemma together with the lemma of \S \ref{genericbound} conclude the proof of Proposition \ref{prop:boundfourthintegral},
hence of Theorem \ref{RSthm}.



 \end{document}